\documentclass[11pt,reqno]{amsart}
\title
[Meromorphic tensor equivalence]
{Meromorphic tensor equivalence for Yangians and quantum loop algebras}

\author[S. Gautam]{Sachin Gautam}
\address{Perimeter Institute for Theoretical Physics,
31 Caroline Street N., Waterloo, ON N2L 2Y5 (Canada)}
\email{sgautam@perimeterinstitute.ca}
\author[V. Toledano Laredo]{Valerio Toledano Laredo}
\address{Department of Mathematics,
Northeastern University,
360 Huntington Avenue,
Boston, MA 02115}
\email{V.ToledanoLaredo@neu.edu}

\usepackage{amsmath,amssymb}
\usepackage{mathrsfs}
\usepackage[usenames,dvipsnames]{color}
\usepackage{hyperref}
\hypersetup{
    colorlinks=true,       			
    linkcolor=blue, 
    citecolor=blue,        			
    filecolor=blue,      				
    urlcolor=blue           			
}

\usepackage[all]{xy}
\xyoption{arc}
\newtheorem*{thm}{Theorem}
\newtheorem*{prop}{Proposition}
\newtheorem*{lem}{Lemma}
\newtheorem*{cor}{Corollary}

\newenvironment{pf}{\paragraph{{\sc Proof}}}{\qed\par\medskip}
\theoremstyle{definition}
\newtheorem*{defn}{Definition}
\newtheorem*{rem}{Remark}

\numberwithin{equation}{section}

\numberwithin{figure}{section}

\newcommand{\Aq}{\mathscr{A}}
\newcommand{\Rq}{\mathscr{R}}
\newcommand{\Rqo}{\mathscr{R}^{0,-}}
\newcommand{\Rqi}{\mathscr{R}^{0,+}}
\newcommand{\Rqio}{\mathscr{R}^{0,\pm}}
\newcommand{\Rqoi}{\mathscr{R}^{0,\mp}}

\newcommand{\calD}{{\mathcal D}}
\newcommand{\calR}{{\mathcal R}}
\newcommand{\Ru}{\mathcal{R}^{0,+}}
\newcommand{\Rd}{\mathcal{R}^{0,-}}
\newcommand{\Rud}{\mathcal{R}^{0,\pm}}
\newcommand{\Rdu}{\mathcal{R}^{0,\mp}}
\newcommand{\Reps}{\mathcal{R}^{0,\veps}}

\newcommand{\bfB}{{\mathbf B}}
\newcommand{\bfD}{{\mathbf D}}
\newcommand{\bfC}{{\mathbf C}}

\newcommand{\spec}{\sigma}

\newcommand{\Fh}[1]{\mathsf{\Gamma}_{#1}}

\newcommand{\id}{\mathbf{1}}

\newcommand{\ve}{\varepsilon}

\newcommand{\lp}{\left(}
\newcommand{\rp}{\right)}
\newcommand{\la}{\left\langle}
\newcommand{\ra}{\right\rangle}

\newcommand{\g}{\mathfrak{g}}
\newcommand{\h}{\mathfrak{h}}

\newcommand{\Sym}{\mathfrak{S}}

\newcommand{\bfA}{\mathbf{A}}

\newcommand{\A}{\mathcal{A}}

\newcommand{\D}{\mathcal{D}}

\newcommand{\J}{\mathcal{J}}

\newcommand{\calO}{\mathcal{O}}

\newcommand{\V}{\mathcal{V}}
\newcommand{\W}{\mathcal{W}}
\newcommand{\X}{\mathcal X}

\newcommand{\C}{\mathbb{C}}
\newcommand{\nC}{\mathbb{C}^{\times}}
\newcommand{\N}{\mathbb{Z}_{\geq 0}}
\newcommand{\Q}{\mathbb{Q}}
\newcommand{\R}{\mathbb{R}}
\newcommand{\Z}{\mathbb{Z}}

\newcommand {\Exp}[1]{e^{2\pi\iota #1}}
\newcommand {\wh}[1]{\widehat{#1}}
\newcommand {\ol}[1]{\overline{#1}}
\newcommand {\ul}[1]{\underline{#1}}

\newcommand{\End}{\operatorname{End}}
\newcommand{\Hom}{\operatorname{Hom}}
\newcommand{\ad}{\operatorname{ad}}
\newcommand{\Ad}{\operatorname{Ad}}

\newcommand{\qKZ}{q\operatorname{KZ}}
\newcommand {\aand}{\qquad\text{and}\qquad}

\renewcommand {\theenumi}{\arabic{enumi}}

\newcommand {\Rep}{\operatorname{Rep}}
\newcommand {\ffd}{\operatorname{fd}}
\newcommand {\fd}{finite--dimensional }
\newcommand {\lhs}{left--hand side }
\newcommand {\rhs}{right--hand side }
\newcommand {\wrt}{with respect to }

\newcommand{\Rloop}{\Rep_{\ffd}(\qloop)}
\newcommand{\Rlsub}{\Rep_{\ffd}^{\Omega}(\qloop)}

\newcommand{\Ryang}{\Rep_{\ffd}(\Yhg)}
\newcommand{\Rync}{\Rep_{\ffd}^{\scriptscriptstyle{\operatorname{NC}}}(\Yhg)}
\newcommand{\Rysub}{\Rep_{\ffd}^{\Pi}(\Yhg)}

\newcommand{\ds}{\displaystyle}
\newcommand{\wt}[1]{\widetilde{#1}}

\newcommand {\Omit}[1]{}

\newcommand{\pp}{2\pi\iota}

\newcommand {\qloop}{U_q(L\g)}

\newcommand {\Yhg}{Y_\hbar(\g)}

\newcommand {\isom}{\stackrel{\sim}{\rightarrow}}

\newcommand{\qbin}[3]{\left[\begin{array}{c} #1 \\ #2\end{array}\right]_{#3}}

\renewcommand {\sl}{\mathfrak{sl}}

\newcommand {\sfA}{{\mathsf A}}

\newcommand {\sfC}{{\mathsf C}}
\newcommand {\sfD}{{\mathsf D}}
\newcommand {\sfE}{{\mathsf E}}

\newcommand {\KM}{Kac--Moody }

\newcommand {\eg}{{\it e.g., }}

\newcommand {\ie}{{\it i.e.}, }

\newcommand {\KD}{Kohno--Drinfeld }

\newcommand{\bfI}{{\mathbf I}}

\newcommand{\DD}[1]{\Delta_{#1}}
\usepackage{graphicx}
\usepackage{esint}

\renewcommand {\theenumi}{\roman{enumi}}

\renewcommand {\Re}{\operatorname{Re}}
\renewcommand {\Im}{\operatorname{Im}}

\newcommand {\IP}{\mathbb P}
\newcommand {\Y}{\mathcal Y}
\newcommand {\QL}{$\mathcal{QL}$}
\newcommand {\veps}{\varepsilon}

\newcommand {\calC}{\mathcal{C}}

\newcommand {\Uqg}{U_q\g}

\newcommand {\nc}{non--congruent }

\newcommand {\Dotimes}[1]{\otimes_{#1}}
\newcommand{\opA}[2]{\mathcal{A}_{#1,#2}}

\newcommand {\T}{T}

\newcommand {\sfX}{{\mathsf X}}
\newcommand {\sfY}{{\mathsf Y}}
\newcommand {\DYhg}{\calD\Yhg}
\newcommand {\<}{\langle}
\renewcommand {\>}{\rangle}

\newcommand {\mmu}{} 
\newcommand {\mma}{} 

\newcommand {\ud}{^\pm}

\newcommand{\qnum}[2]{[#1]_{q^{#2}}}
\newcommand{\whacky}[1]{\{#1\}_q}
\newcommand{\cowt}[1]{\lambda^{\vee}_{#1}(q)}

\newcommand {\calB}{{\mathcal B}}

\newcommand {\KT}{Khoroshkin--Tolstoy }
\newcommand {\DY}{{\mathcal D}\Yhg}

\newcommand {\Ztimes}{\Z_{\neq 0}}
\newcommand {\Ntimes}{\Z_{>0}}

\newcommand {\Ju}{\mathcal J^+}
\newcommand {\Jd}{\mathcal J^-}

\newcommand {\Jeps}{\mathcal J^\veps}
\newcommand {\Jmeps}{\mathcal J^{-\veps}}

\newcommand {\rat}{^{\scriptscriptstyle{\operatorname{rat}}}}
\newcommand {\mer}{^{\scriptscriptstyle{\operatorname{mer}}}}

\newcommand {\logo}{\log_0}
\newcommand {\calS}{{\mathcal S}}
\newcommand {\uls}{\ul{s}}

\newcommand {\FR}{Frenkel--Reshetikhin }
\newcommand {\TV}{Tarasov--Varchenko }
\newcommand {\KL}{Kazhdan--Lusztig }

\newcommand{\Rfloop}{\Rep_{\ffd}(\qloop)}
\newcommand{\Rfyang}{\Rep_{\ffd}(\Yhg)}
\newcommand{\IQ}{\mathbb{Q}}

\newcommand {\n}{\mathfrak n}
\newcommand {\sfn}{n}

\begin{document}

\begin{abstract}
Let $\g$ be a complex semisimple Lie algebra, and $\Yhg$,
$\qloop$ the corresponding Yangian and quantum loop algebra,
with deformation parameters related by $q=e^{\pi\iota\hbar}$.
When $\hbar$ is not a rational number, we constructed in \cite
{sachin-valerio-2} a faithful functor $\Fh{}$ from the category
of \fd representations of $\Yhg$ to those of $\qloop$.
The functor $\Fh{}$ is governed by the additive difference
equations defined by the commuting fields of the Yangian,
and restricts to an equivalence on a subcategory of $\Ryang$
defined by choosing a branch of the logarithm.
In this paper, we construct a tensor structure on $\Fh{}$
and show that, if $|q|\neq 1$, it yields an equivalence of
meromorphic braided tensor categories, when $\Yhg$ and
$\qloop$ are endowed with the deformed Drinfeld coproducts 
and the commutative part of their universal $R$--matrices.
This proves in particular the \KD theorem for the abelian
$q$KZ equations defined by $\Yhg$.
The tensor structure arises from the abelian $q$KZ equations
defined by an appropriate regularisation of the commutative
part of the $R$--matrix of $\Yhg$.
\end{abstract}

\thanks
{The second author was supported in part through the NSF grant DMS--1206305}
\maketitle

\setcounter{tocdepth}{1}
\tableofcontents

\newpage
\section{Introduction}\label{sec: intro}

\subsection{} 

Let $\g$ be a complex semisimple Lie algebra, and $\Yhg$
and $\qloop$ the Yangian and quantum loop algebra of $\g$.
Recall that the latter are deformations of the enveloping algebras
of the current Lie algebra $\g[s]$ of $\g$ and its loop algebra
$\g[z,z^{-1}]$ respectively. We shall assume throughout that
the deformation parameters $\hbar$ and $q$ are related by
$q=e^{\pi\iota\hbar}$, and that $q$ is not a root of unity.

The present paper builds upon the equivalence of categories
of \fd representations $\Fh{}:\Rfyang\longrightarrow\Rfloop$
constructed in \cite{sachin-valerio-2}.\footnote{Strictly speaking,
$\Fh{}$ is defined on a subcategory of $\Rfyang$, and becomes
an equivalence after restricting the source category suitably.
We will gloss over this point here, and refer the reader to \cite
{sachin-valerio-2} or Section \ref{sec: functor} below for the
precise statement.} A natural question stemming from \cite
{sachin-valerio-2} is whether $\Fh{}$ is a tensor functor, that
is admits a family of natural isomorphisms $\J_{V_1,V_2}:
\Fh{}(V_1)\otimes\Fh{}(V_2)\to\Fh{}(V_1\otimes V_2)$ of
$\qloop$--modules which
are associative \wrt triples of representations.\footnote
{Although $\Fh{}(V)=V$ as vector spaces, $\J_{V_1,V_2}=
\operatorname{id}_{V_1\otimes V_2}$ is not the required isomorphism
since the actions of $\qloop$ on $\Fh{}(V_1)\otimes\Fh{}(V_2)$
and $\Fh{}(V_1\otimes V_2)$ do not coincide.} Partial evidence,
pointing towards a positive answer, is obtained in \cite{sachin-valerio-2}
where it is shown that $\Fh{}$ is compatible with taking the
$q$--characters of \FR and Knight, and therefore induces
a homomorphism of Grothendieck rings.

\subsection{} 

The goal of this paper is to show that $\Fh{}$ admits a
tensor structure. Our main result, which will be explained
in more detail below, is that $\Fh{}$ gives rise to an
equivalence of {\it meromorphic} tensor categories.
Moreover, when $|q|\neq 1$, this equivalence also
preserves
the meromorphic braiding arising from the abelianisation
of the universal $R$--matrices of $\Yhg$ and $\qloop$. 

This may be regarded as a meromorphic version of the \KL
equivalence between the category $\calO$ of representations of
the affine algebra $\wh{\g}$ at level $\kappa$ and the category of \fd
representations of the quantum group $\Uqg$, where
$q=e^{\pi\iota/m(\kappa+h^\vee)}$ \cite{KL12,KL3,KL4}.
Here $m$ is the ratio of the square length of the long roots
to the short roots and $h^{\vee}$ is the dual Coxeter number.
More precisely, for $\kappa\notin\IQ$, the central ingredient
of the KL equivalence is the construction of a tensor functor
from the (Drinfeld) category $\calD(\g)$ of \fd $\g$--modules,
with associativity and commutativity constraints given by
the monodromy of the KZ equations with deformation
parameter $\hbar=1/(\kappa+h^\vee)$, to the category of \fd
representations of $\Uqg$ \cite{KL3}.

In the present work, $\calD(\g)$ is replaced by $\Rfyang$,
$\Rep_{\ffd}(\Uqg)$ by $\Rfloop$, and the KZ equations by
an appropriate abelianisation of the additive, difference
$q$KZ equations defined by the universal $R$--matrix
$\calR(s)$ of $\Yhg$ \cite{frenkel-reshetikhin,smirnov}.

\subsection{} 

Our equivalence implies in particular that the monodromy
of these difference equations, a meromorphic function of
the spectral parameter $\zeta=\Exp{s}$, is explicitly expressed in terms
of the abelianisation of the universal $R$--matrix $\Rq(\zeta)$
of $\qloop$. The latter result is a version of the \KD theorem
for abelian $q$KZ equations.

This result was conjectured by \FR \cite{frenkel-reshetikhin}
for the non--abelian $q$KZ equations, and proved in the
rational and trigonometric cases by \TV when $\g=\sl_2$,
and attention is restricted to evaluation representations
with generic highest weights \cite{tarasov-varchenko-trigonometric,
tarasov-varchenko-elliptic}.

One difficulty in addressing the general case is that no functorial
way of relating arbitrary representations of $\Yhg$ and $\qloop$
was known to exist outside of type $\sfA$ prior to \cite{sachin-valerio-2}.
\footnote{For $\g=\sl_2$, evaluation representations of $\Yhg$
can be explicitly deformed to representations of $\qloop$. More
generally, in type $\sfA_n$, Moura proved the \KD Theorem for
the trigonometric $q$KZ equations with values in the vector
representation of $\qloop$ \cite{mouraKD} and, jointly with
Etingof, used this to construct a functor from the \fd
representations of $\qloop$ arising from the RTT construction
to those of Felder's elliptic quantum group \cite{etingof-moura}.}
We shall prove the \KD theorem for the full (non--abelian) $q$KZ
equations for any $\g$ in a sequel to this paper \cite
{sachin-valerio-qKL}.\footnote{A more general result, for the Lie algebras associated 
by Maulik--Okounkov to quivers \cite{maulik-okounkov-qgqc}, was
independently announced by Okounkov \cite{okounkov-neu}.
It includes in particular the \KD theorem for the qKZ equations
corresponding to simply--laced Lie algebras.}

\subsection{} 

A crucial feature of our approach is that the relevant monoidal
structures arise from the {\it deformed Drinfeld coproducts} on
$\Yhg$ and $\qloop$, rather than from the standard (Kac--Moody)
ones.\footnote{The relation to the standard coproduct is discussed
in \ref{ss:conj}.} 
The Drinfeld coproduct was defined for $\qloop$ by Drinfeld
\cite{drinfeld-new2}, and regularised through deformation by
Hernandez \cite{hernandez-affinizations,hernandez-drinfeld-coproduct}. 
Whereas this coproduct has long been understood to arise
from the polarisation of the loop algebra $\g((z))$ given by
\[\g((z))=\Bigl(\n_-((z))\oplus z^{-1}\h[z^{-1}]\Bigr)\oplus
\Bigl(\h[[z]]\oplus\n_+((z))\Bigr),\] 
a proper understanding of the structure it confers \fd
representations has been lacking so far.

We define a similar deformed coproduct on $\Yhg$, and show
that these endow $\Rfyang$ and $\Rfloop$ with the structures
of {\it meromorphic tensor categories} that is, roughly speaking,
categories endowed with a monoidal structure and associativity
constraints depending meromorphically on parameters. This
notion was outlined by \FR who used the term {\it analytic tensor
categories} \cite[p. 49]{frenkel-reshetikhin}, and formalised
by Soibelman to describe the structure of \fd representations
of $\qloop$ corresponding to the standard 
coproduct and the universal $R$--matrix $\Rq(\zeta)$ \cite
{soibelman-meromorphic}. Our observation that the deformed
Drinfeld coproduct fits within, and provides new examples of
such categories seems to be new.

Our first main result may be succintly stated as saying that
$\Fh{}$ is a meromorphic tensor functor.

\subsection{} 

Our second main result is that, when $|q|\neq 1$, $\Fh{}$
is a {\it braided} meromorphic tensor functor. In more detail,
$\Rfyang$ and $\Rfloop$ are known not to be braided
tensor categories. As pointed out, however, the universal
$R$--matrix $\Rq(\zeta)$ of $\qloop$ defines a meromorphic
commutativity constraint on $\Rfloop$ with respect to the
standard tensor product $\otimes$.

We show that the same holds \wrt the deformed Drinfeld
tensor product $\Dotimes{\zeta}$, {\it provided} $\Rq(\zeta)$
is replaced by the diagonal component $\Rq_0(\zeta)$ of
its Gauss decomposition. Thus, $\Rfloop$ may be endowed
with two distinct structures of meromorphic braided tensor
category, namely as
\[(\Rfloop,\otimes,\Rq(\zeta))\aand
(\Rfloop,\Dotimes{\zeta},\Rq_0(\zeta))\]
The latter structure does not seem to have been noticed
before, though it should be closely related to the large volume
limit in quantum cohomology.\footnote{These two structures
are, in fact, meromorphically equivalent, see \ref{ss:conj}.}

We prove a similar result for the Yangian by constructing the
commutative part $\calR_0(s)$ of its universal $R$--matrix,
and showing that it defines commutativity constraints for the
deformed Drinfeld tensor product $\Dotimes{s}$ of $\Yhg$.
The construction of $\calR_0(s)$ is more delicate than that
of $\Rq_0(\zeta)$, and involves a non--trivial analytic regularisation
of the formal infinite product formulae for $\calR_0(s)$ obtained
by \KT \cite{khoroshkin-tolstoy}.

We then show that
\[\Fh{}:(\Rfyang,\Dotimes{s},\calR_0(s))\longrightarrow
(\Rfloop,\Dotimes{\zeta},\Rq_0(\zeta))\]
is compatible with the meromorphic braiding. As mentioned
above, this implies in particular the \KD theorem for the
additive $q$KZ equations defined by $\calR_0(s)$, namely
the fact that their monodromy is expressed in terms of
$\Rq_0(\zeta)$, where $\zeta=\Exp{s}$.

\section{Statement of main results}

This section contains a more detailed description of our main
results, their background, and a sketch of some of their proofs.

\subsection{The deformed Drinfeld coproduct of $\qloop$} 

The Drinfeld coproduct on $\qloop$ was defined by Drinfeld in \cite
{drinfeld-new2}, and involves formal infinite sums of elements in $\qloop^{\otimes 2}$.
Composing with the $\C^\times$--action on the first factor, Hernandez
obtained a deformed coproduct, which is an algebra homomorphism 
\[\Delta_\zeta:\qloop\to\qloop((\zeta^{-1}))\otimes\qloop\]
where $\zeta$ is a formal variable \cite[\S 6]{hernandez-affinizations}.
The map $\Delta_\zeta$ is coassociative, in the sense that $\Delta_
{\zeta_1}\otimes\id\circ\Delta_{\zeta_2}=\id\otimes\Delta_{\zeta_2}
\circ\Delta_{\zeta_1\zeta_2}$ \cite[Lemma 3.4]{hernandez-drinfeld-coproduct}.

When computed on the tensor product of two \fd representations $\V_1,
\V_2$ of $\qloop$, the deformed Drinfeld coproduct $\Delta_{\zeta}$ is
analytically well--behaved. Specifically, the action of $\qloop$ on $\V_1((\zeta
^{-1}))\otimes\V_2$ is the Laurent expansion at $\zeta=\infty$ of a family
of actions of $\qloop$ on $\V_1\otimes\V_2$, whose matrix coefficients
are rational functions of $\zeta$ \cite[Lemma 3.10] {hernandez-drinfeld-coproduct}.
We denote $\V_1\otimes\V_2$ endowed with this action by $\V_1
\Dotimes{\zeta}\V_2$.

\subsection{} 

In Section \ref{sec: dt}, we give simple contour integral formulae
for the action of $\qloop$ on $\V_1\Dotimes{\zeta}\V_2$. These
yield an alternative proof of the rationality of $\Dotimes{\zeta}$,
as well as an explicit determination of its poles as a function of
$\zeta$. 

Specifically, let $\V$ be a \fd representation of $\qloop$, $\bfI$
the set of vertices of the Dynkin diagram of $\g$, $\{\Psi_i(z),\X
_i^\pm(z)\}_{i\in\bfI}$ the $\End(\V)$--valued rational functions
of $z\in\IP^1$ whose Taylor expansion at $z=\infty,0$ give the
action of the generators of $\qloop$ on $\V$ (see Section
\ref{ss:rationality}), and $\spec(\V)\subset\C^\times$ the set
of poles of these functions.

Let $\V_1,\V_2\in\Rloop$, and let $\zeta\in\C^\times$ be such
that $\zeta\sigma(\V_1)$ and $\sigma(\V_2)$ are disjoint. Then,
the action of $\qloop$ on $\V_1\Dotimes{\zeta}\V_2$ is given by
the following formulae for any $m\in\N$ and $k\in\Z$
\begin{align*}
\Delta_{\zeta}(\Psi^\pm_{i,\pm m}) &=
\sum_{p_1+p_2=m} \zeta^{\pm p_1}\Psi^\pm_{i,\pm p_1}\otimes\Psi^\pm_{i,\pm p_2}\\
\Delta_{\zeta}(\X^+_{i,k}) &= 
\zeta^k \X^+_{i,k}\otimes 1 + \oint_{C_2} 
\Psi_i(\zeta^{-1}w)\otimes \X^+_i(w) w^{k-1}dw\\
\Delta_{\zeta}(\X^-_{i,k}) &= 
\oint_{C_1} 
\X^-_i(\zeta^{-1}w)\otimes \Psi_i(w)w^{k-1}dw
+ 1\otimes\X^-_{i,k}
\end{align*}
where
\begin{itemize}
\item $C_1,C_2\subset\C^\times$ are Jordan curves which do not enclose $0$.
\item $C_1$ encloses $\zeta\spec(\V_1)$ and none of the points in $\spec(\V_2)$.
\item $C_2$ encloses $\spec(\V_2)$ and none of the points in $\zeta\spec(\V_1)$.
\end{itemize}

Direct inspection shows that $\Delta_\zeta$ is a rational function
of $\zeta$, with poles contained in $\spec(\V_2)\spec(\V_1)^{-1}$.

\subsection{}

A remarkable feature of the deformed Drinfeld coproduct
$\Dotimes{\zeta}$ is that it endows $\Rloop$ with the structure
of a meromorphic tensor category in the sense of \cite
{soibelman-meromorphic}. This category is strict, in that for
any $\V_1,\V_2,\V_3\in\Rloop$, the identification of vector
spaces
\[(\V_1\Dotimes{\zeta_1}\V_2)\Dotimes{\zeta_2}\V_3=
\V_1\Dotimes{\zeta_1\zeta_2}(\V_2\Dotimes{\zeta_2}\V_3)\]
intertwines the action of $\qloop$.\footnote{Readers unfamiliar
with the associativity identity above may note that it is also
satisfied by the (holomorphic) tensor product defined by
$\V_1\odot_{\zeta}\V_2=\V_1(\zeta)\otimes\V_2$, where
$\otimes$ is the standard tensor product, and $\V_1(\zeta)$
is the pull--back of $\V_1$ by the $\nC$--action on $\qloop$
by dilations.}

Meromorphic braided tensor categories were introduced by
Soibelman in \cite{soibelman-meromorphic} to formalise the
structure of the category of \fd representations of $\qloop$
endowed with the standard (Kac--Moody) tensor product and
the universal $R$--matrix $\Rq(\zeta)$. 

\subsection{The deformed Drinfeld coproduct of $\Yhg$} \label{ss:Drinfeld for Y}

A Drinfeld coproduct was conjecturally defined for the double
Yangian $\DYhg$ by Khoroshkin--Tolstoy \cite{khoroshkin-tolstoy}.
Like its counterpart for $\qloop$, it involves formal infinite sums.
Moreover, the Yangian $\Yhg\subset\DYhg$ is not closed under it.

By degenerating our contour integral formulae for $\Dotimes{\zeta}$,
we obtain in Section \ref{ssec: dt-yangian} a family of actions $V_1
\Dotimes{s}V_2$ of $\Yhg$ on the tensor product of two \fd
representations of $\Yhg$, which is a rational function of a
parameter $s\in\C$. Its expansion at $s=\infty$ should coincide
with the deformation of the Drinfeld coproduct on $\DYhg$ via the
translation action of $\C$ on $\Yhg$, once the negative modes
of $\DYhg$ are reexpressed in terms of the positive ones through
a Taylor expansion of the corresponding generating functions.

Let $\{\xi_{i,r},x_{i,r}^{\pm}\}_{i\in\bfI, r\in\N}$ be the loop generators
of $\Yhg$ (see \cite{drinfeld-yangian-qaffine}, or \S \ref{sec: yqla}
for definitions). On a \fd representation $V$, the generating series
\[\xi_i(u) = 1+\hbar\sum_{r\geq 0} \xi_{i,r}u^{-r-1} \aand
x_i^{\pm}(u) = \hbar\sum_{r\geq 0} x_{i,r}^{\pm} u^{-r-1}\]
are expansions at $u=\infty$ of $\End(V)$--valued rational functions
\cite[Prop. 3.6]{sachin-valerio-2}. Let $\spec(V)\subset\C$ be the
he set of poles of the functions $\{x_i^\pm(u),\xi_i(u)^{\pm 1}\}_
{i\in\bfI}$ on $V$.

Let $s\in\C$ be such that $\spec(V_1)+s$ and $\spec(V_2)$ are
disjoint. Then, the action of $\Yhg$ on $V_1\Dotimes{s}V_2$ is
given by
\begin{align*}
\DD{s}(\xi_{i,r}) &=
\tau_s(\xi_{i,r})\otimes 1+
\hbar\sum_{p_1+p_2=r-1}\tau_s(\xi_{i,p_1})\otimes\xi_{i,p_2}+1\otimes\xi_{i,r}\\
\DD{s}(x_{i,r}^+) &=
\tau_s(x^+_{i,r})\otimes 1
+
\hbar^{-1}\oint_{C_2} 
\xi_i(v-s)\otimes x_i^+(v)v^r dv\\
\DD{s}(x_{i,r}^-) &= \hbar^{-1}\oint_{C_1} x_i^-(v-s)\otimes \xi_i(v)v^r dv
+ 1\otimes x_{i,r}^-
\end{align*}
where $\tau_s$ is the translation automorphism of $\Yhg$ given by
\[\tau_s(\xi_i(u))=\xi_i(u-s)\aand\tau_s(x_i^\pm(u))=x_i^\pm(u-s)\]
and $C_1,C_2$ are Jordan curves such that
\begin{itemize}
\item $C_1$ encloses $\spec(V_1)+s$ and none of the points in $\spec(V_2)$.
\item $C_2$ encloses $\spec(V_2)$ and none of the points in $\spec(V_1)+s$.
\end{itemize}

As for $\qloop$, direct inspection shows that the action of $\Yhg$
on $V_1\Dotimes{s}V_2$ is a rational function of $s$, with poles
contained in $\spec(V_2)-\spec(V_1)$. Moreover, the tensor product
$\Dotimes{s}$ gives $\Ryang$ the structure of a meromorphic
tensor category, which is strict in the sense that, for any $V_1,
V_2,V_3\in\Ryang$
\[(V_1\Dotimes{s_1}V_2)\Dotimes{s_2}V_3=
V_1\Dotimes{s_1+s_2}(V_2\Dotimes{s_2}V_3)\]

\subsection{Meromorphic tensor structure on $\Fh{}$} \label{ss:nc}

Recall the notion of {\it \nc} representation of $\Yhg$ \cite[\S 5.1]{sachin-valerio-2}.
Let $\{\xi_i(u),x_i^\pm(u)\}_{i\in\bfI}$ be the generating functions
defined in \ref{ss:Drinfeld for Y}. $V$ is called \nc if, for any
$i\in\bfI$, the poles of $x_i^+(u)$ (resp. $x_i^-(u)$) do not
differ by non--zero integers. If $V$ is non--congruent, the
monodromy of the difference equations defined by the
commuting fields $\xi_i(u)$ may be used to define an
action of $\qloop$ on the vector space $\Fh{}(V)=V$ \cite
{sachin-valerio-2}.

\subsection{}

If $V_1,V_2\in\Ryang$ are non--congruent, the Drinfeld tensor
product $V_1\Dotimes{s}V_2$ is generically non--congruent in
$s$. Our first main result is the following (see Theorem \ref{thm: tensor}).

\begin{thm}\label{th:main 1}\hfill
\begin{enumerate}
\item There exists a meromorphic $GL(V_1\otimes V_2)$--valued
function $\J_{V_1,V_2}(s)$, which is natural in $V_1,V_2$, and
such that
\[\J_{V_1,V_2}(s):
\Fh{}(V_1)\Dotimes{\zeta}\Fh{}(V_2)\longrightarrow\Fh{}(V_1\Dotimes{s} V_2)\]
is an isomorphism of $\qloop$--modules, where $\zeta=e^{2\pi
\iota s}$.
\item $\J$ is a meromorphic tensor structure on $\Fh{}$. That is,
for any non--congruent $V_1,V_2,V_3\in\Ryang$, the following
is a commutative diagram 
\[\xymatrix{
(\Fh{}(V_1)\Dotimes{\zeta_1}\Fh{}(V_2))\Dotimes{\zeta_2}\Fh{}(V_3)\ar[dd]_{\J_{V_1,V_2}(s_1)\otimes\id}
\ar@{=}[r]&
\Fh{}(V_1)\Dotimes{\zeta_1\zeta_2}(\Fh{}(V_2)\Dotimes{\zeta_2}\Fh{}(V_3))\ar[dd]^{\id\otimes\J_{V_2,V_3}(s_2)}\\ & \\
\Fh{}(V_1\Dotimes{s_1}V_2)\Dotimes{\zeta_2}\Fh{}(V_3)\ar[dd]_{\J_{V_1\Dotimes{s_1}V_2,V_3}(s_2)}
&
\Fh{}(V_1)\Dotimes{\zeta_1\zeta_2}\Fh{}(V_2\Dotimes{s_2}V_3)\ar[dd]^{\J_{V_1,V_2\Dotimes{s_2}V_3}(s_1+s_2)}\\ & \\
\Fh{}((V_1\Dotimes{s_1}V_2)\Dotimes{s_2}V_3)
\ar@{=}[r]&
\Fh{}(V_1\Dotimes{s_1+s_2}(V_2\Dotimes{s_2}V_3))
}\]
where $\zeta_i=\exp(2\pi\iota s_i)$.
\end{enumerate}
\end{thm}

\subsection{}\label{ss:J(s)}

Just as the functor $\Fh{}$ is governed by the abelian, additive
difference equations defined by the commuting fields $\xi_i(u)$
of the Yangian, the tensor structure $\J(s)$ arises from another
such difference equation, namely an appropriate abelianisation
of the $\qKZ$ equations on $V_1\otimes V_2$ \cite{frenkel-reshetikhin,smirnov}.

Specifically, let
\[\calR^0(s)=1+\hbar\frac{\Omega_\h}{s}+\cdots\]
be the diagonal part in the Gauss decomposition of the universal
$R$--matrix of $\Yhg$ acting on $V_1\otimes V_2$, where $\Omega
_\h\in\h\otimes\h$ is the Cartan part of the Casimir tensor of $\g$
\cite{khoroshkin-tolstoy}. Unlike the analogous case of $\qloop$
\cite{kazhdan-soibelman,etingof-moura}, the expansion of $\calR
^0(s)$ does {\it not} converge near $s=\infty$. Indeed, when evaluated
on the tensor product of highest--weight vectors of two \fd irreducible
representations of $\Yhg$, this series is given by the Stirling expansion
of a ratio of Gamma functions \cite[Theorem 7.2]{khoroshkin-tolstoy},
which is known not to converge. We show, however, that $\calR
^0(s)$ possesses two {\it distinct}, meromorphic regularisations
$\Rud(s)$ in \S \ref{sec: R0}. These are asymptotic to $\calR^0(s)$
in the half--planes $\pm\Re(s/\hbar)\gg 0$, and are related by the
unitarity constraint $\Ru(s)\Rd(-s)^{21}=1$.

Each $\Rud(s)$ gives rise to the abelian $\qKZ$ equation
\[\Phi\ud(s+1)=\Rud(s)\Phi\ud(s)\]
where $\Phi\ud$ is an $\End(V_1\otimes V_2)$--valued function
of $s$. This equation admits a canonical right fundamental solution
$\Phi^{\pm}_+(s)$, which is holomorphic and invertible on an obtuse
sector contained inside the half--plane $\Re(s)\gg 0$, 
and possesses an asymptotic expansion of the form $(1+O
(s^{-1}))s^{\hbar\Omega_\h}$ within it (see Proposition \ref
{pr:left right}). The tensor structure $\J_{V_1,V_2}(s)$ may
be taken to be $\Phi^+_+(s+1)^{-1}$ or $\Phi^-_+(s+1)^{-1}$,
and is a regularisation of the infinite product
\[ \cdots\Rud(s+3)\Rud(s+2)\Rud(s+1) \]
Specifically,
\[\J^\pm_{V_1,V_2}(s) = e^{\hbar\gamma\Omega_{\h}}
\stackrel{\longleftarrow}{\prod}_{m\geq 1}
\Rud(s+m)e^{-\frac{\hbar\Omega_{\h}}{m}}\]
where $\gamma=\lim_n(1+1/2+\cdots+
1/n-\log(n))$ is the Euler--Mascheroni constant.

\subsection{Regularisation of $\calR^0(s)$}\label{ss:begin R0}

As mentioned above,  the abelian $R$--matrix $\calR^0(s)$ needs
to be regularised. A conjectural construction of $\calR^0(s)$ as a
formal infinite product with values in the double Yangian $\DYhg$
was given by Khoroshkin--Tolstoy \cite[Thm. 5.2]{khoroshkin-tolstoy}.
To make sense of this product, we notice in Section \ref{sec: R0}
that $\calR^0(s)$ formally satisfies an abelian additive difference
equation whose step is a multiple of $\hbar$.\footnote{This equation
should in turn be a consequence of the (non--linear) difference
equation satisfied by the full $R$--matrix of $\Yhg$ obtained
from crossing symmetry.} We then prove that the coefficient
matrix $\A(s)$ of this equation can be interpreted as a rational
function of $s$, and define $\Rud(s)$ as the canonical
fundamental solutions of the difference equation. Let us outline
this approach in more detail.

\subsection{}\label{ssec: intro-A}

Let $b_{ij}=d_ia_{ij}$ be the entries of the symmetrized Cartan matrix
of $\g$. Let $T$ be an indeterminate, and $\bfB(T)=([b_{ij}]_T)$ the
corresponding matrix of $T$--numbers. Then, there exists an integer
$l=m h^\vee$, which is a multiple of the dual Coxeter number $h^\vee$
of $\g$, and is such that $\bfB(T)^{-1}=[l]_T^{-1}\bfC(T)$, where the
entries of $\bfC(T)$ are Laurent polynomials in $T$ with coefficients
in $\N$ \cite{khoroshkin-tolstoy}.

Consider the following $GL(V_1\otimes V_2)$--valued function of
$s\in\C$
\[\A(s)=\exp\lp-
\sum_{\begin{subarray}{c} i,j\in\bfI\\ r\in\Z\end{subarray}}c_{ij}^{(r)}
\oint_{\calC} t_i'(v)\otimes t_j\lp v+s+\frac{(l+r)\hbar}{2}\rp\, dv \rp\]
where
\begin{itemize}
\item $c_{ij}(T)=\sum_{r\in\Z}c_{ij}^{(r)}T^r$ are the entries of $\bfC(T)$.
\item the contour $\calC$ encloses the poles of $\xi_i(u)^{\pm 1}$ 
on $V_1$.
\item $t_i(u) = \log(\xi_i(u))$ is defined by choosing a branch of the
logarithm.
\item $s\in\C$ is such that $v\to t_j(v+s+(l+r)\hbar/2)$ is analytic on
$V_2$ within $\calC$, for every $j\in\bfI$ and $r\in\Z$ such that $c_
{ij}^{(r)}\neq 0$.
\end{itemize}
We prove in Section \ref{ssec: op-A} that $\A$ extends to a rational
function of $s$ which has the following expansion near $s=\infty$
\[\A(s)=1-l\hbar^2\frac{\Omega_{\h}}{s^2} + O(s^{-3})\]

\subsection{}\label{ssec: intro-R}

The infinite product $\calR^0(s)$ considered in \cite{khoroshkin-tolstoy}
formally satisfies
\[\calR^0(s+l\hbar) = \A(s)\calR^0(s)\]
This difference equation is regular (that is, the coefficient of $s^{-1}$
in the expansion of $\A(s)$ at $s=\infty$ is zero), and therefore admits
two canonical meromorphic fundamental solutions $\Rud(s)$. The
latter are uniquely determined by the requirement that they be
holomorphic and invertible for $\pm\Re(s/\hbar)\gg 0$, and asymptotic
to $1+O(s^{-1})$ as $s\to\infty$ in that domain (see \eg \cite{birkhoff-difference,
borodin,krichever} or \cite[\S 4]{sachin-valerio-2}). Explicitly,
\begin{align*}
\Ru(s) &= \stackrel{\longrightarrow}{\prod}_{n\geq 0} \A(s+nl\hbar)^{-1} \\
\Rd(s) &= \stackrel{\longrightarrow}{\prod}_{n\geq 1} \A(s-nl\hbar) 
\end{align*}

The functions $\Rud(s)$ are distinct regularisations of $\calR^0(s)$,
and are related by the unitarity constraint
\[\Ru_{V_1,V_2}(s)\Rd_{V_2,V_1}(-s)^{21}=1\]
We show in Theorem \ref{thm: R0} that they define meromorphic
commutativity constraints on $\Ryang$ endowed with the Drinfeld
tensor product $\Dotimes{s}$.

\subsection{Kohno--Drinfeld theorem for abelian $q$KZ equations}

Our second main result is a Kohno--Drinfeld theorem for the abelian,
additive $q$KZ equations defined by $\Rud(s)$. Together with Theorem
\ref{th:main 1}, it establishes an equivalence of meromorphic braided
tensor categories between $\Ryang$ and $\Rloop$ akin to the \KL
equivalence between the affine Lie algebra $\wh{\g}$ and corresponding
quantum group $U_q\g$.

Fix $V_1,\ldots,V_n\in\Ryang$. The abelian $q$KZ equations are the
integrable system of additive difference equations for a meromorphic
function $F:\C^n\to\End(V_1\otimes\cdots\otimes V_n)$ which are
given by \cite{frenkel-reshetikhin,smirnov}
\begin{equation}\label{eq:qkz intro}
F(\ul{s}+e_i) = A_i(\ul{s})F(\ul{s})
\end{equation}
where $\ul{s} = (s_1,\ldots,s_n)\in\C^n$, $\{e_i\}_{i=1}^n$ is the
standard basis of $\C^n$, and $A_i(\ul{s})$ is given by
\begin{multline*}
A_i(\ul{s}) = \Rud_{i-1,i}(s_{i-1}-s_i-1)^{-1}\cdots \Rud_{1,i}(s_1-s_i-1)^{-1}\\
\cdot\Rud_{i,n}(s_i-s_n)\cdots \Rud_{i,i+1}(s_i-s_{i+1})
\end{multline*}
with $\Rud_{i,j}=\Rud_{V_i,V_j}$  the regularisations of the
commutative $R$--matrix of $\Yhg$ described in \ref{ss:begin R0}.

These equations admit a set of fundamental solutions $\Phi^\pm_
\sigma$ which generalise the right/left solutions in the $n=2$ case.
They are parametrised by permutations $\sigma\in\Sym_n$, and
have prescribed asymptotic behaviour when $s_i-s_j\to\infty$ for
any $i<j$, in such a way that $\Re(s_{\sigma^{-1}(i)}-s_{\sigma^
{-1}(j)})\gg 0$. By definition, the monodromy of \eqref{eq:qkz intro}
is the 2--cocyle on $\Sym_n$ with values in the group of meromorphic
$GL(V_1\otimes \cdots\otimes V_n)$--valued functions 
of the variables $\zeta_i=\Exp{s_i}$ given by 
\[M^\pm_{\sigma,\tau}(\ul{s})=(\Phi^\pm_\sigma(\ul{s}))^{-1}\cdot\Phi^\pm_\tau(\ul{s})\]

\subsection{}

A Kohno--Drinfeld theorem for the $q$KZ equations determined by
the full (non--abelian) $R$--matrix of $\Yhg$ was conjectured by \FR
\cite[\S 6]{frenkel-reshetikhin}. It states that the monodromy of \eqref
{eq:qkz intro}, with $\calR^0$ replaced by $\calR$, is given by 
the universal $R$--matrix $\Rq(\zeta)$ of
$\qloop$ acting on a tensor product of suitable $q$--deformations
of $V_1,\ldots,V_n$.

Assuming that $|q|\neq 1$, we prove this theorem for the abelian $q$KZ
equations determined by $\Rud$. To this end, we first construct the
commutative part $\Rq^0(\zeta)$ of
the $R$--matrix of $\qloop$ in \S \ref{sec: R0qla} by following a procedure
similar to that described in \ref{ss:begin R0}--\ref{ssec: intro-R}. Namely,
we start from Damiani's formula for $\Rq^0(\zeta)$ \cite{damiani}, show
that if formally satisfies a regular $q$--difference equation with respect
to the parameter $\zeta$, and deduce from this that it is the expansion
at $\zeta=0$ of the corresponding canonical solution (unlike the case
of $\Yhg$, no regularisation of $\Rq^0(\zeta)$ is necessary here).
We also show that $\Rq^0(\zeta)$ defines meromorphic commutativity
constraints on $\Rloop$ endowed with the deformed Drinfeld coproduct.

We then prove the following (Theorem \ref{ss:KD})
\begin{thm}\label{th:main 2}
Assume that $|q|\neq 1$,
and set
\[\veps=\left\{\begin{array}{ll}+&\text{if $|q|<1$}\\[1.1ex]-&\text{if $|q|>1$}\end{array}\right.\]
Let $V_1,\ldots,V_n\in\Ryang$ be non--congruent, and let $\V
_\ell =\Fh{}(V_\ell)$ be the corresponding representations of $\qloop$.

Then, the monodromy of the abelian $q$KZ equations determined by
$\calR^{0,\veps}(s)$ on $V_1\otimes\cdots\otimes V_n$ is given by
$\Rq^0(\zeta)$. Specifically, the following holds for any $\sigma\in\Sym
_n$ and $i=1,\ldots, n-1$ such that $\sigma^{-1}(i)<\sigma^{-1}(i+1)$, 
\[
(\Phi^{\ve}_{\sigma}(\ul{s}))^{-1}\cdot\Phi^{\ve}_{(i\,i+1)\sigma}(\ul{s}) = 
\Rq^{0}
_{\V_i,\V_{i+1}}(\zeta_i\zeta_{i+1}^{-1})
\]
\end{thm}
The same result holds for the monodromy of the $q$KZ equations
determined by $\calR^{0,-\veps}(s)$, provided $\Rq^0(\zeta)$ is
replaced by $\Rq^0_{21}(\zeta^{-1})^{-1}$.\footnote{Theorem
\ref{ss:KD} contains both of these statements in a uniform fashion.
Thus $\Rq^0(\zeta)$ of Theorem \ref{th:main 2} above is $\Rq^{0,
\veps}(\zeta)$ of Theorem \ref{ss:KD} with $\veps = \pm$ according
to the statement above.}

\subsection{Relation to the \KM coproduct}\label{ss:conj}

We conjecture that the twist $\J(s)$ also yields a non--meromorphic
tensor structure on the functor $\Fh{}$, when the categories $\Ryang$
and $\Rloop$ are endowed with the standard monoidal structures
arising from the Kac--Moody coproducts on $\Yhg,\qloop$.

More precisely, the Drinfeld and Kac--Moody coproducts on $\qloop$
are related by a meromorphic twist, given by the lower triangular part
$\Rq^{\qloop}_-(\zeta)$ of the universal $R$--matrix \cite{ekp}. A similar
statement holds for $\Yhg$ \cite{sachin-valerio-qKL}. Composing, we obtain
a meromorphic tensor structure $J(s)$ on $\Fh{}$ relative to the standard
monoidal structures
\[\xymatrix{
\Fh{}(V_1)(\zeta)\otimes\Fh{}(V_2) \ar[rr]^{\calR^{\qloop}_-(\zeta)}\ar[dd]_{J_{V_1,V_2}(s)}
&& \Fh{}(V_1)\Dotimes{\zeta}\Fh{}(V_2)\ar[dd]^{\J_{V_1,V_2}(s)}\\
&&\\
\Fh{}(V_1(s)\otimes V_2) \ar[rr]^{\calR^{\Yhg}_-(s)} && \Fh{}(V_1\Dotimes{s}V_2)
}\]
We conjecture that $J_{V_1,V_2}(s)$ is holomorphic in $s$, and can therefore
be evaluated at $s=0$, thus yielding a tensor structure on $\Fh{}$ with respect
to the standard coproducts. We will return to this in \cite{sachin-valerio-qKL}.

\subsection{Extension to arbitrary \KM algebras}

The results of \cite{sachin-valerio-2} hold for an arbitrary symmetrisable
\KM algebra $\g$, provided one considers the categories of 
representations of $\Yhg$ and $\qloop$ whose restriction to $\g$ and $\Uqg$
respectively are integrable and in category $\calO$. Although we restricted
ourselves to the case of a \fd semisimple $\g$ in this paper, our results
on the Drinfeld coproducts of $\Yhg$ and $\qloop$ are valid for an arbitrary
$\g$, and it seems likely that the same should hold for the construction of
the tensor structure $\J(s)$. The main obstacle in working in this generality
is the construction and regularisation of $\calR^0(s)$ for an arbitrary $\g$.
Once this is achieved, the proof of Theorems \ref{th:main 1} and \ref
{th:main 2} carries over verbatim.

\subsection{Outline of the paper}

In Section \ref{sec: yqla}, we review the definitions of $\Yhg$ and $\qloop$.
Section \ref{sec: dt} is devoted to defining the Drinfeld coproducts on $\qloop$
and $\Yhg$. We give a construction of the diagonal part $R^0$ of the $R$--matrix
of $\Yhg$ in \S \ref{sec: R0}. Section \ref{sec: functor} reviews the definition
of the functor $\Fh{}$ given in \cite{sachin-valerio-2}. The construction of a
meromorphic tensor structure on $\Fh{}$ is given in \S \ref{sec: tensor}. In
Section \ref{sec: R0qla}, we show that, when $|q|\neq 1$, the commutative
part $\Rq^0(\zeta)$ of the $R$--matrix of $\qloop$ defines a meromorphic commutativity
constraint on $\Rloop$. Finally, in Section \ref{se:KD}, we prove a Kohno--Drinfeld
theorem for the abelian $q$KZ equations defined by $\calR^0(s)$. Appendix
\ref{se:KT} gives the inverses of all symmetrised $q$--Cartan matrices
of finite type. 

\subsection{Acknowledgments}

We are grateful to David Hernandez for his comments on an earlier
version of this paper, to Sergey Khoroshkin for correspondence
about the inversion of a $q$--Cartan matrix, and to Alexei Borodin
and Julien Roques for correspondence on the asymptotics of solutions
of difference equations. 
We are also extremely grateful to the anonymous referee for the careful
reading, comments and suggestions which helped improve the exposition.
Part of this paper was written while the first author visited IHES in the
summer of 2013. He is grateful to IHES for its invitation and wonderful
working conditions. 

\section{Yangians and quantum loop algebras}\label{sec: yqla}

\subsection{}\label{ssec: kma}

Let $\g$ be a complex, semisimple Lie algebra and $(\cdot,\cdot)$ the
invariant bilinear form on $\g$ normalised so that the squared length
of short roots is $2$. Let $\h\subset\g$ be a Cartan subalgebra of $\g$,
$\{\alpha_i\}_{i\in\bfI}\subset\h^*$ a basis of simple roots of $\g$ relative
to $\h$ and $a_{ij}=2(\alpha_i,\alpha_j)/(\alpha_i,\alpha_i)$ the entries
of the corresponding Cartan matrix $\bfA$. Set $d_i=(\alpha_i,\alpha_i)
/2\in\{1,2,3\}$, so that $d_ia_{ij}=d_j a_{ji}$ for any $i,j\in\bfI$. 

\subsection{The Yangian $\Yhg$}\label{ssec: yangian}

Let $\hbar\in\C$. The Yangian $\Yhg$ is the unital, associative
$\C$--algebra generated by elements $\{x^{\pm}_{i,r},\xi_{i,r}\}
_{i\in\bfI,r\in\N}$, subject to the following relations

\begin{enumerate}
\item[(Y1)] For any $i,j\in\bfI$, $r,s\in\N$
\[[\xi_{i,r}, \xi_{j,s}] = 0 \]
\item[(Y2)] For $i,j\in\bfI$ and $s\in \N$
\[[\xi_{i,0}, x_{j,s}^{\pm}] = \pm d_ia_{ij} x_{j,s}^{\pm}\]
\item[(Y3)] For $i,j\in\bfI$ and $r,s\in\N$
\[[\xi_{i,r+1}, x^{\pm}_{j,s}] - [\xi_{i,r},x^{\pm}_{j,s+1}] =
\pm\hbar\frac{d_ia_{ij}}{2}(\xi_{i,r}x^{\pm}_{j,s} + x^{\pm}_{j,s}\xi_{i,r})\]
\item[(Y4)] For $i,j\in\bfI$ and $r,s\in \N$
\[
[x^{\pm}_{i,r+1}, x^{\pm}_{j,s}] - [x^{\pm}_{i,r},x^{\pm}_{j,s+1}]=
\pm\hbar\frac{d_ia_{ij}}{2}(x^{\pm}_{i,r}x^{\pm}_{j,s} + x^{\pm}_{j,s}x^{\pm}_{i,r})
\]
\item[(Y5)] For $i,j\in\bfI$ and $r,s\in \N$
\[[x^+_{i,r}, x^-_{j,s}] = \delta_{ij} \xi_{i,r+s}\]
\item[(Y6)] Let $i\not= j\in\bfI$ and set $m = 1-a_{ij}$. For any
$r_1,\cdots, r_m\in \N$ and $s\in \N$
\[\sum_{\pi\in\Sym_m}
\left[x^{\pm}_{i,r_{\pi(1)}},\left[x^{\pm}_{i,r_{\pi(2)}},\left[\cdots,
\left[x^{\pm}_{i,r_{\pi(m)}},x^{\pm}_{j,s}\right]\cdots\right]\right]\right]=0\]
\end{enumerate}

\subsection{Remark}\label{rm:Y6}

By \cite[Lemma 1.9]{levendorskii}, the relation (Y6) follows from
(Y1)--(Y3) and the special case of (Y6) when $r_1=\cdots=r_m=0$.
In turn, the latter automatically holds on \fd representations of the
algebra defined by relations (Y2) and (Y5) alone (see, \eg \cite
[Prop. 2.7]{sachin-valerio-2}). Thus, a finite--dimensional
representation $V$ of $\Yhg$ is given by operators $\{\xi_{i,r},
x^\pm_{i,r}\}_{i\in\bfI,r\in\N}$ in $\End(V)$ satisfying relations
(Y1)--(Y5).

\subsection{}\label{ssec: formal-series-y}

Assume henceforth that $\hbar\neq 0$, and define $\xi_i(u),x^\pm_i(u)
\in\Yhg[[u^{-1}]]$ by
\[\xi_i(u)=1 + \hbar\sum_{r\geq 0} \xi_{i,r}u^{-r-1}\aand
x^{\pm}_i(u)=\hbar\sum_{r\geq 0} x_{i,r}^{\pm} u^{-r-1}\]

For an associative algebra $A$, we denote by $A[u,v; u^{-1},v^{-1}]]$
the algebra of formal series $\sum_{r,s} a_{r,s} u^r v^s$ for which
there exist $M,N\in\Z$ such that $a_{r,s}\neq 0$ implies $r\leq M$
and $s\leq N$.\\

\begin{prop}\cite[Prop. 2.3]{sachin-valerio-2}\label{pr:Y fields} 
The relations (Y1),(Y2)--(Y3),(Y4),(Y5) and (Y6) are respectively
equivalent to the following identities in $\Yhg[u,v;u^{-1},v^{-1}]]$
\begin{enumerate}
\item[($\Y$1)] For any $i,j\in\bfI$, 
\[[\xi_i(u), \xi_j(v)]=0\]\\[-4ex]
\item[($\Y$2)] For any $i,j\in\bfI$,
\[[\xi_{i,0},x^\pm_j(u)]=\pm d_ia_{ij}x^\pm_j(u)\]\\[-4ex]
\item[($\Y$3)] For any $i,j\in \bfI$, and $a = \hbar d_ia_{ij}/2$
\[(u-v\mp a)\xi_i(u)x_j^{\pm}(v)=
(u-v\pm a)x_j^{\pm}(v)\xi_i(u)\mp 2a x_j^{\pm}(u\mp a)\xi_i(u)\]
\\[-4ex]
\item[($\Y$4)] For any $i,j\in \bfI$, and $a = \hbar d_ia_{ij}/2$
\begin{multline*}
(u-v\mp a) x_i^{\pm}(u)x_j^{\pm}(v)\\
= (u-v\pm a)x_j^{\pm}(v)x_i^{\pm}(u)
+\hbar\lp [x_{i,0}^{\pm},x_j^{\pm}(v)] - [x_i^{\pm}(u),x_{j,0}^{\pm}]\rp
\end{multline*}\\[-4ex]
\item[($\Y$5)] For any $i,j\in \bfI$
\[(u-v)[x_i^+(u),x_j^-(v)]=-\delta_{ij}\hbar\left(\xi_i(u)-\xi_i(v)\right)\]
\item[($\Y$6)] For any $i\neq j\in\bfI$, $m=1-a_{ij}$,
\[\sum_{\pi\in\Sym_m}
\left[x^{\pm}_i(u_{\pi_1}),\left[x^{\pm}_i(u_{\pi(2)}),\left[\cdots,
\left[x^{\pm}_i(u_{\pi(m)}),x^{\pm}_j(v)\right]\cdots\right]\right]\right]=0\]
\end{enumerate}
\end{prop}

\begin{rem}\label{rk:Y fields}
Taking the coefficient of $u^0$ in relation ($\Y$3) gives
\[
\hbar\xi_{i,0}x_j^{\pm}(v) - vx_j^{\pm}(v) \mp a x_j^{\pm}(v) = 
\hbar x_j^{\pm}(v)\xi_{i,0} - vx_j^{\pm}(v) \pm a x_j^{\pm}(v)
\]
Thus we get $[\xi_{i,0},x_j^{\pm}(v)] = \pm d_ia_{ij} x_j^{\pm}(v)$ which is relation
($\Y$2).
\end{rem}

\subsection{Shift automorphism}\label{ssec: shift-yangian}

The group of translations of the complex plane acts on
$\Yhg$ by
\[\tau_a(y_r) = \sum_{s=0}^r
\left(\begin{array}{c}r\\s\end{array}\right)
a^{r-s}y_s\]
where $a\in\C$, $y$ is one of $\xi_i,x_i^\pm$. In terms of
the generating series introduced in \ref{ssec: formal-series-y},
\[\tau_a(y(u)) = y(u-a)\]
Given a representation $V$ of $\Yhg$ and $a\in \C$, set
$V(a)=\tau_a^*(V)$.

\subsection{Quantum loop algebra $\qloop$}\label{ssec: qla}

Let $q\in\C^\times$ be of infinite order. For any $i\in\bfI$, set $q
_i=q^{d_i}$. We use the standard notation for Gaussian integers
\begin{gather*}
[n]_q = \frac{q^n - q^{-n}}{q-q^{-1}}\\[.5ex]
[n]_q! = [n]_q[n-1]_q\cdots [1]_q\qquad
\qbin{n}{k}{q} = \frac{[n]_q!}{[k]_q![n-k]_q!}
\end{gather*}

The quantum loop algebra $\qloop$ is the unital, associative $\C$--algebra generated
by elements $\{\Psi_{i,\pm r}^\pm\}_{i\in\bfI,r\in\N}$, $\{\X_{i,k}^\pm\}
_{i\in\bfI,k\in\Z}$, subject to the following relations
\begin{itemize}
\item[(QL1)] For any $i,j\in\bfI$, $r,s\in\N$,
\begin{gather*}
[\Psi_{i,\pm r}^\pm,\Psi_{j,\pm s}^\pm]=0
\qquad
[\Psi_{i,\pm r}^\pm,\Psi_{j,\mp s}^\mp]=0
\qquad
\Psi_{i,0}^{\pm}\Psi_{i,0}^{\mp} = 1
\end{gather*}
\item[(QL2)] For any $i,j\in\bfI$, $k\in\Z$,
\[\Psi_{i,0}^+\X^\pm_{j,k}\Psi_{i,0}^-=q_i^{\pm a_{ij}}\X^\pm_{j,k}\]
\item[(QL3)] For any $i,j\in\bfI$, $\veps\in\{\pm\}$ and $l\in\Z$
\[\Psi^\veps_{i,k+1}\X^\pm_{j,l} - q_i^{\pm a_{ij}}\X^\pm_{j,l}\Psi^\veps_{i,k+1}
=
q_i^{\pm a_{ij}}\Psi^\veps_{i,k}\X^\pm_{j,l+1}-\X^\pm_{j,l+1}\Psi^\veps_{i,k}\]
for any $k\in\Z_{\geq 0}$ if $\veps=+$ and $k\in\Z_{<0}$ if $\veps=-$
\item[(QL4)] For any $i,j\in\bfI$ and $k,l\in \Z$
\[\X^\pm_{i,k+1}\X^\pm_{j,l} - q_i^{\pm a_{ij}}\X^\pm_{j,l}\X^\pm_{i,k+1}=
q_i^{\pm a_{ij}}\X^\pm_{i,k}\X^\pm_{j,l+1}-\X^\pm_{j,l+1}\X^\pm_{i,k}\]
\item[(QL5)] For any $i,j\in\bfI$ and $k,l\in \Z$
\[[\X^+_{i,k},\X^-_{j,l}] = \delta_{ij} \frac{\Psi^+_{i,k+l} - \Psi^-_{i,k+l}}{q_i - q_i^{-1}}\]
where $\Psi^{\pm}_{i,\mp k}=0$ for any $k\geq 1$.
\item[(QL6)] For any $i\neq j\in\bfI$, $m=1-a_{ij}$, $k_1,\ldots, k_m\in\Z$
and $l\in \Z$
\[\sum_{\pi\in \Sym_m} \sum_{s=0}^m (-1)^s\qbin{m}{s}{q_i}
\X^\pm_{i,k_{\pi(1)}}\cdots \X^\pm_{i,k_{\pi(s)}} \X^\pm_{j,l}\X^\pm_{i,k_{\pi(s+1)}}\cdots \X^\pm_{i,k_{\pi(m)}} = 0\]
\end{itemize}

\subsection{Remark}\label{rm:QL6}
By \cite[Lemma 2.12]{sachin-valerio-2}, the relation (QL6) follows
from (QL1)--(QL3) and the special case of (QL6) when $k_1=\cdots
=k_m=0$. In turn, the latter automatically holds on \fd representations
of the algebra defined by relations (QL2) and (QL5) alone (see, \eg
\cite[Prop. 2.13]{sachin-valerio-2}). Thus, a finite--dimensional
representation $\V$ of $\qloop$ is given by operators $\{\Psi_{i,\pm r}
^\pm,\X_{i,k}^\pm\}_{i\in\bfI,r\in\N,k\in\Z}$ in $\End(\V)$ satisfying
relations (QL1)--(QL5).

\subsection{}\label{ssec: qla-fields}
Define $\Psi_i(z)^+,\X^\pm_i(z)^+\in\qloop[[z^{-1}]]$ and
$\Psi_i(z)^-,\X^\pm_i(z)^-\in \qloop[[z]]$ by
\begin{align*}
\Psi_i(z)^+     &= \sum_{r\geq 0}\Psi_{i,r}^+ z^{-r} & \Psi_i(z)^-      &= \sum_{r\leq 0}\Psi_{i,r}^-z^{-r}\\
\X_i^\pm(z)^+   &= \sum_{r\geq 0}\X_{i,r}^\pm z^{-r}& \X_i^\pm(z)^-&=-\sum_{r< 0}\X_{i,r}^\pm z^{-r}
\end{align*}

\begin{prop}\cite[Prop. 2.7]{sachin-valerio-2}\label{pr:qloop fields}
The relations (QL1),(QL2)--(QL3),(QL4), (QL5),(QL6) imply the following
relations in $\qloop[z,w;z^{-1},w^{-1}]]$
\begin{enumerate}
\item[(\QL1)] For any $i,j\in\bfI$, 
\begin{gather*}
[\Psi_i(z)^+,\Psi_j(w)^+]=0
\end{gather*}
\item[(\QL2)] For any $i,j\in\bfI$,
\[\Psi_{i,0}^+\X_j^\pm(z)^+\Psi_{i,0}^-=q_i^{\pm a_{ij}}\X_j^\pm(z)^+\]
\item[(\QL3)] For any $i,j\in\bfI$
\begin{multline*}
(z-q_i^{\pm a_{ij}}w)\Psi_i(z)^+\X_j^\pm(w)^+\\
=(q_i^{\pm a_{ij}}z-w)\X_j^\pm(w)^+\Psi_i(z)^+-
(q_i^{\pm a_{ij}} - q_i^{\mp a_{ij}})q_i^{\pm a_{ij}}w\X_j^\pm(q_i^{\mp a_{ij}}z)^+\Psi_i(z)^+
\end{multline*}
\item[(\QL4)] For any $i,j\in\bfI$
\begin{multline*}
(z-q_i^{\pm a_{ij}}w)\X_i^\pm(z)^+\X_j^\pm(w)^+-
(q_i^{\pm a_{ij}}z-w)\X_j^\pm(w)^+\X_i^\pm(z)^+\\
=z\lp\X_{i,0}^\pm\X_j^\pm(w)^+-q_i^{\pm a_{ij}}\X_j^\pm(w)^+\X_{i,0}^\pm\rp
+ w\lp\X_{j,0}^\pm\X_i^\pm(z)^+-q_i^{\pm a_{ij}}\X_i^\pm(z)^+\X_{j,0}^\pm\rp
\end{multline*}
\item[(\QL5)] For any $i,j\in\bfI$
\[(z-w)[\X^+_i(z)^+,\X^-_j(w)^+] =
\frac{\delta_{ij}}{q_i-q_i^{-1}}\left(z\Psi_i(w)^+-w\Psi_i(z)^+-(z-w)\Psi_{i,0}^-\right)\]
\item[(\QL6)] For any $i\neq j\in\bfI$, and $m=1-a_{ij}$
\begin{multline*}
\sum_{\pi\in \Sym_m} \sum_{s=0}^m (-1)^s\qbin{m}{s}{q_i}
\X^\pm_i(z_{\pi(1)})^+\cdots \X^\pm_i(z_{\pi(s)})^+ \X^\pm_j(w)^+\\
\cdot \X^\pm_i(z_{\pi(s+1)})^+\cdots \X^\pm_i(z_{\pi(m)})^+ = 0
\end{multline*}
\end{enumerate}
\end{prop}

\subsection{Shift automorphism}\label{ssec: shift-qla}

The group $\nC$ of dilations of the complex plane acts on $\qloop$ by
\[\tau_{\alpha}(Y_k)=\alpha^k Y_k\]
where $\alpha\in\nC$, $Y$ is one of $\Psi^\pm_i,\X^\pm_i$. 
In terms of the generating series of \ref{ssec: qla-fields}, we
have
\[\tau_{\alpha}(Y(z)^\pm) = Y(\alpha^{-1}z)^\pm\]
Given a representation $\V$ of $\qloop$ and $\alpha\in\nC$, we
denote $\tau_{\alpha}^*(\V)$ by $\V(\alpha)$.

\subsection{Rationality}\label{ss:rationality}

The following rationality property is due to Beck--Kac \cite{beck-kac} and
Hernandez \cite{hernandez-drinfeld-coproduct} for $\qloop$ and to the
authors for $\Yhg$. In the form below, the result appears in \cite[Prop. 3.6]{sachin-valerio-2}.

\begin{prop}\label{prop: rationality}\hfill
\begin{enumerate}
\item Let $V$ be a $\Yhg$--module on which the operators $\{\xi_{i,0}\}_{i\in\bfI}$ act semisimply
with \fd weight spaces. Then, for every weight $\mu$ of $V$, the
generating series\\

$\ds{\xi_i(u)\in\End(V_{\mu})[[u^{-1}]]
\quad\text{and}\quad
x_i^{\pm}(u)\in\Hom(V_{\mu},V_{\mu\pm\alpha_i})[[u^{-1}]]}$\\

\noindent
defined in \ref{ssec: formal-series-y} are the expansions at $\infty$
of rational functions of $u$. Specifically, let $t_{i,1} = \xi_{i,1} - \ds
\frac{\hbar}{2} \xi_{i,0}^2\in\Yhg^\h$. Then,
\[x_i^{\pm}(u)= 2d_i\hbar  u^{-1}\lp 2d_i \mp \frac{\ad(t_{i,1})}{u}\rp^{-1}x_{i,0}^{\pm}\]
and
\[\xi_i(u) = 1 + [x_i^+(u),x_{i,0}^-]\]
\item Let $\V$ be a $\qloop$--module on which the operators $\{\Psi_{i,0}^{\pm}\}
_{i\in\bfI}$ act
semisimply with \fd weight spaces. Then, for every weight $\mu$ of $\V$
and $\veps\in\{\pm\}$, the generating series\\

$\ds{\Psi_i(z)^\pm\in\End(\V_\mu)[[z^{\mp 1}]]
\quad\text{and}\quad
\X_i^\veps(z)^\pm\in\Hom(\V_\mu,\V_{\mu+\veps\alpha_i}))[[z^{\mp 1}]]}$\\

\noindent
defined in \ref{ssec: qla-fields} are the expansions of rational functions
$\Psi_i(z),\X_i^\veps(z)$ at $z=\infty$ and $z=0$. Specifically, let
$H^\pm_{i,\pm 1}=\pm\Psi_{i,0}^{\mp}\Psi_{i,\pm 1}^{\pm}/(q_i-q_i^{-1})$.
Then,
\[\begin{split}
\X_i^\veps(z)
&=\phantom{-z}\lp 1- \veps\frac{\ad(H^+_{i,1})}{[2]_{q_i}z}\rp^{-1}\X_{i,0}^\veps\\
&=-z\lp 1- \veps z\frac{\ad(H_{i,-1}^-)}{[2]_{q_i}}\rp^{-1}\X_{i,-1}^\veps
\end{split}\]
and
\[\Psi_i(z) = \Psi_{i,0}^- + (q_i-q_i^{-1})[\X^+_i(z),\X^-_{i,0}]\]
\end{enumerate}
\end{prop}

\subsection{Poles of \fd representations}\label{ss:poles}

By Proposition \ref{prop: rationality}, we can define, for a given
$V\in\Ryang$, a subset $\spec(V)\subset\C$ consisting of the
poles of the rational functions $\xi_i(u)^{\pm 1}, x_i^{\pm}(u)$.

Similarly, for any $\V\in\Rloop$, we define a subset $\spec(\V)
\subset\nC$ consisting of the poles of the functions $\Psi_i(z)
^{\pm 1},\X^{\pm}_i(z)$.

\subsection{} 

The following is a direct consequence of Proposition \ref
{prop: rationality} and contour deformation. We set $\oint_\calC f=\frac{1}{2\pi\iota}\int_\calC
f$.
\begin{cor}\label{co:modes}\hfill
\begin{enumerate}
\item Let $V\in\Ryang$ and $\calC\subset\C$ be a Jordan curve
enclosing $\sigma(V)$.\footnote{By a Jordan curve, we shall
mean a disjoint union of simple, closed curves the inner domains
of which are pairwise disjoint.} Then, the following holds on $V$ for
any $r\in\N$
\[x_{i,r}^\pm=\frac{1}{\hbar}\oint_\calC x_i^\pm(u) u^r du
\aand
\xi_{i,r}=\frac{1}{\hbar}\oint_\calC \xi_i(u) u^r du\]
\item Let $\V\in\Rloop$ and $\calC\subset\C^\times$ be a Jordan
curve enclosing $\sigma(\V)$ and not enclosing $0$. Then, the
following holds on $\V$ for any $k\in\Z$ and $r\in\Ntimes$
\[\X_{i,k}^\pm=\oint_\calC\X_i^\pm(z)z^{k-1}dz
\quad\quad 
\Psi^\pm_{i,\pm r}=\pm\oint_\calC\Psi_i(z)z^{\pm r-1}dz\]
and
\[\oint_\calC\Psi_i(z)\frac{dz}{z}=\Psi_{i,0}^+-\Psi_{i,0}^-\]
\end{enumerate}
\end{cor}

\subsection{}\label{sec: imp-lem}

The following result will be needed later.

\begin{lem}\label{lem: imp-lem}
Let $V$ be a \fd representation of $\Yhg$ and $i,j\in\bfI$. If $u_0$ is a pole of
$x_j^{\pm}(u)$, then $\ds u_0 \pm \frac{\hbar d_ia_{ij}}{2}$ are poles of $\xi_i
(u)^{\pm 1}$.
\end{lem}

\begin{pf}
Consider the relation ($\mathcal{Y}3$) of Proposition \ref{pr:Y fields} (here $a = \hbar d_ia_{ij}/2$).
\begin{equation}\label{eq: commt-prime1}
\Ad(\xi_i(u))x_j^+(v) = \frac{u-v+a}{u-v-a}x_j^+(v) - \frac{2a}{u-v-a}x_j^+(u-a)\\
\end{equation}

Set $v=u+a$ to get $\Ad(\xi_i(u))x_j^+(u+a) = x_j^+(u-a)$. Combining this
with equation \eqref{eq: commt-prime1} above we get:

\begin{equation}\label{eq: commt-prime2}
\Ad(\xi_i(u))^{-1}x_j^+(v) = \frac{u-v-a}{u-v+a}x_j^+(v) + \frac{2a}{u-v+a}x_j^+(u+a)
\end{equation}

Differentiating \eqref{eq: commt-prime2} \wrt $v$ and then setting $v=u-a$ yields 

\begin{equation}\label{eq: commt-prime3}
2a \lp\Ad(\xi_i(u))\rp^{-1} \lp \frac{d}{du} x_j^+(u-a)\rp
 = x_j^+(u+a) - x_j^+(u-a)
\end{equation}

Differentiating \eqref{eq: commt-prime1} \wrt $u$, and combining equations
\eqref{eq: commt-prime2}, \eqref{eq: commt-prime3} with the following fact 
\[ \frac{d}{du}\Ad(\xi_i(u))x_j^+(v)=\Ad(\xi_i(u))\,[\xi_i(u)^{-1}\xi_i'(u),x_j^+(v)]\]
shows that
\begin{align}
[\xi_i(u)^{-1}\xi_i'(u),x_j^+(v)] &= \lp\frac{1}{u-v+a} - \frac{1}{u-v-a}\rp
x_j^+(v) \nonumber\\
&+\frac{1}{u-v-a}x_j^+(u-a) - \frac{1}{u-v+a}x_j^+(u+a) \label{eq: imp-lem}
\end{align}

Thus, if $x_j^+(v)$ has a pole at $u_0$ of order $N$, then multiplying
both sides by $(v-u_0)^N$ and letting $v\to u_0$ we get:
\[
[\xi_i(u)^{-1}\xi_i'(u),X] = \lp\frac{1}{u-u_0+a} - \frac{1}{u-u_0-a}\rp X
\]
where $X = \left. (v-u_0)^Nx_j^+(v)\right|_{v=u_0}$. Hence the logarithmic
derivative of $\xi_i(u)$ has poles at $u_0\pm a$, which implies that 
$u_0\pm a$ must be poles of $\xi_i(u)^{\pm 1}$. The argument for
$x_j^-(v)$ is same as above, upon replacing $a$ by $-a$.
\end{pf}

\section{The Drinfeld coproduct}\label{sec: dt}

In this section, we review the definition of the deformed Drinfeld
coproduct on $\qloop$ following \cite{hernandez-affinizations,
hernandez-drinfeld-coproduct}. We then express it in terms of
contour integrals, and use these to determine the poles of the
coproduct as a function of the deformation parameter. By
degenerating the integrals, we obtain a deformed Drinfeld
coproduct for the Yangian $\Yhg$. We also point out that these
coproducts define a meromorphic tensor product on the category
of \fd representations of $\qloop$ and $\Yhg$.

\subsection{Drinfeld coproduct on $\qloop$}\label{ssec: dt-qla1}

Let $\V,\W\in\Rloop$. Twisting Drinfeld's coproduct on $\qloop$
by the $\C^\times$--action on the first factor yields an action of
$\qloop$ on $\V((\zeta^{-1}))\otimes\W$, where $\zeta$ is a formal
variable \cite{hernandez-affinizations,hernandez-drinfeld-coproduct}.
This action is given on the generators of $\qloop$ by\footnote
{We use a different convention than \cite{hernandez-affinizations,
hernandez-drinfeld-coproduct}. The coproduct $\Delta^{(H)}_
{\zeta}$ given in \cite{hernandez-affinizations,hernandez-drinfeld-coproduct}
yields an action on $\V\otimes\W((\zeta))$ obtained by twisting
the Drinfeld coproduct by the $\C^\times$--action on the second
tensor factor. The above action is equal to $\Delta^{(H)}_{\zeta^
{-1}}(\tau_{\zeta}(X))$.}
\begin{align*}
\Psi^{\pm}_{i,\pm m}&\longrightarrow
\sum_{p_1+p_2=m} \zeta^{\pm p_1}\Psi^\pm_{i,\pm p_1}\otimes\Psi^\pm_{i,\pm p_2}\\
\X_{i,k}^+&\longrightarrow
\zeta^k \X^+_{i,k}\otimes 1 + 
\sum_{l\geq 0} \zeta^{-l}\Psi^-_{i,-l}\otimes \X^+_{i,k+l}\\
\X^-_{i,k}&\longrightarrow
\sum_{l\geq 0} \zeta^{k-l}\X^-_{i,k-l}\otimes \Psi^+_{i,l}+1\otimes \X^-_{i,k}
\end{align*}
Hernandez proved that the above formulae are the Laurent expansions
at $\zeta=\infty$ of a family of actions of $\qloop$ on $\V\otimes\W$ the
matrix coefficients of which are rational functions of $\zeta$ \cite[Lemma 3.10]
{hernandez-drinfeld-coproduct}.

\subsection{}\label{ssec: dt-qla2}

Let $\V,\W\in\Rloop$ be as above, and $\spec(\V),\spec(\W)\subset
\C^\times$ be their sets of poles (see \ref{ss:poles}). Let $\zeta\in\C^
\times$ be such that $\zeta\sigma(\V)$ and $\sigma(\W)$ are disjoint,
and define an action of the generators of $\qloop$ on $\V\otimes\W$
as follows
\begin{align*}
\Delta_{\zeta}(\Psi^\pm_{i,\pm m}) &=
\sum_{p_1+p_2=m} \zeta^{\pm p_1}\Psi^\pm_{i,\pm p_1}\otimes\Psi^\pm_{i,\pm p_2}\\
\Delta_{\zeta}(\X^+_{i,k}) &= 
\zeta^k \X^+_{i,k}\otimes 1 + \oint_{C_2} 
\Psi_i(\zeta^{-1}w)\otimes \X^+_i(w) w^{k-1}dw\\
\Delta_{\zeta}(\X^-_{i,k}) &= 
\oint_{C_1} 
\X^-_i(\zeta^{-1}w)\otimes \Psi_i(w)w^{k-1}dw
+ 1\otimes\X^-_{i,k}
\end{align*}
where
\begin{itemize}
\item $C_1,C_2\subset\C^\times$ are Jordan curves which do not
enclose $0$.
\item $C_1$ encloses $\zeta\spec(\V)$ and none of the points in $\spec(\W)$.
\item $C_2$ encloses $\spec(\W)$ and none of the points in $\zeta\spec(\V)$.
\end{itemize}

The above operators are holomorphic functions of $\zeta\in\C^\times
\setminus\spec(\W)\spec(\V)^{-1}$. The corresponding generating series $\Delta_\zeta(\Psi_i(z)^\pm),
\Delta_\zeta(\X_i^\veps(z)^\pm)$ are the expansions at $z=\infty,0$
of the $\End(\V\otimes\W)$--valued holomorphic functions
\begin{align*}
\Delta_{\zeta}(\Psi_i(z)) &= \Psi_i(\zeta^{-1}z)\otimes \Psi_i(z)\\
\Delta_{\zeta}(\X^+_i(z)) &= 
\X^+_i(\zeta^{-1}z)\otimes 1 + \oint_{C_2} \frac{zw^{-1}}{z-w} 
\Psi_i(\zeta^{-1}w)
\otimes \X^+_i(w)\, dw\\
\Delta_{\zeta}(\X^-_i(z)) &=
\oint_{C_1} \frac{zw^{-1}}{z-w} 
\X^-_i(\zeta^{-1}w)\otimes \Psi_i(w)\, dw
+ 1\otimes \X^-_i(z)
\end{align*}
where the integrals are understood to mean the function of
$z$ defined for $z$ outside of $C_1,C_2$. 
Throughout this paper, inside/outside of a Jordan curve $C$ refers to
the bounded/unbounded components of the complement $\C\setminus C$,
and thus they exclude $C$ itself.
We shall prove
below that their dependence in both $\zeta$ and $z$ is rational.

\subsection{}

\begin{thm}\label{th:D U}\hfill
\begin{enumerate}
\item The Laurent expansion of $\Delta_\zeta$ at $\zeta=\infty$ is
given by the deformed Drinfeld coproduct of Section \ref{ssec: dt-qla1}.
\item $\Delta_\zeta$ defines an action of $\qloop$ on $\V\otimes \W$.
The resulting representation is denoted by $\V\Dotimes{\zeta}\W$.
\item The action of $\qloop$ on $\V\Dotimes{\zeta}\W$ is a rational
function of $\zeta$, with poles contained in $\spec(\W)\spec(\V)^{-1}$.
\item The identification of vector spaces
\[\lp\V_1\Dotimes{\zeta_1}\V_2\rp\Dotimes{\zeta_2}\V_3 = 
\V_1\Dotimes{\zeta_1\zeta_2}\lp\V_2\Dotimes{\zeta_2}\V_3\rp\]
intertwines the action of $\qloop$.
\item If $\V\cong\C$ is the trivial representation of $\qloop$, then
\[\V\Dotimes{\zeta}\W=\W\aand\W\Dotimes{\zeta}\V=\W(\zeta)\]
\item The following holds for any $\zeta,\zeta'\in\C^\times$
\[\V\Dotimes{\zeta\zeta'}\W=\V(\zeta)\Dotimes{\zeta'}\W
\aand
V(\zeta')\Dotimes{\zeta}\W(\zeta')=(\V\Dotimes{\zeta}\W)(\zeta')
\]
\item The following holds for any $\zeta\in\C^\times$
\[\spec(\V\Dotimes{\zeta}\W)\subseteq (\zeta\spec(\V))\cup\spec(\W)\]
\end{enumerate}
\end{thm}
\begin{pf}
(i) Expanding $\Delta_\zeta(\Psi^\pm_{i,m})$ and $\Delta_\zeta(\X^
\pm_{i,k})$ as Laurent series in $\zeta^{-1}$ yields the following for
any $m\in\N$ and $k\in\Z$
\begin{align*}
\DD{\zeta}(\Psi^{\pm}_{i,\pm m}) &= \sum_{n=0}^m \zeta^{\pm n}
\Psi^{\pm}_{i,\pm n}\otimes \Psi^{\pm}_{\pm (m-n)}\\
\DD{\zeta}(\X_{i,k}^+) 
&= 
\zeta^k \X^+_{i,k}\otimes 1 + 
\sum_{l\geq 0}\zeta^{-l} \oint_{C_2}\Psi_{i,-l}^-\otimes \X^+_i(w) w^{k+l-1}dw\\
&=
\zeta^k \X^+_{i,k}\otimes 1 + 
\sum_{l\geq 0} \zeta^{-l}\Psi^-_{i,-l}\otimes \X^+_{i,k+l} \\
\DD{\zeta}(\X^-_{i,k})
&=
\oint_{\zeta^{-1}C_1} 
\X^-_i(w)\otimes \Psi_i(\zeta w) \zeta^k w^{k-1}dw
+ 1\otimes\X^-_{i,k}\\
&=
\sum_{l\geq 0}\zeta^{k-l}\oint_{\zeta^{-1}C_1} 
\X^-_i(w)\otimes \Psi_{i,l}^+ w^{k-l-1}dw+1\otimes\X^-_{i,k}\\
&=
\sum_{l\geq 0} \zeta^{k-l}\X^-_{i,k-l}\otimes \Psi^+_{i,l}+1\otimes \X^-_{i,k}
\end{align*}
where the third and sixth equalities follow by Corollary \ref{co:modes}, and
the fourth by a change of variables. Note that $C_1$ is assumed to enclose
$\zeta\sigma(\V_1)$, thus $\zeta^{-1}C_1$ in the computation of $\Delta_{\zeta}
(\X^-_{i,k})$ above encloses $\sigma(\V_1)$.

(ii) By Remark \ref{rm:QL6}, it suffices to check the relations (QL1)--(QL5).
These follow from (i) and \cite[Prop. 6.3]{hernandez-affinizations}, since
it is sufficient to prove them when $\zeta$ is a formal variable. Alternatively,
a direct proof can be given along the lines of Theorem \ref{th:D Y} below.

(iii) The rationality of $\V\Dotimes{\zeta}\W$ follows from (i) and \cite
[Lemma 3.10]{hernandez-drinfeld-coproduct}. Alternatively, let $\{w_j\}_{j\in J}
\subset\C^\times$ be the poles of $\X_i^+(w)$ on $\W$, and
\[\X_i^+(w)=\X_{i,0}^++\sum_{j\in J,n\geq 1}\X^+_{i;j,n}(w-w_j)^{-n}\]
its corresponding partial fraction decomposition. Since $C_2$ encloses
all $w_j$, and $\Psi_i(\zeta^{-1} w)w^{k-1}$ is regular inside $C_2$,
we get 
\[\Delta_{\zeta}(\X^+_{i,k})=
\zeta^k \X^+_{i,k}\otimes 1+
\sum_{j,n}\partial_w^{(n-1)}\left.\left(\Psi_i(\zeta^{-1}w)w^{k-1}\right)\right|_{w=w_j}
\otimes\X^+_{i;j,n}
\]
where $\partial^{(p)}=\partial^p/p!$. This is clearly a rational function of
$\zeta$, whose poles are a subset of the points $\zeta=w_j{w'_k}^{-1}$,
where $w'_k$ is a pole of $\Psi_i(w)$ on $\V$. A similar argument shows
that $\Delta_{\zeta} (\X^-_{i,k})$ is also a rational function whose poles
are contained in $\spec(\W)\spec(\V)^{-1}$. 

(iv) Follows from (i) and \cite[Lemma 3.4]{hernandez-drinfeld-coproduct}.

(v), (vi) and (vii) are clear.
\end{pf}

\subsection{Degeneration}\label{ssec: deg}

The formulae for the Drinfeld coproduct on $\Yhg$ given in \ref
{ssec: dt-yangian} below can be formally obtained by degenerating
those for the Drinfeld coproduct of $\qloop$ given in \ref{ssec: dt-qla2}.
This amounts to setting $z=e^{2\pi\iota\epsilon u}$, $w=e^{2\pi\iota
\epsilon v}$, and letting $\epsilon\to 0$. Under this limit, the $1$--form
$\frac{zw^{-1}}{z-w}dw$ goes to $\frac{dv}{u-v}$.
In addition, we replace the trigonometric functions $\Psi_i(z),
\X_i^\pm(z)$ coming from $\qloop$ by their rational counterparts
$\xi_i(u), x_i^{\pm}(u)$. This method is solely a heuristic, and
a proof that the formulae given in \ref {ssec: dt-yangian} satisfy
the relations of the Yangian $\Yhg$ is provided in \ref{ssec: pf-dt-y23}--\ref{ssec: pf-dt-y-coass}.

\subsection{Drinfeld coproduct on $\Yhg$}\label{ssec: dt-yangian}

Let now $V,W\in\Ryang$, and $\spec(V),\spec(W)\subset\C$
be their sets of poles. Let $s\in\C$ be such that $\spec(V)+s$ and
$\spec(W)$ are disjoint, and define an action of the generators
of $\Yhg$ on $V\otimes W$ via
\begin{align*}
\DD{s}(\xi_i(u)) &= \xi_i(u-s)\otimes\xi_i(u)\\
\DD{s}(x_i^+(u)) &= x_i^+(u-s)\otimes 1 + \oint_{C_2} \frac{1}{u-v}
\xi_i(v-s)\otimes x_i^+(v)\, dv\\
\DD{s}(x_i^-(u)) &= \oint_{C_1} \frac{1}{u-v} x_i^-(v-s)\otimes \xi_i(v)\, dv
+ 1\otimes x_i^-(u)
\end{align*}
where
\begin{itemize}
\item $C_2$ encloses $\spec(W)$ and none of the points in $\spec(V)+s$.
\item $C_1$ encloses $\spec(V)+s$ and none of the points in $\spec(W)$.
\item The integrals are understood to mean the holomorphic functions of
$u$ they define in the domain where $u$ is outside of $C_1,C_2$.\\
\end{itemize}

In terms of the generators $\{\xi_{i,r},x_{i,r}^\pm\}$, the above formulae
read
\begin{align*}
\DD{s}(\xi_{i,r}) &=
\tau_s(\xi_{i,r})\otimes 1+
\hbar\sum_{p_1+p_2=r-1}\tau_s(\xi_{i,p_1})\otimes\xi_{i,p_2}+1\otimes\xi_{i,r}\\
\DD{s}(x_{i,r}^+) &=
\tau_s(x^+_{i,r})\otimes 1
+
\hbar^{-1}\oint_{C_2}
\xi_i(v-s)\otimes x_i^+(v)v^r dv\\
\DD{s}(x_{i,r}^-) &= \hbar^{-1}\oint_{C_1} x_i^-(v-s)\otimes \xi_i(v)v^r dv
+ 1\otimes x_{i,r}^-
\end{align*}

\subsection{}

\begin{thm}\label{th:D Y}\hfill
\begin{enumerate}
\item The formulae in \ref{ssec: dt-yangian} define an action of $\Yhg$ on
$V\otimes W$. The resulting representation is denoted by $V\Dotimes{s}W$. 
\item The action of $\Yhg$ on $V\Dotimes{s} W$ is a rational function of
$s$, with poles contained in $\spec(W)-\spec(V)$.
\item The identification of vector spaces
\[\lp V_1\Dotimes{s_1}V_2\rp\Dotimes{s_2}V_3 = 
V_1\Dotimes{s_1+s_2}\lp V_2\Dotimes{s_2}V_3\rp\]
intertwines the action of $\Yhg$.
\item If $V\cong\C$ is the trivial representation of $\Yhg$, then
\[V\Dotimes{s}W=W\aand W\Dotimes{s}V=W(s)\]
\item The following holds for any $s,s'\in\C$,
\[V\Dotimes{s+s'}W=V(s)\Dotimes{s'}W \aand
V(s')\Dotimes{s}W(s')=(V\Dotimes{s}W)(s')\] 
\item The following holds for any $s\in\C$,
\[\sigma(V\Dotimes{s}W)\subset (\spec(V)+s)\cup\spec(W)\]
\end{enumerate}
\end{thm}

\begin{pf}
(ii) is proved as in Theorem \ref{th:D U}, and (iv)--(vi) are clear.

To prove (i), it suffices by Remark \ref{rm:Y6} to check that
relations (Y1)--(Y5) hold on $V\Dotimes{s}W$. By (v), we
may assume that $\spec(V)\cap\spec(W)=\emptyset$, and
that $s=0$. We choose the contours $C_1$ and $C_2$
enclosing $\spec(V)$ and $\spec(W)$ respectively, such
that they do not intersect. The relation (Y1) holds trivially.
The relations (Y2) and (Y3) are checked in \ref{ssec: pf-dt-y23},
(Y4) in \ref{ssec: pf-dt-y4} and (Y5) in \ref{ssec: pf-dt-y5}.

(iii) is proved in \ref{ssec: pf-dt-y-coass}.
\end{pf}

\subsection{Proof of (Y2) and (Y3)}\label{ssec: pf-dt-y23}

We prove these relations for the $+$ case only. By Proposition
\ref{pr:Y fields} and Remark \ref{rk:Y fields}, it is equivalent to show that $\Delta_0$ preserves
the relation
\[\xi_i(u_1)x_j^+(u_2)\xi_i(u_1)^{-1}
=
\frac{u_1-u_2+a}{u_1-u_2-a}
x_j^+(u_2)-\frac{2a}{u_1-u_2-a} x_j^+(u_1-a)\]
where $a=\hbar d_ia_{ij}/2$. It suffices to prove this for $u_1,
u_2$ large enough, and we shall assume that $u_2$ lies
outside of $C_2$, and that $u_1$ lies outside of $C_2+a$.

Applying $\Delta_0$ to the \lhs gives
\begin{multline*}
\xi_i(u_1)x_j^+(u_2)\xi_i(u_1)^{-1}\otimes 1
+ \oint_{C_2}\frac{1}{u_2-v}\xi_i(v)\otimes 
\xi_i(u_1)x_j^+(v)\xi_i(u_1)^{-1}\, dv\\
= \xi_i(u_1)x_j^+(u_2)\xi_i(u_1)^{-1}\otimes 1
+\oint_{C_2}\frac{u_1-v+a}{(u_2-v)(u_1-v-a)}\xi_i(v)\otimes 
x_j^+(v)\, dv\\
-\oint_{C_2}\frac{2a}{(u_2-v)(u_1-v-a)}\xi_i(v)\otimes 
x_j^+(u_1-a)\, dv
\end{multline*}
where the third summand is equal to zero since the integrand
is regular inside $C_2$.

Applying now $\Delta_0$ to the \rhs yields
\begin{multline*}
\xi_i(u_1)x_j^+(u_2)\xi_i(u_1)^{-1}\otimes 1\\
+\frac{1}{u_1-u_2-a}
\oint_{C_2}\lp\frac{u_1-u_2+a}{u_2-v}
 - \frac{2a}{u_1-a-v}\rp
\xi_i(v)\otimes 
x_j^+(v)\, dv
\end{multline*}
The equality of the two expressions now follows from the identity
\[\frac{u_1-u_2+a}{u_2-v}
 - \frac{2a}{u_1-a-v}=
\frac{(u_1-u_2-a)(u_1+a-v)}{(u_2-v)(u_1-a-v)}\]

\subsection{Proof of (Y4)}\label{ssec: pf-dt-y4}

We check this relation for the $+$ case only. We need to prove
that $\Delta_0$ preserves the relation
\begin{equation}\label{eq: y4-eqform}
x^+_{i,r+1}x^+_{j,s} - x^+_{i,r}x^+_{j,s+1} - ax^+_{i,r}x^+_{j,s}
 = x^+_{j,s}x^+_{i,r+1} - x^+_{j,s+1}x^+_{i,r} + ax^+_{j,s}x^+_{i,r}
\end{equation}
where $a=\hbar d_ia_{ij}/2$. Note that $\DD{0}(x^+_{i,m}x^+_{j,n})$ is
equal to
\begin{multline*}
x^+_{i,m}x^+_{j,n}\otimes 1
+\frac{1}{\hbar} \oint_{C_2}v^nx^+_{i,m}\xi_j(v)\otimes x^+_j(v)\, dv
+\frac{1}{\hbar} \oint_{C_2} v^m\xi_i(v)x^+_{j,n}\otimes x^+_i(v)\, dv\\
+\frac{1}{\hbar^2} \varoiint_{C_2} v_1^mv_2^n\xi_i(v_1)\xi_j(v_2)\otimes
x^+_i(v_1)x^+_j(v_2)\, dv_1dv_2
\end{multline*}
We now apply $\DD{0}$ to both sides of relation \eqref{eq: y4-eqform}, and consider
the four summands of $\DD{0}(x^+_{i,m}x^+_{j,n})$ separately.

{\em The first summand} of $\DD{0}$ of the left and right--hand
sides of \eqref{eq: y4-eqform} are, respectively
\begin{gather*}
\lp x^+_{i,r+1}x^+_{j,s} - x^+_{i,r}x^+_{j,s+1} - ax^+_{i,r}x^+_{j,s}\rp\otimes 1 \\
\lp x^+_{j,s}x^+_{i,r+1} - x^+_{j,s+1}x^+_{i,r} + ax^+_{j,s}x^+_{i,r}\rp\otimes 1
\end{gather*}
which cancel because of \eqref{eq: y4-eqform}.

{\em The second summand} of the \lhs and the {\em third summand}
of the \rhs are, respectively
\begin{gather*}
\frac{1}{\hbar}\oint_{C_2} v^s(x^+_{i,r+1} - vx^+_{i,r} - ax^+_{i,r})\xi_j(v)\otimes x^+_j(v)\, dv\\
\frac{1}{\hbar}\oint_{C_2} v^s\xi_j(v)(x^+_{i,r+1} - vx^+_{i,r} + ax^+_{i,r})\otimes x^+_j(v)\, dv
\end{gather*}
which cancel because of the following version of (Y2) and (Y3)
\[(x^+_{i,r+1} - vx^+_{i,r} - ax^+_{i,r})\xi_j(v)
 = \xi_j(v)(x^+_{i,r+1} - vx^+_{i,r}+ax^+_{i,r})\]
Similarly the third summand of the \lhs and the second summand
of the \rhs cancel.

\noindent {\em The fourth summands} of the left and right--hand
sides of \eqref{eq: y4-eqform} are, respectively
\begin{gather*}
\frac{1}{\hbar^2}\varoiint_{C_2}v_1^rv_2^s(v_1-v_2-a)
\xi_i(v_1)\xi_j(v_2)\otimes x^+_i(v_1)x^+_j(v_2)\, dv_1dv_2 \\
\frac{1}{\hbar^2}\varoiint_{C_2}v_1^rv_2^s(v_1-v_2+a)
\xi_j(v_2)\xi_i(v_1)\otimes x^+_j(v_2)x^+_i(v_1)\, dv_1dv_2
\end{gather*}
By ($\Y$4), their difference is equal to
\[\frac{1}{\hbar}\varoiint_{C_2}
v_1^rv_2^s\,\xi_i(v_1)\xi_j(v_2)\otimes\lp[x^+_{i,0},x^+_j(v_2)]-[x^+_i(v_1),x^+_{j,0}]\rp\, dv_1dv_2\]
which is equal to zero because the first (resp. second) summand
is regular when $v_1$ (resp. $v_2$) lies inside $C_2$.

\subsection{Proof of (Y5)}\label{ssec: pf-dt-y5}

We need to check that $\Delta_0$ preserves the relation
\[[x_i^+(u_1), x_j^-(u_2)] =
-\hbar\delta_{ij}
\frac{\xi_i(u_1)-\xi_i(u_2)}{u_1-u_2}\]
As in Section \ref{ssec: pf-dt-y23} above, it suffices to prove this for $u_1,
u_2$ large enough, and we shall assume that $u_1, u_2$ lie
outside of $C_1,C_2$ respectively.
Applying $\DD{0}$ to the \lhs yields
\begin{multline*}
\oint_{C_1}\frac{1}{u_2-v}
[x_i^+(u_1),x_j^-(v)]\otimes \xi_j(v)\, dv\\
+\oint_{C_2}\frac{1}{u_1-v} \xi_i(v)\otimes
[x^+_i(v),x_j^-(u_2)]\, dv+\calB
\end{multline*}
where
\[\calB = \oint_{C_1}\oint_{C_2} 
\frac{1}{(u_1-v_2)(u_2-v_1)}
[\xi_i(v_2)\otimes x_i^+(v_2),x_j^-(v_1)\otimes\xi_j(v_1)]
dv_2dv_1\]

We shall prove below that $\calB=0$. Thus, by relation (Y5)
the above is equal to zero if $i\neq j$. If $i=j$, it is equal to
\begin{align*}
&-\oint_{C_1}\frac{\hbar}{(u_2-v)(u_1-v)}
\lp\xi_i(u_1)-\xi_i(v)\rp\otimes \xi_i(v)\, dv\\
&
-\oint_{C_2}\frac{\hbar}{(u_1-v)(v-u_2)}
\xi_i(v)\otimes\lp\xi_i(v)-\xi_i(u_2)\rp \, dv \\
&= \oint_{C_1\sqcup C_2} \frac{\hbar}{(u_1-v)(u_2-v)}\xi_i(v)\otimes
\xi_i(v)\, dv\\
&= \frac{\hbar}{u_1-u_2}\lp \xi_i(u_2)\otimes\xi_i(u_2)
 - \xi_i(u_1)\otimes\xi_i(u_1)\rp
\end{align*}
where the first equality follows because $\xi_i(u_1)\otimes\xi_i(v)$
(resp. $\xi_i(v)\otimes\xi_i(u_2)$) is regular when $v$ is inside $C
_1$ (resp. $C_2$), and the second by deformation of contours and
the fact that $\xi_i(v)\otimes\xi_i(v)$ is regular outside $C_1\sqcup
C_2$.

{\em Proof that $\calB=0$}. We shall need the following variant of
relation ($\mathcal{Y}3$) of Proposition \ref{pr:Y fields}.
\begin{equation}\label{eq: dt-pf-y3}
(u-v)[\xi_i(u),x_j^{\pm}(v)] = \pm a\{\xi_i(u), x_j^{\pm}(v)-x_j^{\pm}(u)\}
\end{equation}
where $a = \hbar d_ia_{ij}/2$ and $\{x,y\} = xy+yx$. The integrand
of $\calB$ can be simplified in two different ways. First we write
\begin{multline*}
[\xi_i(v_2)\otimes x_i^+(v_2),x_j^-(v_1)\otimes\xi_j(v_1)]\\
=[\xi_i(v_2),x_j^-(v_1)]\otimes x_i^+(v_2)\xi_j(v_1) 
+x_j^-(v_1)\xi_i(v_2)\otimes [x_i^+(v_2),\xi_j(v_1)]
\end{multline*}
Using \eqref{eq: dt-pf-y3}, we get
\[\begin{split}
\calB
=\oint_{C_1}\oint_{C_2}
&
\frac{a}{(u_1-v_2)(u_2-v_1)(v_1-v_2)}\\
&\lp\{\xi_i(v_2),x_j^-(v_1)-x^-_j(v_2)\}\otimes x_i^+(v_2)\xi_j(v_1)\right.\\
&\left.- x_j^-(v_1)\xi_i(v_2)\otimes\{\xi_j(v_1),x_i^+(v_2)-x_i^+(v_1)\}\rp\,
dv_2dv_1\\
=\oint_{C_1}\oint_{C_2}
&
\frac{a}{(u_1-v_2)(u_2-v_1)(v_1-v_2)}
\lp\{\xi_i(v_2),x_j^-(v_1)\}\otimes x_i^+(v_2)\xi_j(v_1)\right.\\
&\left.- x_j^-(v_1)\xi_i(v_2)\otimes\{\xi_j(v_1),x_i^+(v_2)\}\rp\,
dv_2dv_1\\
=\oint_{C_1}\oint_{C_2}
&
\frac{a}{(u_1-v_2)(u_2-v_1)(v_1-v_2)}
\lp\xi_i(v_2)x_j^-(v_1)\otimes x_i^+(v_2)\xi_j(v_1)\right.\\
&\left.- x_j^-(v_1)\xi_i(v_2)\otimes\xi_j(v_1)x_i^+(v_2)\rp\,
dv_2dv_1\\
\end{split}\]
where the second equality follows from the fact that $\{\xi_i(v_2),x^-_j(v_2)\}
\otimes x_i^+(v_2)\xi_j(v_1)$ (resp. $x_j^-(v_1)\xi_i(v_2)\otimes\{\xi_j(v_1),
x_i^+(v_1)\}$) is regular when $v_1$ is inside $C_1$ (resp. $v_2$ is inside
$C_2$).

Now if we write instead
\begin{multline*}
[\xi_i(v_2)\otimes x_i^+(v_2),x_j^-(v_1)\otimes\xi_j(v_1)]\\
=
\xi_i(v_2)x_j^-(v_1)\otimes [x_i^+(v_2),\xi_j(v_1)]
+ [\xi_i(v_2),x_j^-(v_1)]\otimes \xi_j(v_1)x_i^+(v_2)
\end{multline*}
and use relation \eqref{eq: dt-pf-y3} as before, we obtain
\begin{multline*}
\calB = \oint_{C_1}\oint_{C_2} 
\frac{-a}{(v_1-v_2)(u_1-v_2)(u_2-v_1)}\lp
\xi_i(v_2)x_j^-(v_1)\otimes x_i^+(v_2)\xi_j(v_1)\right.\\
 - \left.x_j^-(v_1)\xi_i(v_2)\otimes \xi_j(v_1)x_i^+(v_2)\rp\,
dv_2dv_1
\end{multline*}
Thus $\calB = -\calB$, whence $\calB=0$.

\subsection{Coassociativity}\label{ssec: pf-dt-y-coass}

We need to show that the generators of $\Yhg$ act
by the same operators on
\[\lp V_1\Dotimes{s_1}V_2\rp\Dotimes{s_2}V_3\aand
V_1\Dotimes{s_1+s_2}\lp V_2\Dotimes{s_2}V_3\rp\]
The action of $\xi_i(u)$ on both modules is given by
$\xi_i(u-s_1-s_2)\otimes \xi_i(u-s_2)\otimes \xi_i(u)$.

To compute the action of $x_i^+(u)$, we shall assume
that $s_1$ and $s_2$ are such that $\spec(V_1)+s_1+
s_2,\spec(V_2)+s_2$ and $\spec(V_3)$ are all disjoint.
By (vi), this implies in particular that $\spec(V_1\Dotimes
{s_1}V_2)+s_2$ and $\spec(V_3)$ are disjoint, and that
so are $\spec(V_1)+s_1+s_2$ and $\spec(V_2\Dotimes
{s_2}V_3)$, so that the above tensor products are defined.

Under these assumptions, the action of $x_i^+(u)$ on
$(V_1\Dotimes{s_1}V_2)\Dotimes{s_2}V_3$ is given by
\[\begin{split}
&\Delta_{s_1}(x_i^+(u-s_2))\otimes 1+
\oint_{C_3}\frac{1}{u-v_3}\Delta_{s_1}(\xi_i(v_3-s_2))\otimes x_i^+(v_3)\,dv_3\\
&=
x_i^+(u-s_2-s_1)\otimes 1\otimes 1+
\oint_{C_2}\frac{1}{u-s_2-v_2}\xi_i(v_2-s_1)\otimes x_i^+(v_2)\otimes 1\,dv_2
\\
&\phantom{=}+\oint_{C_3}\frac{1}{u-v_3}\xi_i(v_3-s_2-s_1)\otimes\xi_i(v_3-s_2)\otimes x_i^+(v_3)\,dv_3
\end{split}\]
where $C_3$ encloses $\spec(V_3)$ and none of the points
of $\spec(V_1)+s_1+s_2$ and $\spec(V_2)+s_2$, $C_2$
encloses $\spec(V_2)$ and none of the points of $\spec(V_1)
+s_1$, and $u$ is assumed to be outside of $C_3$ and $C_2
+s_2$.

The action of $x_i^+(u)$ on $V_1\Dotimes{s_1+s_2}\lp V_2
\Dotimes{s_2}V_3\rp$ is given by
\[\begin{split}
&x_i^+(u-s_1-s_2)\otimes 1\otimes 1+
\oint_{C_{23}}\frac{1}{u-v_{23}}\xi_i(v_{23}-s_1-s_2)\otimes \Delta_{s_2}(x_i^+(v_{23}))\,dv_{23}\\
&=
x_i^+(u-s_1-s_2)\otimes 1\otimes 1\\
&\phantom{=}+
\oint_{C_{23}}\frac{1}{u-v_{23}}\xi_i(v_{23}-s_1-s_2)\otimes x_i^+(v_{23}-s_2)\otimes 1 \,dv_{23}\\
&\phantom{=}+
\oint_{C_{23}}\oint_{C_3^\prime}\frac{1}{u-v_{23}}\frac{1}{v_{23}-v_3^\prime}
\xi_i(v_{23}-s_1-s_2)\otimes\xi(v_3^\prime-s_2)\otimes x_i^+(v_3^\prime)\,dv_3^\prime dv_{23}
\end{split}\]
where $C_{23}$ encloses $\spec(V_2)+s_2\cup\spec(V_3)$
and none of the points of $\spec(V_1)+s_1+s_2$, $C_3^\prime$
is chosen inside $C_{23}$ and encloses $\spec(V_3)$ and none
of the points of $\spec(V_2)+s_2$, and $u$ is assumed to be
outside of $C_{23}$.

Since the singularites of the first integrand which are contained
in $C_{23}$ lie in $\spec(V_2)+s_2$, the corresponding integral
is equal to
\[\oint_{C_2^\prime}\frac{1}{u-v_2^\prime}\xi_i(v_2^\prime-s_1-s_2)\otimes x_i^+(v_2^\prime-s_2)\otimes 1 \,dv_2^\prime\]
where $C_2^\prime$ contains $\spec(V_2)+s_2$ and none of
the points of $\spec(V_1)+s_1+s_2$. On the other hand, 
integrating in $v_{23}$ in the second integral yields
\[\oint_{C_3^\prime}\frac{1}{u-v_3^\prime}
\xi_i(v_3^\prime-s_1-s_2)\otimes\xi(v_3^\prime-s_2)\otimes x_i^+(v_3^\prime)\,dv_3^\prime\]
so that the two actions of $x_i^+(u)$ agree. The proof for $x_i^-
(u)$ is similar.

\section{The commutative $R$-matrix of the Yangian}\label{sec: R0}

In this section, we construct the commutative part $\calR^0(s)$
of the $R$--matrix of the Yangian, and show that it defines 
meromorphic commutativity constraints on $\Ryang$, when the
latter is equipped with the Drinfeld tensor product defined in \S
\ref{sec: dt}.

A conjectural formula expressing $\calR^0(s)$ as a formal infinite
product with values in the {\it double Yangian} $\D\Yhg$ was given
by Khoroshkin--Tolstoy \cite[Thm. 5.2]{khoroshkin-tolstoy}. We
review their formula in \S \ref{ssec: q-cartan}--\ref{ss:KT}, and 
outline our own construction in \ref{ss:our strategy}.
Our starting point is the observation that $\calR^0(s)$ formally
satisfies an additive difference equation whose coefficient matrix
$\A(s)$ we show to be a rational function on \fd representations
of $\Yhg$. By taking the left and right canonical fundamental
solutions of this equation, we construct two regularisations
$\Rud(s)$ of $\calR^0(s)$ which are meromorphic functions
of the parameter $s$, and then show that they have the required
intertwining properties with respect to the Drinfeld coproduct.

Note that Sections \ref{ss:KT} and \ref{ss:our strategy} are
included solely to motivate our construction, and that the
definition of $\Rud(s)$ and the proofs of its properties
are independent of the results of \cite{khoroshkin-tolstoy}.
In particular, we do not work with the double Yangian. 

\subsection{The $T$--Cartan matrix of $\g$}\label{ssec: q-cartan}

Let $\bfA = (a_{ij})$ be the Cartan matrix of $\g$ and $\bfB=
(b_{ij})$ its symmetrization, where $b_{ij}=d_ia_{ij}$. Let $\T$
be an indeterminate, and let $\bfB(\T) = ([b_{ij}]_\T)\in GL_\bfI
(\C[\T^{\pm 1}])$ the corresponding matrix of $T$--numbers.
Then, there exists an integer $l=mh^\vee$, which is a multiple
of the dual Coxeter number $h^\vee$ of $\g$, and is such that
\begin{equation}\label{eq:B inv}
\bfB(T)^{-1}=\frac{1}{[l]_\T}\bfC(\T)
\end{equation}
where the entries of $\bfC(\T)$ are Laurent polynomials in $\T$
with positive integer coefficients.\footnote{This result is stated
without proof in \cite[p. 391]{khoroshkin-tolstoy}, and proved
for $\g$ simply--laced in \cite[Prop. 2.1]{hernandez-leclerc2}.
We give a proof in Appendix \ref{se:KT}, which also corrects
the values of the multiple $m$ tabulated in \cite{khoroshkin-tolstoy}
for the $\sfC_n$ and $\sfD_n$ series. With those corrections,
the value of $m$ for any $\g$ is the ratio of the squared length
of long roots and short ones.} We denote the entries of the
matrix $\bfC(T)$ by $c_{ij}(\T) = \sum_{r\in\Z}c_{ij}^{(r)} \T^r$,
and note that $c_{ji}(\T)=c_{ij}(\T)=c_{ij}(T^{-1})$.

\subsection{The \KT construction}\label{ss:KT}

The starting point of \cite{khoroshkin-tolstoy} is a conjectural
presentation of the Drinfeld double $\DY$ of the Yangian
$\Yhg$. $\DY$ is generated by
$\{\xi_{i,r},x_{i,r}^{\pm}\}_{i\in\bfI,r\in\Z_{\geq 0}}$
and $\{\xi_{i,r}, x_{i,r}^{\pm}\}_{i\in\bfI,r\in\Z_{<0}}$,
where the first are the generators of $\Yhg$. 
We will not need the complete presentation of $\DY$. For
our purposes, it is sufficient to know that $\DY$
contains the following two sets of commuting elements: 
$\{\xi_{i,r}\}_{i\in\bfI,r\in\Z_{\geq 0}}$ and $\{\xi
_{i,r}\}_{i\in\bfI,r\in\Z_{<0}}$. Let $Y_0^\pm\subset\DY$ be the
subalgebras they generate. The Hopf pairing $\<-,-\>$
on $\DY$ restricts to a perfect pairing $Y_0^+\otimes Y_0^-
\to\C$, and the commutative part of the $R$--matrix of $\Yhg$
is given by
\begin{equation}\label{eq:R0}
\calR^0=\exp\left(\sum_{i\in\bfI,r\in\N} a^+_{i,r}\otimes a^-_{i,-r-1}\right)
\end{equation}
where $\{a^+_{i,r}\}_{i\in\bfI,r\in\Z_{\geq 0}}$ and $\{a^-_{i,r}\}
_{i\in\bfI,r\in\Z_{< 0}}$ are generators of $Y_0^+,Y_0^-$
respectively, which are primitive modulo elements which
pair trivially with $Y_0^\pm$, and such that $\<a^+_{i,r},
a^-_{j,-s-1}\>=\delta_{ij}\delta_{rs}$.

Constructing these generators amounts to finding formal
power series
\[a^+_i(u)=\sum_{r\geq 0} a^+_{i,r}u^{-r-1}\in Y_0^+[[u^{-1}]]
\quad\text{and}\quad
a^-_i(v)=\sum_{r< 0} a^-_{i,r}v^{-r-1}\in Y_0^-[[v]]\]
such that
$\<a^+_i(u),a^-_j(v)\>=\delta_{ij}/(u-v)$. To this end,
introduce the generating series
\[\xi_i^+(u)=1+\hbar\sum_{r\geq 0}\xi_{i,r}u^{-r-1}
\aand
\xi_i^-(v)=1-\hbar\sum_{r<0}\xi_{i,r}v^{-r-1}\]
Then, by definition of $\DY$, we have
\[\<\xi^+_i(u),\xi^-_j(v)\>=\frac{u-v+a}{u-v-a}\in\C[[u^{-1},v]]\]
where $a=\hbar b_{ij}/2$.
Define now
\begin{equation}\label{eq:tipm}
t_i^+(u)=\log(\xi^+_i(u))\in Y_0^+[[u^{-1}]]
\quad\text{and}\quad
t_i^-(v)=\log(\xi^-_i(v))\in Y_0^-[[v]]
\end{equation}
Then, it follows that
\[\<t_i^+(u),t_j^-(v)\>=\log\left(\frac{u-v+a}{u-v-a}\right)\]
Indeed, $\xi_i^\pm(u)$ are group--like modulo terms which pair
trivially with $Y_0^+,Y_0^-$, and if $a,b$ are primitive elements
of a Hopf algebra endowed with a Hopf pairing $\<-,-\>$, then
$\<e^a,e^b\>=e^{\<a,b\>}$. Differentiating \wrt $u$ yields
\[\<\frac{d}{du}t_i^+(u),t_j^-(v)\>=\frac{1}{u-v+a}-\frac{1}{u-v-a}\]

Let $T$ be the shift operator acting on functions of $v$ as $Tf(v)
=f(v-\hbar/2)$. Then, the above identity may be rewritten as
\[\<\frac{d}{du}t_i^+(u),t_j^-(v)\>=
(T^{b_{ij}}-T^{-b_{ij}})\frac{1}{u-v}=
(T-T^{-1})\bfB(T)_{ij}\frac{1}{u-v}\]
where $\bfB(T)$ is the matrix introduced in \ref{ssec: q-cartan}.
It follows that if $\bfD(T)$ is an $\bfI\times\bfI$ matrix with entries
in $\C[[T,T^{-1}]]$, then
\[\sum_k\bfD(T)_{jk}\<\frac{d}{du}t_i^+(u),t_k^-(v)\>=
(T-T^{-1})(\bfD(T)\bfB(T))_{ji}\frac{1}{u-v}\]

By \eqref{eq:B inv}, choosing $\bfD(T)=(T^l-T^{-l})^{-1}\bfC(T)$,
and setting
\begin{equation}\label{eq:a's}
a_i^+(u)=\frac{d}{du}t^+_i(u)\aand
a_j^-(v)=\sum_{k\in\bfI} (T^l-T^{-l})^{-1}\bfC(T)_{jk}t_k^-(v)
\end{equation}
gives the sought for generators, provided one can interpret $(T^l-T^{-l})
^{-1}$. This can be done by expanding in powers of $T^l$ or of
$T^{-l}$, and leads to two distinct formal expressions for $\calR
^0$ \cite[(5.27)--(5.28)]{khoroshkin-tolstoy}.

\subsection{} \label{ss:our strategy}

To make sense of the above construction of $\calR^0$ on the
tensor product $V_1\otimes V_2$ of two \fd representations of
$\Yhg$, we proceed as follows.
\renewcommand {\theenumi}{\arabic{enumi}}
\begin{enumerate}
\item By \ref{prop: rationality}, $a_i^+(u)$ acting on $V_1$:
\[a_i^+(u)=\frac{d}{du} t_i^+(u)=\xi_i^+(u)'\xi_i^+(u)^{-1}\]
is a rational $\End(V_1)$--valued function of $u$, regular near $\infty$.
\item If $a_j^-(v)$ defined by \eqref{eq:a's} can be shown to be
a meromorphic function of $v$, we may interpret the sum over
$r$ in \eqref{eq:R0} as the contour integral $\oint_C a_i^+(u)
\otimes a_i^-(u)\,du$, where $C$ encloses all poles of $a_i^+(u)$
and none of those of $a_i^-(u)$.
\item The action of $\calR^0$ on $V_1(s)\otimes V_2$ would
then be given by
\[\calR^0(s)=
\exp\left(\sum_i\oint_{C+s}a_i^+(u-s)\otimes a_i^-(u)\,du\right)
=
\exp\left(\sum_i\oint_{C}a_i^+(u)\otimes a_i^-(u+s)\,du\right)\]
where $C$ encloses all poles of $a_i^+(u)$ on $V_1$ and
none of those of $a_i^-(u)$ on $V_2$.
\item We show in \ref{ss:log} that, on any \fd representation of
$\Yhg$, $t_i^+(u)$ is the expansion near $u=\infty$ of a meromorphic
function of $u$ defined on the complement
of a compact cut--set $0\in\sfX\subset\C$, and interpret $t_i^-(v)$ as
the corresponding analytic continuation of $t_i^+(u)$.
This resolves
in particular the ambiguity in the definition \eqref{eq:tipm} of $t_i^-(v)$ as a formal
power series in $v$, since the constant term of $\xi_i^-(v)$ is not
equal to $1$, and allows to apply the shift operator $T$ to $t_j^-(v)$,
since $T$ does not act on formal power series of $v$.
Moreover, since we work with $\Yhg$, we do not have the operators
$\{\xi_{i,r}\}_{i\in\bfI,r\in\Z_{<0}}$ at our disposal. This makes
the reinterpretation of $t_i^-(v)$ as a meromorphic function essential
for our purposes.
\item To interpret $a^-_j(v)$, we note that it formally satisfies the
difference equation $a^-_j(v+l\hbar)-a^-_j(v)=b^-_j(v)$, where
\[b_j^-(v)=
-\sum_{k\in\bfI} T^{-l}\bfC(T)_{jk}t^-_k(v)=
-\sum_{k\in\bfI,r\in\Z} c_{jk}^{(r)}t^-_k(v+(l+r)\frac{\hbar}{2})
\]
and we used the fact that $\bfC(T)=\bfC(T^{-1})$. This implies
that $\calR^0(s)$ formally satisfies
\begin{equation}\label{eq:formal ADE}
\calR^0(s+l\hbar)\calR^0(s)^{-1}=
\exp\left(\sum_i\oint_C a_i^+(u)\otimes b^-_i(u+s)\ du\right)
\end{equation}
\item We show in \ref{ssec: op-A}--\ref{ssec: pf-A2} that the
operator ${\mathcal A}(s)$ given by the \rhs of \eqref{eq:formal ADE} is a rational function
of $s$ such that ${\mathcal A}(\infty)=1$. We then define two
regularisations $\Rud(s)$ of $\calR^0(s)$ as the canonical
right and left fundamental solutions of the difference equation
\eqref{eq:formal ADE}, and show in \ref{ssec: R0-main} that
they define meromorphic commutativity constraints on $\Ryang$
endowed with the deformed Drinfeld coproduct.
\end{enumerate}
\renewcommand {\theenumi}{\roman{enumi}}

\subsection{Matrix logarithms}\label{ss:log}

We shall need the following result

\begin{prop}\label{pr:matrix log}
Let $V$ be a complex, \fd vector space, and $\xi:\C\to\End(V)$ a rational
function such that
\begin{itemize}
\item $\xi(\infty)=1$.
\vskip .1cm 
\item $[\xi(u),\xi(v)]=0$ for any $u,v\in\C$.
\end{itemize}
Let $\sigma(\xi)\subset\C$ be the set of poles of $\xi(u)^{\pm 1}$, and define
the cut--set $\sfX(\xi)$ by
\begin{equation}\label{eq:star}
\sfX(\xi)=\bigcup_{a\in\sigma(\xi)}[0,a]
\end{equation}
where $[0,a]$ is the line segment joining $0$ and $a$. Then, there is a
unique single--valued, holomorphic function $t(u)=\log(\xi(u)):\C\setminus
\sfX(\xi)\to\End(V)$ such that
\begin{equation}\label{eq:t}
\exp(t(u))=\xi(u)\aand t(\infty)=0
\end{equation}
Moreover, $[t(u),t(v)]=0$ for any $u,v\in\C$, and $t(u)'=\xi(u)^{-1}\xi'(u)$.
\end{prop}

\begin{pf}
The equation \eqref{eq:t} uniquely defines $t(u)$ as a holomorphic function
near $u=\infty$. 
To continue $t(u)$ meromorphically, note first that the semisimple
and unipotent factors $\xi_S(u),\xi_U(u)$ of the multiplicative Jordan
decomposition of $\xi(u)$ are rational functions of $u$ since $[\xi(u),
\xi(v)]=0$ for any $u,v$ (see \eg \cite[Lemma 4.12]{sachin-valerio-2}). Thus,
\[t_N(u)=\log(\xi_U(u))=\sum_{k\geq 1} (-1)^{k-1}\frac{(\xi_U(u)-1)^k}{k}\]
is a well--defined rational function of $u\in\C$ whose poles are contained in 
those of $\xi(u)$.

To define $\log(\xi_S(u))$ consistently, note that the eigenvalues of
$\xi(u)$ are rational functions of the form $\prod_j(u-a_j)(u-b_j)^{-1}$.
Since, for $a\in\C^\times$, the function $\log(1-au^{-1})$ is single--valued
on the complement of the interval $[0,a]$, where $\log$ is the standard
determination of the logarithm, we may define a single--valued, holomorphic
function $\log(\xi_S(u))$ on the complement of the intervals $[0,a]$,
where $a$ ranges over the (non--zero) zeros and poles of the
eigenvalues of $\xi(u)$.

Finally, we set
\[t(u)=t_N(u)+t_S(u)\]
The fact that $[t(u),t(v)]=0$ is clear from the construction, or from the fact
that it clearly holds for $u,v$ near $\infty$. Finally, the derivative of $t(u)$
can be computed by differentiating the identity $\exp(t(u))=\xi(u)$, and using
the formula for the left--logarithmic derivative of the exponential function
(see, \eg \cite{faraut}).
\end{pf}

\begin{defn}
If $V$ is a \fd representation of $\Yhg$, and $\xi_i(u)$ is the rational function
$\xi_i(u)=1+\hbar\sum_{r\geq 0}\xi_{i,r}u^{-r-1}$ given by Proposition \ref
{prop: rationality}, the corresponding logarithm will be denoted by $t_i(u)$.
\end{defn}

\subsection{The operator $\opA{V_1}{V_2}(s)$}\label{ssec: op-A}

Let $V_1,V_2$ be two \fd representations of $\Yhg$. Let $\calC_1$
be a contour enclosing the set of poles of the operators $\xi_i(u)^{\pm 1}$ on $V_1$,
and consider the following operator on $V_1\otimes V_2$
\[\opA{V_1}{V_2}(s)=
\exp\lp -\sum_{\begin{subarray}{c} i,j\in\bfI\\ r\in\Z\end{subarray}} c_{ij}^{(r)}
\oint_{\calC_1} t_i'(v)\otimes t_j\lp v+s+\frac{(l+r)\hbar}{2}\rp\, dv \rp\]
where $s\in\C$ is such that $t_j(v+s+\hbar(l+r)/2)$ is an analytic function
of $v$ within $\calC_1$ for every $j\in\bfI$ and $r\in\Z$ such that $c_{ij}^
{(r)}\neq 0$ for some $i\in\bfI$.\\

Let $\Omega_{\h}\in\h\otimes\h\subset \Yhg\otimes \Yhg$ be the Cartan part of the Casimir tensor. 
Explicitly, 
\begin{equation}\label{eq: omega-cartan}
\Omega_{\h} = \sum_{i\in\bfI} d_ih_i\otimes \varpi_i^{\vee} = \sum_{i\in\bfI} \varpi_i^{\vee}\otimes d_ih_i
\end{equation}
where $d_ih_i = \xi_{i,0}$, and $\varpi_i^{\vee}$ are the fundamental coweights,
which are defined by $(\varpi_i^{\vee},d_jh_j) = \delta_{ij}$. By definition of the
bilinear form  $(\cdot,\cdot)$ on $\h\times \h$, we have $\varpi_i^{\vee} = \sum_{j\in\bfI}
(\bfB^{-1})_{ij} d_jh_j$.

\begin{thm}\label{thm: coeff}\hfill
\begin{enumerate}
\item $\opA{V_1}{V_2}(s)$ extends to a rational function of $s$ which
is regular at $\infty$, and such that
\[\opA{V_1}{V_2}(s)=1-l\hbar^2\frac{\Omega_\h}{s^2}+O(s^{-3})\]
The poles of $\opA{V_1}{V_2}(s)^{\pm 1}$ are contained in
\[\spec(V_2)-\spec(V_1)-\frac{\hbar}{2}\{l+r\}\]
where $r$ ranges over the integers such that $c_{ij}^{(r)}\neq 0$ for
some $i,j\in\bfI$.
\item For any $s,s'$ we have $[\opA{V_1}{V_2}(s), \opA{V_1}{V_2}(s')]=0$.
\item For any $V_1,V_2,V_3\in\Ryang$, we have
\begin{align*}
\opA{V_1\Dotimes{s_1}V_2}{V_3}(s_2) &= \opA{V_1}{V_3}(s_1+s_2)\opA{V_2}{V_3}
(s_2)\\
\opA{V_1}{V_2\Dotimes{s_2}V_3}(s_1+s_2) &= \opA{V_1}{V_3}(s_1+s_2)
\opA{V_1}{V_2}(s_1)
\end{align*}
\item The following shifted unitarity condition holds
\[\sigma\circ \opA{V_1}{V_2}(-s)\circ \sigma^{-1}=\opA{V_2}{V_1}(s-l\hbar)\]
where $\sigma:V_1\otimes V_2\to V_2\otimes V_1$ is the flip of the tensor
factors.
\item For every $a,b\in\C$ we have
\[
\opA{V_1(a)}{V_2(b)}(s) = \opA{V_1}{V_2}(s+a-b)
\]
\end{enumerate}
\end{thm}

\begin{pf}

Properties (ii),(iii) and (v) follow from the definition of $\A$, and the fact
that $t_i(u)$ are primitive with respect to the Drinfeld coproduct. To prove
(i) and (iv), we work in the following more general situation.

Let $V,W$ be complex, \fd vector spaces, $A,B:\C\to\End(V)$ rational
functions satisfying the assumptions of Proposition \ref{pr:matrix log}, and
let $\log A(v),\log B(v)$ be the corresponding logarithms. Let $\spec(A),
\spec(B)$ denote the set of poles of $A(v)^{\pm 1}$ and $B(v)^{\pm 1}$
respectively. Set
\[X(s) = \exp\lp\oint_{\calC_1} A(v)^{-1}A'(v)\otimes\log(B(v+s))\, dv\rp\]
where $\calC_1$ encloses $\spec(A)$, and $s$ is such that $\log(B(v+s))$
is analytic within $\calC_1$.\\

\noindent
{\bf Claim 1.} The operator $X(s)\in\End(V\otimes W)$ is a rational function
of $s$, regular at $\infty$, and has the following Taylor series expansion
near $\infty$
\[X(s) = 1 + (A_0\otimes B_0)s^{-2} + O(s^{-3})\]
where $A(s)=1+A_0s^{-1}+O(s^{-2})$ and $B(s)=1+B_0s^{-1}+O(s^{-2})$.
Moreover, the poles of $X(s)^{\pm 1}$ are contained in $\spec(B)-\spec(A)$.\\

\noindent
Note that this claim implies the first part of Theorem \ref{thm: coeff} (i), since
\[\begin{split}
\opA{V_1}{V_2}(s)
&=
\prod_{\begin{subarray}{c} i,j\in\bfI\\ r\in\Z\end{subarray}}
\exp\lp \oint_\calC t_i'(v)\otimes t_j\lp v+s+\frac{(l+r)\hbar}{2}\rp\, dv \rp^{-c_{ij}^{(r)}}\\
&=
1-\hbar^{2}s^{-2}\sum_{\begin{subarray}{c} i,j\in\bfI\\ r\in\Z\end{subarray}}
c_{ij}^{(r)}\xi_{i,0}\otimes\xi_{j,0}+O(s^{-3})\\
&= 1 - l\hbar^2 \Omega_{\h} s^{-2} + O(s^{-3})
\end{split}\]
since $\sum_{r\in\Z} c_{ij}^{(r)} = \left.c_{ij}(T)\right|_{T=1}$ is the $(i,j)$ entry of $l\cdot\bfB^{-1}$.\\

Part (iv) of Theorem \ref{thm: coeff} is a consequence of the following
claim, together with the fact that $c_{ji}^{(r)}=c_{ij}^{(r)}=c_{ij}^{(-r)}$.

\noindent
{\bf Claim 2.} $\ds X(s) = \exp\lp\oint_{\calC_2} \log(A(v-s))\otimes 
B(v)^{-1}B'(v)\, dv\rp$, where $\calC_2$ encloses $\spec(B)$ and
$s\in\C$ is such that $\log(A(v-s))$ is analytic within $\calC_2$.\\

We prove these claims in \S \ref{ssec: pf-A1} and \ref{ssec: pf-A2}
respectively.
\end{pf}

\subsection{Proof of Claim 1}\label{ssec: pf-A1}

Since $A(v)$ commutes with itself for different values of $v$,
the semisimple and unipotent parts $A(v) = A_S(v)A_U(v)$
of the Jordan decomposition of $A(v)$ are rational functions
of $v$ \cite[Lemma 4.12]{sachin-valerio-2}. Since the
logarithmic derivative of $A(v)$ separates the two additively,
we can treat the semisimple and unipotent cases separately.

The semisimple case reduces to the scalar case, \ie when $V$
is one--dimensional and
\[A(v)=\prod_j\frac{v-a_j}{v-b_j}=1+(\sum_j b_j-a_j)v^{-1}+O(v^{-2})\]
for some $a_j,b_j\in\C$. In this case,
\[\begin{split}
X(s)
&=
\exp\lp\sum_j\oint_{\calC_1}
\left(\frac{1}{v-a_j}-\frac{1}{v-b_j}\right)\otimes\log(B(v+s))dv\rp\\
&=
\exp\lp\sum_j
1\otimes\left(\log(B(s+a_j))-\log(B(s+b_j))\right)\rp\\
&=
\prod_j 1\otimes B(s+a_j) B(s+b_j)^{-1}
\end{split}\]
which is a rational function of $s$ such that the poles of $X(s)^{\pm 1}$
are contained in $\spec(B)-\spec(A)$. Moreover,
\[X(s) = 1 + s^{-2} \left(\sum_j b_j-a_j\right)\otimes B_0  + O(s^{-3})\]

Assume now that $A(v)$ is unipotent. In this case,
\[\log(A(v)) = \sum_{k\geq 1} (-1)^{k-1} \frac{(A(v)-1)^k}{k}=A_0 v^{-1}+O(v^{-2})\]
is given by a finite sum, and is therefore a rational function of $v$. Decomposing
it into partial fractions yields
\[\log(A(v)) = \sum_{\begin{subarray}{c} j\in J\\ n\in\N\end{subarray}}
\frac{N_{j,n}}{(v-a_j)^{n+1}}\]
where $J$ is a finite indexing set, $a_j\in\C$ and $\sum_j N_{j,0}
=A_0$. In this case we obtain
\[X(s) = \exp\lp\sum_{\begin{subarray}{c} j\in J\\ n\in\N\end{subarray}}
-(n+1) N_{j,n}\otimes \left.\frac{\partial_v^{n+1}}{(n+1)!}\log(B(v)) \right|_{v=s+a_j}
\rp\]
This is again a rational function of $s$ since the $N_{j,n}$ are nilpotent
and pairwise commute, such that the poles of $X(s)^{\pm 1}$ are contained
in $\spec(B)-\spec(A)$. Moreover,
\[X(s)=1+s^{-2}\sum_j N_{j,0}\otimes B_0+O(s^{-3})\]

\subsection{Proof of Claim 2}\label{ssec: pf-A2}

Let $\sfX(A),\sfX(B)\subset\C$ be defined by \eqref{eq:star}, and $\calC_1,
\calC_2$ be two contours enclosing $\sfX(A)$ and $\sfX(B)$ respectively.
For each $s\in\C$ such that $\calC_1+s$ is outside of $\calC_2$, we have
\[\begin{split}
\oint_{\calC_1} A(v)^{-1}A'(v)\otimes \log(&B(v+s))\, dv\\
&= 
-\oint_{\calC_1} \log(A(v))\otimes B(v+s)^{-1}B'(v+s)\, dv\\
&=\phantom{-}
\oint_{\calC_2-s} \log(A(v))\otimes B(v+s)^{-1}B'(v+s)\, dv\\
&=\phantom{-}
\oint_{\calC_2} \log(A(w-s))\otimes B(w)^{-1}B'(w)\, dw\\
\end{split}\]
where the first equality follows by integration by parts, the second by a
deformation of contour since the integrand is regular at $v=\infty$ and
has zero residue there, and the third by the change of variables $w=v+s$.

\subsection{The abelian $R$--matrix of $\Yhg$}\label{ssec: R0-prop}

Let $V_1,V_2\in\Ryang$, and let $\opA{V_1}{V_2}(s)\in GL(V_1\otimes
V_2)$ be the operator defined in Section \ref{ssec: op-A}. Consider the additive
difference equation
\begin{equation}\label{eq:A h ADE}
\calR_{V_1,V_2}(s+l\hbar) = \opA{V_1}{V_2}(s)\calR_{V_1,V_2}(s)
\end{equation}
where $l = mh^{\vee}$ was defined in Section \ref{ssec: q-cartan}.

This equation is regular, in that $\opA{V_1}{V_2}(s)=1+O(s^{-2})$ by
Theorem \ref{thm: coeff}. In particular, it admits two canonical meromorphic
fundamental solutions \[\Rud_{V_1,V_2}:\C\to GL(V_1\otimes V_2)\]
which are uniquely determined by the following requirements (see \eg
\cite{birkhoff-difference,borodin,krichever} or \cite[\S 4]{sachin-valerio-2})
\begin{itemize}
\item $\Rud_{V_1,V_2}(s)$ is holomorphic and invertible for $\pm\Re(s/\hbar)\gg 0$.
\item $\Rud_{V_1,V_2}(s)$ possesses an asymptotic expansion
of the form
\[\Rud_{V_1,V_2}(s)\sim 1+\calR^\pm_0 s^{-1}+\calR^\pm_1s^{-2}+\cdots\]
in any half--plane $\pm\Re(s/\hbar)>m$, $m\in\R$. In other words, 
we can find $R>0$ so that for any $N\geq 0$, there is a constant
$C_N$ such that
\[
\left\|
\Rud_{V_1,V_2}(s) - \left(1+\sum_{k=0}^{N-1} \calR^{\pm}_k s^{-k-1}\right)\right\|
< \frac{C_N}{|s|^{N+1}}
\]
for $|s|>R$ in the corresponding domain, where $\|\cdot \|$ is a
fixed norm on $\End(V_1\otimes V_2)$.
\end{itemize}

Explicitly,
\begin{align*}
\Ru_{V_1,V_2}(s)&=\prod^{\rightarrow}_{n\geq 0} \opA{V_1}{V_2}(s+nl\hbar)^{-1}\\
\Rd_{V_1,V_2}(s)&=\prod^{\rightarrow}_{n\geq 1} \opA{V_1}{V_2}(s-nl\hbar)
\end{align*}
where the products converge uniformly on compact sets of $\pm
\Re(s/\hbar)\gg 0$ since $\opA{V_1}{V_2}(s) = 1 + O(s^{-2})$. Note that
the order of products indicated above is immaterial, since $\opA{V_1}
{V_2}(s)$ takes values in a commutative subalgebra of $\End(V_1\otimes
V_2)$.

\subsection{}\label{ssec: R0-main}

The following is the main result of this section.

\begin{thm}\label{thm: R0}
$\Rud_{V_1,V_2}(s)$ have the following properties
\begin{enumerate}
\item The map
\[\sigma\circ \Rud_{V_1,V_2}(s) :
V_1(s)\Dotimes{0}V_2 \rightarrow V_2\Dotimes{0} V_1(s)\]
where $\sigma$ is the flip of tensor factors, is a morphism of $\Yhg
$--modules, which is natural in $V_1$ and $V_2$.
\item For any $V_1,V_2,V_3\in\Ryang$ we have
\begin{align*}
\Rud_{V_1\Dotimes{s_1}V_2,V_3}(s_2) &= \Rud_{V_1,V_3}(s_1+s_2)
\Rud_{V_2,V_3}(s_2)\\
\Rud_{V_1,V_2\Dotimes{s_2}V_3}(s_1+s_2) &= \Rud_{V_1,V_3}(s_1+s_2)
\Rud_{V_1,V_2}(s_1)
\end{align*}
\item The following unitary condition holds
\[\sigma\circ \Rud_{V_1,V_2}(-s)\circ \sigma^{-1}=
\Rdu_{V_2,V_1}(s)^{-1}\]
\item For $a,b\in\C$ we have
\[\Rud_{V_1(a),V_2(b)}(s) = \Rud_{V_1,V_2}(s+a-b)\]
\item For any $s,s'$, \[[\Rud_{V_1,V_2}(s),\Rud_{V_1,V_2}(s')]=0=[\Rud_{V_1,V_2}(s),\Rdu_{V_1,V_2}(s')]\] 
\item $\Rud_{V_1,V_2}(s)$ have the same asymptotic expansion,
which is of the form
\begin{equation}\label{eq: R0-asym}
\Rud_{V_1,V_2}(s) \sim 1 + \hbar\Omega_{\h} s^{-1} + O(s^{-2})
\end{equation}
\item There is a $\rho>0$ such that the asymptotic expansion of
$\Rud_{V_1,V_2}(s)$ is valid on any domain
\[\left\{\pm\Re(s/\hbar)>m\right\}\cup\left\{|\Im(s/\hbar)|>\rho,\,\arg(\pm s/\hbar)\in(-\pi+\delta,\pi-\delta)\right\}\]
where $m\in\R$ and $\delta\in(0,\pi)$ are arbitrary.
\item The poles of $\Ru_{V_1,V_2}(s)^{\pm 1}$ and $\Rd_{V_1,V
_2}(s)^{\pm 1}$ are contained in
\[\spec(V_2)-\spec(V_1)-\Z_{\geq 0}l\hbar-\frac{\hbar}{2}\{l+r\}
\quad\text{and}\quad
\spec(V_2)-\spec(V_1)+\Z_{>0}l\hbar-\frac{\hbar}{2}\{l+r\}\]
where $r$ ranges over the integers such that $c_{ij}^{(r)}\neq 0$
for some $i,j\in\bfI$.

\end{enumerate}
\end{thm}
\begin{pf}
Part (i) is proved in \ref{ssec: commutativity} after some preparatory
results. Properties (ii)--(vi) and (viii) follow from Theorem \ref{thm: coeff}
and Section \ref{ssec: R0-prop}.
(vii) is proved in \cite[Lemma 8.1]{vanderput-singer}. 
\end{pf}

\subsection{Commutation relations with $\A_{V_1,V_2}(s)$}\label{ssec: op-A-comm}

Let $\calC\subset\C$ be a contour,
and $a_\ell:\C\to\End(V_\ell)$, $\ell=1,2$ two meromorphic functions
which are analytic within $\calC$ and commute with the operators
$\{\xi_{i,r}\}_{i\in\bfI,r\in\N}$. For any $k\in\bfI$, define operators
$X_k^{\pm,\ell}\in\End(V_1\otimes V_2)$ by
\[X_k^{\pm, 1} = \oint_{\calC}
a_1(v)x_k^{\pm}(v)\otimes a_2(v)\, dv
\quad\text{and}\quad
X_k^{\pm, 2} = \oint_{\calC}
a_1(v)\otimes a_2(v)x_k^{\pm}(v)\, dv \label{eq: Xk}
\]

\begin{prop}\label{prop: A-comm}
The following commutation relations hold
\begin{align*}
\Ad(\A_{V_1,V_2}(s))X_k^{\pm,1} &=
\oint_{\calC}
a_1(v)x_k^{\pm}(v)\otimes a_2(v)
\xi_k(v+s+l\hbar)^{\pm 1}\xi_k(v+s)^{\mp 1}\, dv\\
\Ad(\A_{V_1,V_2}(s))X_k^{\pm,2} &=
\oint_{\calC}
a_1(v)\xi_k(v-s)^{\pm 1}\xi_k(v-s-l\hbar)^{\mp 1}\otimes a_2(v)x_k^{\pm}(v)\, dv
\end{align*}
\end{prop}

\begin{pf}
We only prove the first relation. The second one follows from the first and
the unitarity property of Theorem \ref{thm: coeff}. We begin by computing
the commutation between $X_k^{\pm,1}$ and a typical summand in $\log
\A_{V_1,V_2}(s)$. Set $b=\pm\hbar d_ia_{ik}/2$. 
Note that the definition of $X_k^{\pm, 1}$ does not change if we replace
the contour $\calC$ by a smaller one $\calC'$, as long as both $\calC$
and $\calC'$ enclose the same set of poles of $x_k^{\pm}(v)$. Let $\calC_1$
be the contour chosen for the definition of $\A_{V_1,V_2}(s)$ given in 
Section \ref{ssec: op-A}. According
to Lemma \ref{lem: imp-lem}, if $v_0$ is a pole of $x_k^{\pm}(v)$ then 
$\calC_1$ must enclose $v_0\pm b$. Combining these observations, we will
assume, in the calculation below, that $\calC_1$ encloses $\calC$ and its
translates by $\pm b$. By \eqref{eq: imp-lem},
\[\begin{split}
[\oint_{\calC_1} t_i'(u)\otimes &t_j\lp u+s\rp\, du, X_k^{\pm,1}]\\
&=
\oint_{\calC_1}
\oint_\calC
a_1(v)[t_i'(u),x_k^\pm(v)]\otimes t_j(u+s)a_2(v)\, dv du\\
&=
\oint_{\calC_1}
\oint_\calC
\frac{1}{u-v+b} a_1(v)x_k^\pm(v)\otimes t_j(u+s)a_2(v)\, dv du\\
&
-\oint_{\calC_1}
\oint_\calC
\frac{1}{u-v-b} a_1(v)x_k^\pm(v)\otimes t_j(u+s)a_2(v)\, dv du\\
&+
\oint_{\calC_1}
\oint_\calC
\frac{1}{u-v-b} a_1(v)x_k^\pm(u-b)\otimes t_j(u+s)a_2(v)\, dv du\\
&-
\oint_{\calC_1}
\oint_\calC
\frac{1}{u-v+b} a_1(v)x_k^\pm(u+b)\otimes t_j(u+s)a_2(v)\, dv du\\
&=
\oint_\calC
a_1(v)x_k^\pm(v)\otimes (t_j(v-b+s)-t_j(v+b+s))a_2(v)\, dv
\end{split}\]
where the third equality follows from the fact that $s$ is such that
$t_j(u+s)$ is holomorphic inside ${\calC_1}$. Note that the third and
the fourth terms on the \rhs of the second equality vanish since their
integrands are holomorphic in the variable $v$.

Let the indeterminate $T$ of Section \ref{ssec: q-cartan} act as
the difference operator $Tt_j(v)=t_j (v-\hbar/2)$. Then,
\[\begin{split}
\sum_{i,j\in\bfI}&\left[\oint_{\calC_1} t_i'(u)\otimes c_{ij}(T)t_j\lp u+s\rp\, du, X_k^{\pm,1}\right]\\
&=
\sum_{i,j\in\bfI}
\oint_\calC
a_1(v)x_k^\pm(v)\otimes a_2(v) c_{ij}(T) (T^{\pm b_{ik}}-T^{\mp b_{ik}})t_j(v+s)\, dv\\
&=
\pm\oint_\calC
a_1(v)x_k^\pm(v)\otimes a_2(v) (T^l-T^{-l})t_k(v+s)\, dv\\
\end{split}\]
where the second equality follows from \eqref{eq:B inv}. The claimed
identity easily follows from this.
\end{pf}

\subsection{}\label{ssec: R0-comm}

Let $X_k^{\pm, 1},X_k^{\pm, 2}$ be the operators defined in \ref
{ssec: op-A-comm}. The following is a corollary of Proposition \ref
{prop: A-comm} and the definition of $\Rud(s)$.

\begin{prop}\label{prop: R0-comm}
The following commutation relations hold for any $\veps\in\{\pm\}$
\begin{align*}
\Ad(\Reps_{V_1,V_2}(s))X_k^{\pm,1} &= 
\oint_\calC
a_1(v)x_k^\pm(v)\otimes a_2(v)\xi_k(v+s)^{\pm 1}\, dv \\
\Ad(\Reps_{V_1,V_2}(s))X_k^{\pm,2} &= 
\oint_\calC
a_1(v)\xi_k(v-s)^{\mp 1}\otimes a_2(v)x_k^\pm(v)\, dv 
\end{align*}
\end{prop}

\subsection{Proof of (i) of Theorem \ref{thm: R0}}\label{ssec: commutativity}

We first rewrite the Drinfeld coproduct in a more symmetric way.
Let $V$ be a \fd representation of $\Yhg$ and $\calC^\pm\subset
\C$ a contour containing the poles of $x_i^\pm(u)$ on $V$. Then,
a simple contour deformation shows that, for any $u$ not contained
inside $\calC^\pm$,
\[\oint_{\calC^\pm} x_i^\pm(v)\frac{dv}{u-v}=x_i^\pm(u)\]
It follows that
\begin{align*}
\DD{s}(x_i^+(u)) &= 
\oint_{\calC_1}x_i^+(v-s)\otimes 1\,\frac{dv}{u-v}+
\oint_{\calC_2}\xi_i(v-s)\otimes x_i^+(v)\,\frac{dv}{u-v}\\
\DD{s}(x_i^-(u)) &= 
\oint_{\calC_1}x_i^-(v-s)\otimes\xi_i(v)\,\frac{dv}{u-v}+
\oint_{\calC_2} 1\otimes x_i^-(v)\,\frac{dv}{u-v}
\end{align*}
where $\calC_1,\calC_2$ are as in \ref{ssec: dt-yangian}.

We need to show that $\sigma\circ \Reps_{V_1,V_2}(s) : V_1(s)\otimes_0
V_2 \to V_2\otimes_0 V_1(s)$ intertwines the action of $\Yhg$.
This is obvious for $\xi_i(u)$, since $\xi_i(u)$ is group--like and
commutes with $\Reps_{V_1,V_2}(s)$. 
Denote now by $x_i^+(u)'$ and $x_i^+(u)''$, the action of $x_i^+(u)$
on $V_1(s)\otimes_0 V_2$ and $V_2\otimes_0 V_1(s)$ respectively.
By above formulas, we have
\begin{align*}
x_i^+(u)' &= \oint_{C_1} x_i^+(v-s)\otimes 1\, \frac{dv}{u-v} + 
\oint_{C_2}\xi_i(v-s)\otimes x_i^+(v)\, \frac{dv}{u-v} \\
x_i^+(u)'' &= \oint_{C_2} x_i^+(v)\otimes 1\, \frac{dv}{u-v} + 
\oint_{C_1}\xi_i(v)\otimes x_i^+(v-s)\, \frac{dv}{u-v}
\end{align*}
Using Proposition \ref{prop: R0-comm}, we can compute
$\Ad\lp\sigma\circ\Reps_{V_1,V_2}(s)\rp x_i^+(u)'$ as follows
\begin{align*}
&\sigma\lp \Reps_{V_1,V_2}(s)\lp
\oint_{C_1} x_i^+(v-s)\otimes 1\, \frac{dv}{u-v} + 
\oint_{C_2}\xi_i(v-s)\otimes x_i^+(v)\, \frac{dv}{u-v}
\rp \Reps_{V_1,V_2}(s)^{-1}\rp\sigma \\
&=\sigma\lp\oint_{C_1} x_i^+(v-s)\otimes \xi_i(v)\, \frac{dv}{u-v}
+ \oint_{C_2} 1\otimes x_i^+(v)\, \frac{dv}{u-v}
\rp\sigma \\
&= \oint_{C_1}\xi_i(v)\otimes x_i^+(v-s)\, \frac{dv}{u-v}
+\oint_{C_2} x_i^+(v)\otimes 1\, \frac{dv}{u-v} 
\end{align*}
This implies that $\Ad\lp\sigma\circ\Reps_{V_1,V_2}(s)\rp x_i^+(u)'=x_i^+(u)''$
and the result follows. The proof for $x_i^-(u)$ is identical.

\section{The functor $\Fh{}$}\label{sec: functor}

We review below the main construction of \cite{sachin-valerio-2}.
Assume henceforth that $\hbar\in\C\setminus\Q$, and that $q=
e^{\pi\iota\hbar}$.

\subsection{Difference equations}

Consider the abelian, additive difference equations, for unknown
functions $\phi_i : \C \to GL(V)$ 
\begin{equation}\label{eq: diff-eq-functor}
\phi_i(u+1)\mmu = \xi_i(u)\mmu \phi_i(u)\mmu
\end{equation}
defined by the commuting fields $\xi_i(u)\mmu=1+\hbar\xi_{i,0}u^{-1}+\cdots$
on a \fd representation $V$ of $\Yhg$.

Let $\phi^\pm_i(u)\mmu:\C\to GL(V)$ be the canonical fundamental solutions
of \eqref{eq: diff-eq-functor}. $\phi^\pm_i(u)\mmu$ are uniquely determined
by the requirement that they be holomorphic and invertible for $\pm\Re(u)\gg 0$,
and admit an asymptotic expansion of the form
\[\phi^\pm_i(u)\mmu\sim
(1+\varphi^\pm_0 u^{-1}+\varphi^\pm_1 u^{-2}\cdots)(\pm u)^{\hbar\xi_{i,0}}\]
in any right (resp. left) half--plane $\pm\Re(s)>m$, $m\in\R$ (see \eg \cite
{birkhoff-difference,borodin,krichever} or \cite[\S 4]{sachin-valerio-2}).
$\phi_i^+(u),\phi_i^-(u)$ are regularisations of the formal infinite products
\[ \xi_i(u)^{-1}\xi_i(u+1)^{-1}\xi_i(u+2)^{-1}\cdots
\quad\text{and}\quad
\xi_i(u-1)\xi_i(u-2)\xi_i(u-3)\cdots\]
respectively. 

Let $S_i(u)\mmu=(\phi^+_i(u)\mmu)^{-1}\phi^-_i(u)$ be the connection matrix
of \eqref{eq: diff-eq-functor}. Thus, $S_i(u)$ is 1--periodic in $u$, and therefore
a function of $z=\exp(2\pi\iota u)$. It is moreover regular at $z=0,\infty$
\cite[Prop. 4.8]{sachin-valerio-2}, and
therefore a rational function of $z$ such that
\[S_i(0)=e^{-\pi\iota\hbar\xi_{i,0}}=S_i(\infty)^{-1}\]
Explicitly,
\[S_i(u)=\lim_{n\to\infty}\xi_i(u+n)\cdots \xi_i(u+1)\xi_i(u)\xi_i(u-1)\cdots\xi_i(u-n)\]

\subsection{Non--congruent representations}

We shall say that $V\in\Ryang$ is {\em non--congruent} if,
for any $i\in\bfI$, the poles of $x^+_i(u)$ (resp. $x^-_i(u)$)
are not congruent modulo $\Ztimes$. Let $\Rync$ be the full
subcategory of $\Ryang$ consisting of non--congruent
representations.

\subsection{The functor $\Fh{}$}\label{ssec: functor}

Given $V\in\Rync$, define the action of the generators of $\qloop$
on  $\Fh{}(V)=V$ as follows.\\
\begin{enumerate}
\item For any $i\in\bfI$, the generating series $\Psi_i(z)^+\mmu$
(resp. $\Psi_i(z)^-\mmu$) of the commuting generators of $\qloop$
acts as the Taylor expansions at $z=\infty$ (resp. $z=0$) of the
rational function
\[\Psi_i(z)\mmu=\left.S_i(u)\mmu\right|_{e^{2\pi\iota u}=z}\]
\end{enumerate}

\noindent
To define the action of the remaining generators of $\qloop$,
let $g_i^\pm(u)\mmu:\C\to GL(V)$ be given by $g_i^+(u)\mmu
=\phi_i^+(u+1)\mmu^{-1}$ and $g_i^-(u)\mmu=\phi_i^-(u)\mmu$.
Explicitly,
\begin{equation}\label{eq:gipm}
\begin{split}
g_i^+(u) &= \lp \prod_{n\geq 1}^{\leftarrow}
\xi_i(u+ n) \,e^{-\hbar\xi_{i,0}/n}\rp e^{\gamma\hbar\xi_{i,0}}\\
g_i^-(u) &= e^{-\gamma\hbar\xi_{i,0}}\lp\prod_{n\geq 1}^{\rightarrow}
\xi_i(u - n) \,e^{\hbar\xi_{i,0}/n}\rp
\end{split}
\end{equation}
where $\ds{\gamma=\lim_{n\to\infty}(1+\cdots+1/n-\log n)}$
is the Euler--Mascheroni constant, are regularisations of the infinite
products
\[\cdots\xi_i(u+2)\xi_i(u+1)
\aand
\xi_i(u-1)\xi_i(u-2)\cdots\]
Note also that, by definition of $g_i^\pm(u)$
\begin{equation}\label{eq:3 terms}
S_i(u)=g_i^+(u)\cdot\xi_i(u)\cdot g_i^-(u)
\end{equation}

\noindent
Let $c_i^{\pm}\in\nC$ be scalars such that $c_i^-c_i^+=d_i\Gamma
(\hbar d_i)^2$.\\ 

\begin{enumerate}
\item[(ii)] For any $i\in\bfI$ and $k\in \Z$, $\X_{i,k}^\pm$ acts as
the operator
\[\X^{\pm}_{i,k}\mmu
=
c_i^\pm\oint_{\calC^\pm_{i\mmu}}
e^{\pp ku}g_i^{\pm}(u)\mma x_i^{\pm}(u)\mmu\, du\]
where the Jordan curve $\calC^\pm_{i\mmu}$ encloses the poles
of $x^\pm_i(u)\mmu$ and none of their $\Ztimes$--translates.
\footnote{Note that such a curve exists for any $i\in\bfI$ 
since $V$ is non--congruent.} The corresponding generating series
are the expansions at $z=\infty,0$ of the $\End(V)$--valued rational 
function given by
\[\X^\pm_i(z)=
c_i^\pm\oint_{\calC^\pm_{i\mmu}}
\frac{z}{z-e^{2\pi\iota u}}g_i^{\pm}(u)\mma x_i^{\pm}(u)\mmu\, du\]
where $z$ lies outside of $\exp(2\pi\iota\calC^\pm_{i\mmu})$.
\end{enumerate}

\subsection{}

Let $\Pi\subset\C$ be a subset such that $\Pi\pm\frac{\hbar}{2}\subset\Pi$.
Let
\[\Rysub\subset\Ryang\]
be the full subcategory of  consisting of the representations $V$ such that
$\spec(V)\subset\Pi$.

Similarly, let $\Omega\subset\nC$ be a subset stable under multiplication by
$q^{\pm 1}$. We define $\Rlsub$ to be the full subcategory of $\Rloop$ consisting
of those $\V$ such that $\spec(\V)\subset\Omega$.

\subsection{}

\begin{thm}\label{thm: functor}\cite[Thm. 5.4, Thm. 6.3, Prop. 7.7]{sachin-valerio-2}\hfill\break
\begin{enumerate}
\item The above operators give rise to an action of $\qloop$ on
$V$. They therefore define an exact, faithful functor
\[\Fh{}:\Rync\longrightarrow\Rloop\]
\item The functor $\Fh{}$ is compatible with shift automorphisms.
That is, for any $V\in\Rync$ and $a\in\C$,
\[\Fh{}(V(a)) = \Fh{}(V)(\Exp{a})\]
\item Let $\Pi\subset\C$ be a non--congruent subset such that
$\Pi\pm\frac{1}{2}\hbar\subset\Pi$. Then, $\Rysub$ is a subcategory
of $\Rync$, and $\Fh{}$ restricts to an isomorphism of abelian
categories.
\[\Fh{\Pi} : \Rysub\isom\Rlsub\]
where $\Omega=\exp(2\pi\iota\Pi)$.
\item $\Fh{\Pi}$ is compatible with the $q$--characters of Knight and
Frenkel--Reshetikhin.
\end{enumerate}
\end{thm}

\section{Meromorphic tensor structure on $\Fh{}$}\label{sec: tensor}

\subsection{The abelian $q$KZ equations}\label{ssec:qKZ}

Let $V_1,V_2$ be \fd representations of $\Yhg$, choose $\veps
\in\{\pm\}$, and let $\Reps_{V_1,V_2}(s)$ be the corresponding
$R$--matrix defined in Section \ref{ssec: R0-prop}. Consider the abelian,
additive $q$KZ equation for an unknown function $f : \C\to\End(
V_1\otimes V_2)$
\begin{equation}\label{eq:ab qkz}
f(s+1)=\Reps_{V_1,V_2}(s)f(s)
\end{equation}

Note that this equation does not fit the usual assumptions in the
study of difference equations since $\Reps_{V_1,V_2}(s)$ is not
rational. Moreover, $\Reps_{V_1,V_2}(s)$ may not have a Laurent
expansion at $\infty$ but, by Theorem \ref{thm: R0}, only an asymptotic
expansion of the form $1+\hbar\Omega_\h/s+O(s^{-2})$ valid in
any domain of the form
\[\left\{\Re(s/\veps\hbar)>m\right\}
\quad\cup\quad
\left\{|\Im(s/\hbar)|>\rho,\,\arg(s/\veps\hbar)\in(-\pi+\delta,\pi-\delta)\right\}\]
where $\rho>0$ is fixed, and $m\in\R,\delta\in(0,\pi)$ are arbitrary.\footnote
{For the $q$KZ equations determined by the full $R$--matrix, these issues
are usually addressed by proving the existence of factorisation $R_{V_1,V_2}
(s)=R\rat_{V_1,V_2}(s)\cdot R\mer_{V_1,V_2}(s)$, where $R\rat
_{V_1,V_2}(s)$ is a rational function of $s$ which intertwines the
Kac--Moody coproduct $\Delta$ and its opposite,
and the meromorphic factor $R\mer_{V_1,V_2}(s)$ intertwines
$\Delta$ (see \cite{kazhdan-soibelman} for the case of $\qloop$), and then working with $R\rat_{V_1,V_2}(s)$ instead
of $R_{V_1,V_2}(s)$. A similar factorisation can be obtained
for the abelian $R$--matrices $\Rud(s)$. We shall, however,
prove in \cite{sachin-valerio-qKL} that neither of these factorisations
are natural \wrt $V_1,V_2$, which is why we work with the
meromorphic $R$--matrices $\Rud(s)$.}
Nevertheless, these asymptotics and the fact that the poles of
$\Reps(s)^{\pm 1}$ are contained in the complement of a domain
of the above form, are sufficient to carry over the standard proofs
(see, \eg \cite[\S 4]{sachin-valerio-2}) and yield the following. 

\begin{prop}\label{pr:left right}
Let $\sfn\in\C^\times$ be perpendicular to $\hbar$ and such that
$\Re(\sfn)\geq 0$.
\begin{enumerate}
\item If $\veps\hbar\notin\R_{<0}$, the equation \eqref{eq:ab qkz}
admits a canonical right meromorphic solution $\Phi^\veps_+:\C\to
GL(V_1\otimes V_2)$, which is uniquely determined by the following
requirements
\begin{itemize}
\item $\Phi^\veps_+$ is holomorphic and invertible for $\Re(s)\gg 0$
if $\Re(\veps\hbar)\geq 0$, and otherwise on a sector of the form
\begin{equation}\label{eq:angle sector 1}
\Re(s) \gg 0\aand\Re(s/\sfn)\gg 0
\end{equation}
\item $\Phi^\veps_+$ has an asymptotic expansion of the form $(1+
O(s^{-1}))s^{\hbar\Omega_\h}$ in any right half--plane if $\Re(
\veps\hbar)>0$, and otherwise in a sector of the form \eqref
{eq:angle sector 1}.
\end{itemize}
\item If $\veps\hbar\notin\R_{>0}$, the equation \eqref{eq:ab qkz}
admits a canonical left meromorphic solution $\Phi^\veps_-:\C\to
GL(V_1\otimes V_2)$, which is uniquely determined by the following
requirements
\begin{itemize}
\item $\Phi^\veps_-$ is holomorphic and invertible for $\Re(s)\ll 0$
if $\Re(\veps\hbar)\leq 0$, and otherwise on a sector of the form
\begin{equation}\label{eq:angle sector 2}
\Re(s) \ll 0\aand\Re(s/\sfn)\ll 0
\end{equation}
\item $\Phi^\veps_-$ has an asymptotic expansion of the form
$(1+O(s^{-1}))(-s)^{\hbar\Omega_\h}$ in any left half--plane
if $\Re(\veps\hbar)<0$, and otherwise in a sector of the form
\eqref{eq:angle sector 2}.
\end{itemize}
\end{enumerate}
\end{prop}

\begin{figure}[h!]\label{fig: zones}
\includegraphics[width=5in]{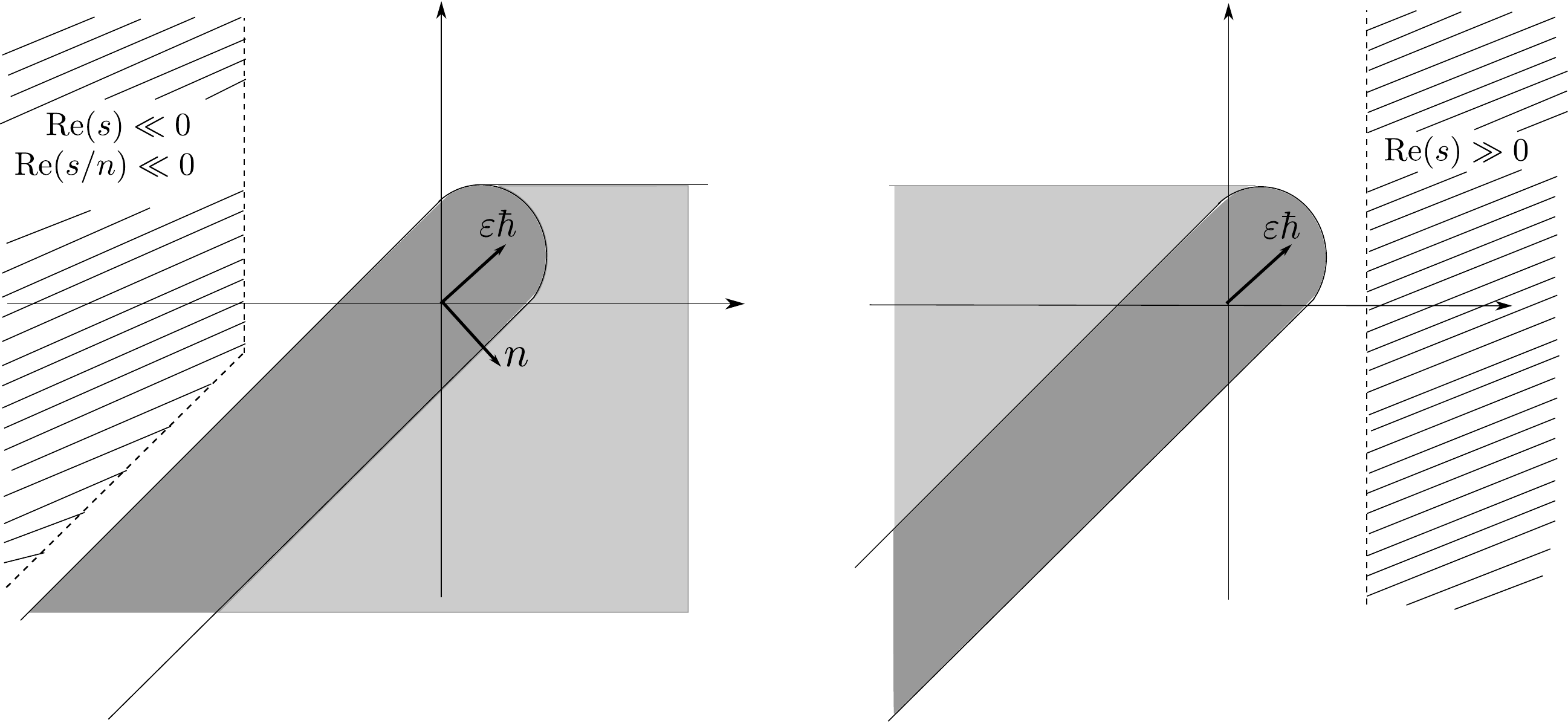}
\caption{Domains of holomorphy and invertibility of $\Phi^{\veps}_+$
(resp. $\Phi^{\veps}_-$) given by the ruled region in the right (resp.
left) picture, when $\Re(\veps\hbar)>0$. The darker grey region
contains poles of $\Reps(s)^{\pm 1}$.}
\end{figure}

The right and left solutions, when defined, are given by the products
\begin{align}
\Phi_+^\veps(s)&=
e^{-\hbar\gamma\Omega_\hbar}\Reps_{V_1,V_2}(s)^{-1}\stackrel{\longrightarrow}{\prod}_{m\geq 1}\Reps_{V_1,V_2}(s+m)^{-1}\,e^{\hbar\Omega_\h/m}
\label{eq:ab qkz+}\\
\Phi_-^\veps(s)&=
e^{-\hbar\gamma\Omega_\hbar}\stackrel{\longrightarrow}{\prod}_{m\geq 1}\Reps_{V_1,V_2}(s-m)\,e^{\hbar\Omega_\h/m}
\label{eq:ab qkz-}
\end{align}

\subsection{Proof of Proposition \ref{pr:left right}}

As mentioned before, the proof follows the same strategy as in \cite[\S 4]
{sachin-valerio-2}. More precisely, we use the fact that $\hbar\Omega_{\h}$
commutes with $\Reps_{V_1,V_2}(s)$ to regularize \eqref{eq:ab qkz}, that
is set (as in \cite[\S 4.6]{sachin-valerio-2}) 
\[
\ol{\Reps_{V_1,V_2}}(s) := (1-\hbar\Omega_{\h}s^{-1})\Reps_{V_1,V_2}(s)
\]
The auxiliary equation $f(s+1) = (1-\hbar\Omega_{\h}s^{-1})f(s)$ can be solved
using the $\Gamma$--function (see \cite[\S 4.5, 4.6]{sachin-valerio-2}), while
the regularized equation (equation \eqref{eq:ab qkz} with $\Reps_{V_1,V_2}$ replaced
by $\ol{\Reps_{V_1,V_2}}$) is solved by taking the infinite products \cite[\S 4.4]{sachin-valerio-2}:
\begin{align*}
\ol{\Phi^{\veps}_+}(s) &= \prod^{\rightarrow}_{n\geq 0} \ol{\Reps_{V_1,V_2}}(s+n)^{-1} \\
\ol{\Phi^{\veps}_-}(s) &= \prod^{\rightarrow}_{n\geq 1} \ol{\Reps_{V_1,V_2}}(s-n)
\end{align*}

This is the only point of departure from the rational case. In order to prove
the convergence of these infinite products, we only need the asymptotics
of $\ol{\Reps_{V_1,V_2}}$ up to the second order in the desired zones, as stated in the
following lemma. Its proof is standard and hence omitted. 

\begin{lem}\label{lem: infinitesumsconverge}
Let $\Omega\subset\C$ be an open set, $W$ a \fd complex vector space,
and $f : \Omega\to\End(W)$ a holomorphic and invertible function such
that the following assumptions hold.
\begin{itemize}
\item[(a)] For each $n\in\Z_{\geq 0}$, $\Omega+n \subset \Omega$.
\item[(b)] There exists a constant $C\in\R_{>0}$ such that
\[
\|f(s)-1\| < \frac{C}{|s|^2} \quad\text{as}\quad s\to\infty, s\in\Omega
\]
for some norm $\|\cdot\|$ on $\End(W)$. 
\end{itemize}
Then the sequence of functions $\{f(s)f(s+1)\cdots f(s+n)\}_{n\geq 1}$
converges uniformly on compact sets in $\Omega$ and hence defines a
holomorphic function $F(s)$ on $\Omega$.

If, in addition,
\begin{itemize}
\item[(c)] $f(s)$ extends to a meromorphic function on $\C$.
\item[(d)] $\Omega$ contains a fundamental domain for $s\mapsto s+1$.
\end{itemize}
then $F(s)$ can be extended to a meromorphic function on $\C$ by
using the equation $F(s) = f(s)F(s+1)$.

The same assertions hold for the infinite product $f(s-1)f(s-2)\cdots$
after changing $\Z_{\geq 0}$ to $\Z_{\leq 0}$ in condition (a) above.
\end{lem}

This, in particular, explains that we have to consider sectors given
in Figure \ref{fig: zones} in order to avoid the poles of 
$\Reps_{V_1,V_2}(s)^{\pm 1}$. Thus we obtain the solutions
$\Phi^{\veps}_{\pm}$ of the difference equation \eqref{eq:ab qkz},
which are explicitly given in \eqref{eq:ab qkz+} and \eqref{eq:ab qkz-},
and whose asymptotics can be computed using the calculation in \cite
[\S 4.7]{sachin-valerio-2}.

\subsection{The tensor structure $\Jeps_{V_1,V_2}(s)$}\label{ssec: twist-defn}

Let $\veps\in\{\pm\}$ be such that $\veps\hbar\notin\R_{<0}$, and
$\Phi^\veps_+(s)$ the right fundamental solution of the abelian $q$KZ
equation \eqref{eq:ab qkz}. Define a meromorphic function
\[\Jeps_{V_1,V_2}:\C\to GL(V_1\otimes V_2)\]
by $\Jeps_{V_1,V_2}(s)=\Phi^\veps_+(s+1)^{-1}$. Thus,
\begin{equation}\label{eq:J}
\Jeps_{V_1,V_2}(s)=
e^{\hbar\gamma\Omega_{\h}}\stackrel{\longleftarrow}{\prod}_{m\geq 1} \Reps_{V_1,V_2}(s+m) 
e^{-\frac{\hbar\Omega_{\h}}{m}}
\end{equation}

\begin{thm}\label{thm: tensor}\hfill
\begin{enumerate}
\item $\Jeps_{V_1,V_2}(s)$ is natural in $V_1,V_2$.
\item If $V_1$ and $V_2$ are non--congruent, and $\zeta = e^{2\pi\iota s}$,
\[
\Jeps_{V_1,V_2}(s) : \Fh{}(V_1)\Dotimes{\zeta} \Fh{}(V_2)\longrightarrow\Fh{}(V_1\Dotimes{s} V_2)
\]
is an isomorphism of $\qloop$--modules for any $s\not\in\spec(V_2)-\spec(V_1)+\Z$.
\item For any non--congruent $V_1,V_2,V_3\in\Ryang$, the following
is a commutative diagram
\[\xymatrix{
(\Fh{}(V_1)\Dotimes{\zeta_1}\Fh{}(V_2))\Dotimes{\zeta_2}\Fh{}(V_3)\ar[dd]_{\Jeps_{V_1,V_2}(s_1)\otimes\id}
\ar@{=}[r]&
\Fh{}(V_1)\Dotimes{\zeta_1\zeta_2}(\Fh{}(V_2)\Dotimes{\zeta_2}\Fh{}(V_3))\ar[dd]^{\id\otimes\Jeps_{V_2,V_3}(s_2)}\\ & \\
\Fh{}(V_1\Dotimes{s_1}V_2)\Dotimes{\zeta_2}\Fh{}(V_3)\ar[dd]_{\Jeps_{V_1\Dotimes{s_1}V_2,V_3}(s_2)}
&
\Fh{}(V_1)\Dotimes{\zeta_1\zeta_2}\Fh{}(V_2\Dotimes{s_2}V_3)\ar[dd]^{\Jeps_{V_1,V_2\Dotimes{s_2}V_3}(s_1+s_2)}\\ & \\
\Fh{}((V_1\Dotimes{s_1}V_2)\Dotimes{s_2}V_3)
\ar@{=}[r]&
\Fh{}(V_1\Dotimes{s_1+s_2}(V_2\Dotimes{s_2}V_3))
}\]
where $\zeta_i=\exp(2\pi\iota s_i)$.
\item The poles of $\Ju_{V_1,V_2}(s)^{\pm 1}$ and $\Jd_{V_1,
V_2}(s)^{\pm 1}$ are contained in
\[\spec(V_2)-\spec(V_1)-\Z_{\geq 0}l\hbar-\frac{\hbar}{2}\{l+r\}-\Z_{> 0}
\quad\text{and}\quad
\spec(V_2)-\spec(V_1)+\Z_{>0}l\hbar-\frac{\hbar}{2}\{l+r\}-\Z_{> 0}\]
where $r$ ranges over the integers such that $c_{ij}^{(r)}\neq 0$
for some $i,j\in\bfI$.
\end{enumerate}
\end{thm}

\begin{rem}
Note that the condition $s\not\in\spec(V_2)-\spec(V_1)+\Z$ implies that 
$V_1\Dotimes{s} V_2$ exists and is non--congruent, which is required in order
to define $\Fh{}(V_1\Dotimes{s} V_2)$.
\end{rem}

\begin{pf}
(i) and (iii)--(iv) follow from \eqref{eq:J} and Theorem \ref{thm: R0}. 
(ii) is proved in \ref{ss:ii}.
\end{pf}

\subsection{}\label{ss:ii}

%
Given an element $X\in \qloop$, we denote its action on $\Fh{}
(V_1)\Dotimes{\zeta} \Fh{}(V_2)$ and $\Fh{}(V_1\Dotimes{s}V
_2)$ by $X'$ and $X''$ respectively. We need to prove that
\[\Jeps_{V_1,V_2}(s)X'\Jeps_{V_1,V_2}(s)^{-1} = X''\]
Since $\xi_i(u)$ are group--like \wrt the Drinfeld coproduct, so
are the fundamental solutions and the connection matrix of the
difference equation $\phi_i(u+1) = \xi_i(u)\phi_i(u)$, which implies
that $\Psi_i(z)' = \Psi_i(z)''$. Since $\Rud_{V_1,V_2}(s)$ and hence $\Jeps_{V_1,V_2}(s)$
commute with these elements, this proves the required relation
for $\{\Psi_i(z)\}_{i\in\bfI}$.

We now prove the relation for $\X_{i,k}^+$. The proof for $\X_{i,k}
^-$ is similar. By \ref{ssec: dt-qla2} and \ref{ssec: functor}, the action
of $(c_i^+)^{-1}\X^+_{i,k}$ on $\Fh{}(V_1)\Dotimes{\zeta}\Fh{}(V_2)$
is given by
\[\begin{split}
\zeta^k\oint_{C_1}e^{2\pi\iota k u}g_i^+(u)&x_i^+(u)\otimes 1\,du\\
&+
\oint_{\calC_2}\Psi_i(\zeta^{-1}w)\otimes\oint_{C_2} g_i^+(u)x_i^+(u)\frac{w}{w-e^{2\pi\iota u}}w^{k-1}\,dwdu\\
=
\zeta^k\oint_{C_1}e^{2\pi\iota k u}g_i^+(u)&x_i^+(u)\otimes 1\,du\\
&+
\oint_{C_2}e^{2\pi\iota k u} g_i^+(u-s)\xi_i(u-s)g_i^-(u-s)\otimes g_i^+(u)x_i^+(u)\,du
\end{split}\]
where
\begin{itemize}
\item $C_\ell$ encloses $\spec(V_\ell)$ and none of its $\Ztimes$--translates.
\item $\calC_2$ encloses $C_2$, $\exp(2\pi\iota\spec(V_2))$
and none of the points in $\exp(2\pi\iota(s+\spec(V_1)))$. Note that these 
sets contain $\spec(\Fh{}(V_2))$ and $\zeta\spec(\Fh{}(V_1))$ by definition.
We also remark that we are assuming $s\not\in \spec(V_2)-\spec(V_1)+\Z$
in (ii) of Theorem \ref{thm: tensor} which makes such a choice of contours
possible.
\end{itemize}
and we used \eqref{eq:3 terms}.

On the other hand, the action of $(c_i^+)^{-1}\X^+_{i,k}$ on $\Fh{}(V_1
\Dotimes{s}V_2)$ is given by
\[\begin{split}
&\oint_{C_{12}} e^{2\pi\iota k u}g_i^+(u-s)\otimes g_i^+(u)
\left(x_i^+(u-s)\otimes 1+\oint_{C_2'}\xi_i(v-s)\otimes x_i^+(v)\frac{dv}{u-v}\right)\,du\\[1.2ex]
=&
\zeta^k\oint_{C_1} e^{2\pi\iota k u} g_i^+(u)x_i^+(u)\otimes g_i^+(u+s)\,du\\
&\phantom{\zeta^k\oint_{C_1} e^{2\pi\iota k u} g_i^+(u)}
+
\oint_{C_2} e^{2\pi\iota k v} g_i^+(v-s)\xi_i(v-s)\otimes g_i^+(v)x_i^+(v)\,dv
\end{split}\]
where
\begin{itemize}
\item $C_{12}$ encloses $\lp\spec(V_1)+s\rp\cup\spec(V_2)$
(which contains $\spec(V_1\Dotimes{s}V_2)$)
and none of its $\Ztimes$--translates. Again it is possible
thanks to our assumption on $s$ imposed in (ii) of Theorem \ref{thm: tensor}
above.
\item $C_2'$ encloses $\spec(V_2)$ and none of the points of $\spec(V_1)+s$.
\end{itemize}
$C_1$ is as above, and we assumed that $C_{12}$ encloses $C_2'$, and that
$C_2'=C_2$.

Let us compute the action of $\Ad(\Jeps_{V_1,V_2}(s))$ on the first
summand of $(c_i^+)^{-1}(\X^+_{i,k})'$. 
Note that $\ad(\Omega_{\h})\,x_i^+(v)\otimes 1 = \sum_{j\in\bfI}
[\varpi_j^{\vee},x_i^+(v)]\otimes \xi_{j,0} = x_i^+(v)\otimes \xi_{i,0}$,
by equation \eqref{eq: omega-cartan}.
Therefore, for any $a\in\C$
\[\Ad(e^{a\Omega_{\h}})\,x_i^+(v)\otimes 1= e^{\ad(a\Omega_{\h})}\, x_i^+(v)\otimes 1
=x_i^+(v)\otimes e^{a\xi_{i,0}}\]
Using this and Proposition \ref{prop: R0-comm}
we get
\[\begin{split}
\Ad(\Jeps_{V_1,V_2}(s))
&\left(\zeta^k\oint_{C_1} e^{2\pi\iota ku}g_i^+(u)x_i^+(u)\otimes 1\, du\right)\\ 
&=\zeta^k\oint_{C_1} e^{2\pi\iota ku}g_i^+(u)x_i^+(u)\otimes e^{\gamma\hbar\xi_{i,0}}
\prod_{n\geq 1} \xi_i(u+s+n)e^{-\hbar\xi_{i,0}/n} \, du\\
&=\zeta^k\oint_{C_1} e^{2\pi\iota ku}g_i^+(u)x_i^+(u)\otimes g_i^+(u+s)\, du
\end{split}\]
by the definition of $g_i^+(u)$ given in \eqref{eq:gipm}. This yields
the first term on the \rhs of $(c_i^+)^{-1}(\X^+_{i,k})''$. A similar 
computation can be carried out for the second summand of
$(c_i^+)^{-1}(\X^+_{i,k})'$ which proves that
\[\Jeps_{V_1,V_2}(s)(\X_{i,k}^+)'\Jeps_{V_1,V_2}(s)^{-1} = (\X_{i,k}^+)'' \]

\section{The commutative $R$--matrix of the quantum loop algebra}\label{sec: R0qla}

In this section, we review the construction of the commutative part
$\Rq^0(\zeta)$ of the $R$--matrix of the quantum loop algebra. We
prove that if $|q|\neq 1$, $\Rq^0(\zeta)$ defines a meromorphic
commutativity constraint on $\Rloop$, when  the latter is equipped
with the Drinfeld tensor product studied in \S\ref{sec: dt}.

\subsection{Drinfeld pairing}

The Drinfeld pairing for the quantum loop algebra was computed
in terms of the loop generators by Damiani \cite{damiani}. Its
restriction to $U^0$ is given in \cite[Corollary 9]{damiani}. Define
$\{H_{i,r}\}_{i\in\bfI, r\in\Z_{\neq 0}}$ by
\begin{equation}\label{eq:psi-H}
\Psi_i^{\pm}(z) = \Psi^{\pm}_{i,0}\,\exp\lp\pm(q_i-q_i^{-1})
\sum_{r\geq 1} H_{i,\pm r} z^{\mp r}\rp
\end{equation}
Then, for each $m,n\geq 1$
\begin{equation}\label{eq:pairingHmodes}
\la H_{i,m}, H_{j,-n}\ra = -\delta_{m,n}\frac{q^{mb_{ij}}-q^{-mb_{ij}}}
{m(q_i-q_i^{-1})(q_j-q_j^{-1})} 
\end{equation}
where $b_{ij} = d_ia_{ij} = d_ja_{ji}$. Define $H_i^{\pm}(z)\in z^{\mp 1}
U^0[[z^{\mp 1}]]$ by
\[
H_i^{\pm}(z) = \pm (q_i-q_i^{-1})\sum_{r\geq 1} H_{i,\pm r}z^{\mp r}
\]
Then, by \eqref{eq:pairingHmodes},
\begin{equation}\label{eq:pairingHfields}
\la H_i^+(z), H_j^-(w)\ra = 
\sum_{m\geq 1}\frac{q^{mb_{ij}}-q^{-mb_{ij}}}{m}\left(\frac{w}{z}\right)^m=
\log\lp\frac{z-q^{-b_{ij}}w}{z-q^{b_{ij}}w}\rp
\end{equation}

\subsection{Construction of $\Rq^0$}

We now follow the argument outlined in \S \ref{ss:KT} to
construct the canonical element $\Rq^0$ of this pairing. Namely,
differentiating \eqref{eq:pairingHfields} with respect to $z$ yields
\[
\la \frac{d}{dz}H_i^+(z), H_j^-(w)\ra = \frac{1}{z-q^{-b_{ij}}w}
- \frac{1}{z-q^{b_{ij}}w} = (T^{b_{ij}}-T^{-b_{ij}})\frac{1}{z-w}
\]
where $Tf(z,w) = f(z,q^{-1}w)$ is the multiplicative shift operator
\wrt $w$. Hence, if we define
\begin{equation}\label{eq:Hdual}
H^{j,-}(w) = (T^l-T^{-l})^{-1}\sum_{k\in\bfI} c_{jk}(T)H_k^-(w)
\,\in wU^0[[w]]
\end{equation}
where $(T^l-T^{-l})^{-1}$ acts on $w^k$, $k\neq 0$, as multiplication
by $(q^{-lk}-q^{lk})^{-1}$, then
\[
\la \frac{d}{dz}H_i^+(z), H^{j,-}(w)\ra = \delta_{ij}\frac{1}{z-w}
\]
Note that $H^{j,-}(w)$ is explicitly given by
\[H^{j,-}(w)=
\sum_{k\in\bfI}(q_k-q_k^{-1})
\sum_{n\geq 1} \lp 
\frac{c_{jk}(q^n)}{q^{nl}-q^{-nl}}
H_{k,-n} \rp w^n\]
so that $\Rq^0$ is equal to
\[\Rq^0 = q^{-\Omega_{\h}} \exp\lp
-\sum_{\begin{subarray}{c} i,j\in\bfI\\ m\geq 1\end{subarray}}
\frac{m(q_i-q_i^{-1})(q_j-q_j^{-1})c_{ij}(q^m)}
{q^{ml}-q^{-ml}} H_{i,m}\otimes H_{j,-m}
\rp
\]

\subsection{$q$--difference equation for $\Rq^0$}\label{ssec: qDER0}

Set $\Rq^0(\zeta) = (\tau_{\zeta}\otimes 1)\Rq^0$, so that
\begin{multline}\label{eq:R0qla-formal}
\Rq^0(\zeta)=\\
q^{-\Omega_{\h}} \exp\lp
-\sum_{\begin{subarray}{c} i,j\in\bfI\\ m\geq 1\end{subarray}}
\frac{m(q_i-q_i^{-1})(q_j-q_j^{-1})c_{ij}(q^m)}
{q^{ml}-q^{-ml}}\zeta^m H_{i,m}\otimes H_{j,-m}
\rp
\end{multline}
It is easy to see that $\Rq^0(\zeta)$ satisfies the following $q
$--difference equation
\begin{multline}\label{eq: qDER0-formal}
\Rq^0(q^{2l}\zeta)\Rq^0(\zeta)^{-1}\\
= \exp\lp
-\sum_{\begin{subarray}{c} i,j\in\bfI \\ m\geq 1\end{subarray}}
m(q_i-q_i^{-1})(q_j-q_j^{-1})c_{ij}(q^m)q^{ml}\zeta^m
H_{i,m}\otimes H_{j,-m}
\rp
\end{multline}

\subsection{}\label{ss:interpret H+ H-}

Using the method employed in \ref{ss:our strategy}--\ref{ssec: op-A},
we will show that the \rhs of \eqref{eq: qDER0-formal} is the expansion
of a rational function at $\zeta=0$, once it is evaluated on a tensor
product of \fd representations.

We note first that a typical summand may be interpreted as a contour
integral as follows
\begin{multline*}
\sum_{m\geq 1} m(q_i-q_i^{-1})(q_j-q_j^{-1})c_{ij}(q^m)q^{ml}\zeta^m
H_{i,m}\otimes H_{j,-m}\\
 = \sum_{r\in\Z} c_{ij}^{(r)} \oint_{\calC} \frac{dH_i^+(w)}{dw} 
\otimes H_j^-(q^{l+r}\zeta w)\, dw
\end{multline*}

On a tensor product of two \fd representations $\V_1,\V_2$, the
first tensor factor is a rational function of $w$ since
\[\frac{dH_i^+(w)}{dw}  = \Psi_i(w)^{-1}\frac{d\Psi_i(w)}{dw}\]
The second tensor factor $H_j^-(q^{l+r}\zeta w)$ can be viewed
as a single--valued function defined outside of a set of cuts radiating
from $\zeta=\infty$. To see this, note that 
$H_j^-(w)$ is a logarithm of the rational $\End(\V_2)$--valued function
$\Psi^+_{j,0}\Psi_j(w)$, and that the latter is regular at $w=0,\infty$
and takes the value $1$ at $w=0$. The result then follows from
the variant of Proposition \ref{pr:matrix log} below.

\begin{prop}\label{pr:matrix log-trig}
Let $\V$ be a complex, \fd vector space, and $\psi:\C\to\End(\V)$
a rational function, regular at $0$ and $\infty$ such that
\begin{itemize}
\item $\psi(0) = 1$.
\item $[\psi(w),\psi(w')]=0$ for any $w,w'\in\C$.
\end{itemize}
Let $\spec(\psi)\subset\C^\times$ be the set of poles of $\psi(w)^
{\pm 1}$, and define the cut--set $\sfY(\psi)$ by
\[
\sfY(\psi) = \bigcup_{\alpha\in\spec(\psi)} [\alpha,\infty)
\]
where $[\alpha,\infty) = \{t\alpha: t\in\R_{\geq 1}\}$.
Then, there is a unique single--valued, holomorphic function
$H(w) = \logo(\psi(w)) : \C\setminus\sfY(\psi)\to\End(\V)$ such that
\[
\exp(H(w)) = \psi(w) \aand H(0) = 0
\]
Moreover, $[H(w),H(w')]=0$ for any $w,w'\in\C$ and $\frac{dH}{dw}
= \psi^{-1}\frac{d\psi}{dw}$.
\end{prop}

The proof of this proposition is analogous to that of Proposition \ref
{pr:matrix log}.

\subsection{The operator $\Aq_{\V_1,\V_2}(\zeta)$}\label{ssec: op-B}

Let $\V_1,\V_2$ be two \fd representations of $\qloop$. Let $\calC_1$
be a contour enclosing the set of poles $\spec(\V_1)$ of $\V_1$, and
consider the following operator on $\V_1\otimes\V_2$
\[
\Aq_{\V_1,\V_2}(\zeta) = \exp\lp
-\sum_{\begin{subarray}{c} i,j\in\bfI\\ r\in\Z\end{subarray}}
c_{ij}^{(r)} \oint_{\calC_1} \frac{dH_i^+(w)}{dw}
\otimes H_j^-(q^{l+r}\zeta w)\, dw
\rp
\]
where
\begin{itemize}
\item $\ds{\frac{dH_i^+}{dw}:\C\to\End(\V_1)}$ is the rational function
$\ds{\Psi_i^{-1}\frac{d\Psi_i}{dw}}$,
\item $H_j^-
:\C\setminus\sfY(\Psi_{j,0}^+\Psi_j(w))
\to\End(\V_2)$ is given by Proposition \ref{pr:matrix log-trig},
\item $\zeta\in\C$ is small enough so that $H_j^-(q^{l+r}\zeta w)$ is
an analytic function of $w$ within $\calC_1$, for every $j\in\bfI$ such
that $c_{ij}^{(r)}\neq 0$ for some $i\in\bfI$.
\end{itemize}
Note that, for $\zeta$ small, the cut--set $\zeta^{-1}q^{-l-r}\sfY(\Psi_{j,0}
^+\Psi_j(w))$ of $H_j^-(q^{l+r}\zeta w)$ can be made to avoid the contour $\calC_1$.
In particular, the \rhs of the equation above defines a holomorphic function
of $\zeta$ in a neighborhood of $\zeta=0$.

The following is the counterpart for $\qloop$ of Theorem \ref{thm: coeff}.

\begin{thm}\label{thm: coeff-trig}\hfill
\begin{enumerate}
\item $\Aq_{\V_1,\V_2}(\zeta)$ is a rational function of $\zeta$, which is
regular at $0$ and $\infty$, and such that $\Aq_{\V_1,\V_2}(0)=1=\Aq_
{\V_1,\V_2}(\infty)$.
\item The poles of $\Aq_{\V_1,\V_2}(\zeta)^{\pm 1}$ are contained in
$\bigcup_r \sigma(\V_2)\sigma(\V_1)^{-1}q^{-l-r}$, where $r$ ranges over the
integers such that $c_{ij}^{(r)}\neq 0$ for some $i,j\in\bfI$.
\item For any $\zeta,\zeta'$ we have 
$[\Aq_{\V_1,\V_2}(\zeta), \Aq_{\V_1,\V_2}(\zeta')]=0$.
\item For any $\V_1,\V_2,\V_3\in\Rloop$, we have
\begin{align*}
\Aq_{\V_1\Dotimes{\zeta_1}\V_2,\V_3}(\zeta_2) &= 
\Aq_{\V_1,\V_3}(\zeta_1\zeta_2)\Aq_{\V_2,\V_3}
(\zeta_2)\\
\Aq_{\V_1,\V_2\Dotimes{\zeta_2}\V_3}(\zeta_1\zeta_2) &= 
\Aq_{\V_1,\V_3}(\zeta_1\zeta_2)
\Aq_{\V_1,\V_2}(\zeta_1)
\end{align*}
\item The following shifted unitarity condition holds
\[\sigma\circ \Aq_{\V_1,\V_2}(\zeta^{-1})\circ \sigma^{-1}=
\Aq_{\V_2,\V_1}(q^{-2l}\zeta)\]
where $\sigma:\V_1\otimes \V_2\to \V_2\otimes \V_1$ is the flip of the tensor
factors.
\item For every $\alpha,\beta\in\nC$ we have
\[
\Aq_{\V_1(\alpha),\V_2(\beta)}(\zeta) = \Aq_{\V_1,\V_2}(\zeta\alpha\beta^{-1})
\]
\end{enumerate}
\end{thm}

\subsection{}

The following result is needed to prove (i) of the theorem above. The rest of 
the theorem follows from the same reasoning as Theorem \ref{thm: coeff}.

\begin{lem}\label{lem:Aq-operator-dual}
$\Aq_{\V_1,\V_2}(\zeta)$ extends to a rational function of $\zeta$.
Let $\calC_1$ be a contour enclosing the set of poles
$\spec(\V_1)\subset\nC$, and such that $0$ is outside of $\calC_1$.
Then we have the following 
\[
\Aq_{\V_1,\V_2}(\zeta) = \exp\lp
-\sum_{\begin{subarray}{c} i,j\in\bfI\\ r\in\Z\end{subarray}}
c_{ij}^{(r)} \oint_{\calC_1} \frac{dH_i^+(w)}{dw}
\otimes H_j^+(q^{l+r}\zeta w)\, dw
\rp
\]
where $H_j^+(w) = \log(\Psi_{j,0}^-\Psi_j(w))$ is defined using
Proposition \ref{pr:matrix log}, and $\zeta\in\C$ is large enough such that
$H_j^+(q^{l+r}\zeta w)$ is an analytic function of $w$
within $\calC_1$ for every $j\in\bfI$ such that
$c_{ij}^{(r)}\neq 0$ for some $i\in\bfI$.
\end{lem}

\begin{pf}
The proof of this lemma is analogous to that of Theorem \ref{thm: coeff} (i)
(see Claim $1$ in the proof of that theorem). Again we revert to a more
general set up as follows. Let $V,W$ be complex, \fd vector spaces, $A:\C
\to \End(V)$, $B:\C\to\End(W)$ rational functions such that 
\begin{itemize}
\item $A(z),B(z)$ are regular and invertible at $z=0$ and $z=\infty$. 
\item $A(\infty) \in GL(V)$ is a semisimple operator.
\item $[A(z),A(w)]=0=[B(z),B(w)]$.
\end{itemize}
Let $b_0(z) = \logo(B(0)^{-1}B(z))$ be defined according to Proposition \ref{pr:matrix log-trig}
and $b_{\infty}(z) = \log(B(\infty)^{-1}B(z))$ using Proposition \ref{pr:matrix log}. 
Let $\calC_1$ denote a contour in $\nC$ enclosing all the poles of $A(z)^{\pm 1}$ and not
enclosing $0$.
Define
\begin{align*}
X_0(\zeta) &= \exp\lp \oint_{\calC_1} A(w)^{-1}A'(w)\otimes b_0(\zeta w)\, dw\rp \\
X_{\infty}(\zeta) &= \exp\lp \oint_{\calC_1} A(w)^{-1}A'(w)\otimes b_{\infty}(\zeta w)\, dw\rp
\end{align*}
where, for $X_0$ we need to take $\zeta$ small enough so that $\zeta^{-1}\sfY(B(0)^{-1}B(w))$ is 
outside of $\calC_1$ and hence $b_0(\zeta w)$ is analytic within $\calC_1$, and for $X_{\infty}$
we need to take $\zeta$ large enough so that $\zeta^{-1}\sfX(B(\infty)^{-1}B(w))$ is outside
of $\calC_1$ and hence $b_{\infty}(\zeta w)$ is analytic within $\calC_1$.\\

We need to prove that both $X_0(\zeta)$ and $X_{\infty}(\zeta)$ extend to the same rational function
of $\zeta$, taking values in $\End(V\otimes W)$. 
For this we consider the Jordan decomposition $A(z)=A_S(z)A_U(z)$. By \cite[Lemma 4.12]{sachin-valerio-2},
$A_S(z)$ and $A_U(z)$ are again rational funtions of $z$. Furthermore, $A_U(\infty) = 1$ by our assumption
that $A(\infty)$ is semisimple.
Since logarithmic derivative of $A(z)$ splits the two additively,
we can treat the semisimple and unipotent cases separately, analogous to the proof of Claim $1$ in
Theorem \ref{thm: coeff} given in Section \ref{ssec: pf-A1}.\\

The semisimple case reduces to the scalar case, \ie when $V$ is one--dimensional and
\[
A(z) = A(\infty) \prod_j \frac{z-\alpha_j}{z-\beta_j}
\]
for some $\alpha_j,\beta_j\in\nC$. Following the computation given in Section
\ref{ssec: pf-A1}, we get
\[
X_0(\zeta) = \prod_j 1\otimes B(\zeta\alpha_j)B(\zeta\beta_j)^{-1} = X_{\infty}(\zeta)
\]
Now assuming $A(z)$ is unipotent and $A(\infty)=1$, we get that $\log(A(z))$ is again
a rational function of $z$, vanishing at $z=\infty$. Decomposing it into partial fractions yields
\[
\log(A(z)) = \sum_{\begin{subarray}{c} j\in J\\ n\in\N\end{subarray}}
\frac{N_{j,n}}{(z-\alpha_j)^{n+1}}
\]
where $J$ is a finite indexing set and $\alpha_j\in\nC$. We obtain
\begin{align*}
X_0(\zeta) &= \exp\lp \sum_{\begin{subarray}{c} j\in J\\ n\in\N\end{subarray}}
-(n+1)N_{j,n} \otimes \left.\frac{\partial_w^{n+1}}{(n+1)!} b_0(\zeta w) \right|_{w=\alpha_j}\rp \\
X_{\infty}(\zeta) &= \exp\lp \sum_{\begin{subarray}{c} j\in J\\ n\in\N\end{subarray}}
-(n+1)N_{j,n} \otimes \left.\frac{\partial_w^{n+1}}{(n+1)!} b_{\infty}(\zeta w) \right|_{w=\alpha_j}\rp 
\end{align*}
which are both rational functions, since the $N_{j,n}$ are nilpotent and pairwise commute. As rational
functions, the two are the same since $b_0'(w) = b_{\infty}'(w) = B(w)^{-1}B'(w)$.
\end{pf}

\subsection{Commutation relation with $\Aq_{\V_1,\V_2}(\zeta)$}
\label{ssec: op-B-comm}
Let $\calC\subset\C$ be a contour,
and $a_\ell:\C\to\End(\V_\ell)$, $\ell=1,2$ two meromorphic functions
which are analytic within $\calC$ and commute with the operators
$\{\Psi^{\pm}_{i,\pm r}\}_{i\in\bfI,r\in\N}$. For any $k\in\bfI$, define operators
$X_k^{\pm,\ell}\in\End(\V_1\otimes \V_2)$ by
\[X_k^{\pm, 1} = \oint_{\calC}
a_1(w)\X_k^{\pm}(w)\otimes a_2(w)\, dw
\quad\text{and}\quad
X_k^{\pm, 2} = \oint_{\calC}
a_1(w)\otimes a_2(w)\X_k^{\pm}(w)\, dw \label{eq: Xk-trig}
\]

\begin{prop}\label{prop: B-comm}
The following commutation relations hold
\begin{align*}
\Ad(\Aq_{\V_1,\V_2}(\zeta))X_k^{\pm,1} &=
\oint_{\calC}
a_1(w)\X_k^{\pm}(w)\otimes a_2(w)
\Psi_k(q^{2l}\zeta w)^{\pm 1}\Psi_k(\zeta w)^{\mp 1}\, dw\\
\Ad(\Aq_{\V_1,\V_2}(\zeta))X_k^{\pm,2} &=
\oint_{\calC}
a_1(w)\Psi_k(\zeta^{-1}w)^{\pm 1}\Psi_k(q^{-2l}\zeta^{-1}w)^{\mp 1}
\otimes a_2(w)\X_k^{\pm}(w)\, dw
\end{align*}
\end{prop}

The proof of this proposition is identical to that of Proposition \ref{prop: A-comm},
except that the following version of relation \eqref{eq: imp-lem} is needed.
For each $i,k\in\bfI$, 
\begin{align*}
\left[ \Psi_i(z)^{-1}\Psi_i'(z),\X^{\pm}_k(w)\right] &=
\lp \frac{1}{z-q^{\mp b_{ik}}w} - \frac{1}{z-q^{\pm b_{ik}}w}\rp\X_k^{\pm}(w)\\
&+\frac{wq^{\pm b_{ik}}}{z(z-q^{\pm b_{ik}}w)} \X_k^{\pm}(q^{\mp b_{ik}}z)
-\frac{w}{z(q^{\pm b_{ik}}z-w)}\X_k^{\pm}(q^{\pm b_{ik}}z)
\end{align*}

One can derive this relation easily from (\QL3) of Proposition \ref{pr:qloop fields}
following the computation given in the proof of Lemma \ref{lem: imp-lem}.

\subsection{Regular $q$--difference equations}\label{ssec: qdiffSauloy}

We review below the existence and uniqueness of solutions of regular
$q$--difference equations. Let $p\in\nC$ be such that $|p|\neq 1$, $W$
a complex, \fd vector space, and consider the difference equation
\begin{equation}\label{eq:p-difference}
F(pz)= B(z)F(z)
\end{equation}
with values in $\End(W)$. Here, $B(z)$ is a meromorphic, $\End(W)
$--valued function. We shall assume that the equation is {\it regular},
that is that $B$ is holomorphic near $z=0$, and such that $B(0)=1$.
The following result is well--known (see, \eg\cite[\S 1.2.2]{sauloy}).

\begin{lem}
There is a unique formal series $F(z)\in \End(W)[[z]]$ which satisfies
\eqref{eq:p-difference} and $F(0)=1$. Moreover, $F$ converges
near $z=0$ to a meromorphic function defined on $\C$. Any meromorphic
solution $G(z)$ of \eqref{eq:p-difference} which is holomorphic
in a neighborhood of $z=0$ is of the form
$F(z)C$ where $C = G(0) \in \End(W)$ is a constant matrix.
\end{lem}

Let us remark that the existence of the formal series is automatic, since
the equation \eqref{eq:p-difference} is equivalent to the recurrence relation
$(p^n-1)F_n=\sum_{m=0}^{n-1}B_{n-m}F_m$, where $F=\sum_{n\geq 0}F_nz^n$
and $B=\sum_{n\geq
0}B_nz^n$, with $F_0=1=B_0$. The convergence of $F$ is proved in
\cite[\S 1.2.2 Lemme 1]{sauloy}. The uniqueness is also clear, since
the ratio $F(z)^{-1}G(z)$ is a holomorphic function on the elliptic
curve $\nC/p^{\Z}$, and hence a constant.\\

 One gets a similar assertion
if $z=0$ is changed to $z=\infty$ and one considers formal series in
$z^{-1}$.

\subsection{The abelian $R$--matrix of $\qloop$}\label{ssec: R0-qla}

Assume now that $|q|\neq 1$. Let $\V_1,\V_2\in\Rloop$, and
let $\Aq_{\V_1,\V_2}(\zeta)$ be the operator defined in \ref{ssec: op-B}.
Consider the $q$--difference equation
\[\ol{\Rq}_{\V_1,\V_2}(q^{2l}\zeta) = \Aq_{\V_1,\V_2}(\zeta)
\ol{\Rq}_{\V_1,\V_2}(\zeta)\]
This equation is regular at $0$ and $\infty$ since $\Aq_{\V_1,
\V_2}(0)=1=\Aq_{\V_1,\V_2}(\infty)$. By Lemma \ref{ssec: qdiffSauloy},
it admits two unique formal solutions $\ol{\Rq}^{\pm}(\zeta)$
near $q^{\pm\infty}$, 
which are normalized by
\[\ol{\Rq}^+(q^\infty) = 1 = \ol{\Rq}^-(q^{-\infty})\]
These solutions converge in a neighborhood of $q^{\pm \infty}$,
and extend to meromorphic functions on the entire complex plane 
which are given by the products
\[\ol{\Rq}^+(\zeta)=\stackrel{\longrightarrow}{\prod}_{n\geq 0}\Aq_{\V_1,\V_2}(q^{2ln}\zeta)^{-1}
\aand
 \ol{\Rq}^-(\zeta)=\stackrel{\longrightarrow}{\prod}_{n\geq 1}\Aq_{\V_1,\V_2}(q^{-2ln}\zeta)\]

Set
\[
\Rq^{0,\pm}_{\V_1,\V_2}(\zeta)=
\left\{
\begin{array}{lcl}
q^{\mp\Omega_{\h}}\ol{\Rq}^\pm(\zeta) & \text{ if } & |q|<1 \\
q^{\pm\Omega_{\h}}\ol{\Rq}^\pm(\zeta) & \text{ if } & |q|>1
\end{array}
\right.
\]
By uniqueness, the evaluation on $\V_1\otimes\V_2$ of the operator
$\Rq^0(\zeta)$ given by \eqref{eq:R0qla-formal}, is the expansion at
$\zeta=0$ of $\Rq^{0,\veps}_{\V_1,\V_2}(\zeta)$, where $\veps\in\{\pm\}$
is such that $q^{\veps\infty}=0$.  

The following is the analog of Theorem \ref{thm: R0} for $\qloop$.

\begin{thm}\label{thm: R0-qla}
The operators $\Rqio_{\V_1,\V_2}(\zeta)$ have the following properties
\begin{enumerate}
\item The map
\[\sigma\circ \Rqio_{\V_1,\V_2}(\zeta) :
\V_1(\zeta)\Dotimes{1}\V_2 \rightarrow \V_2\Dotimes{1} \V_1(\zeta)\]
where $\sigma$ is the flip of tensor factors, is a morphism of $\qloop
$--modules, which is natural in $\V_1$ and $\V_2$.
\item For any $\V_1,\V_2,\V_3\in\Rloop$ we have
\begin{align*}
\Rqio_{\V_1\Dotimes{\zeta_1}\V_2,\V_3}(\zeta_2) &= \Rqio_{\V_1,\V_3}(\zeta_1\zeta_2)
\Rqio_{\V_2,\V_3}(\zeta_2)\\
\Rqio_{\V_1,\V_2\Dotimes{\zeta_2}\V_3}(\zeta_1\zeta_2) &= \Rqio_{\V_1,\V_3}
(\zeta_1\zeta_2)
\Rqio_{\V_1,\V_2}(\zeta_1)
\end{align*}
\item The following unitary condition holds
\[\sigma\circ \Rqio_{\V_1,\V_2}(\zeta^{-1})\circ \sigma^{-1}=
\Rqoi_{\V_2,\V_1}(\zeta)^{-1}\]
\item For $\alpha,\beta\in\nC$, we have
\[\Rqio_{\V_1(\alpha),\V_2(\beta)}(\zeta) 
= \Rqio_{\V_1,\V_2}(\zeta\alpha\beta^{-1})\]
\item For any $\zeta,\zeta'$, 
\[[\Rqio_{\V_1,\V_2}(\zeta),\Rqio_{\V_1,\V_2}(\zeta')]=
0=[\Rqio_{\V_1,\V_2}(\zeta),\Rqoi_{\V_1,\V_2}(\zeta')]\] 
\item $\Rqio_{\V_1,\V_2}(\zeta)$ is holomorphic near $q^{\pm\infty}$,
and
\[\Rqio_{\V_1,\V_2}(q^{\pm \infty})=
\left\{\begin{array}{ll}
q^{-\Omega_\h}&\text{if $q^{\pm\infty}=0$}\\
q^{\Omega_\h}&\text{if $q^{\pm\infty}=\infty$}
\end{array}\right.\]
\item The poles of $\Rqi_{\V_1,\V_2}(\zeta)^{\pm 1}$ and $\Rqo_{\V_1,
\V_2}(\zeta)^{\pm 1}$ are contained in
\[\sigma(\V_2)\sigma(\V_1)^{-1}q^{-l-r}q^{-2l\Z_{\geq 0}}
\quad\text{and}\quad
\sigma(\V_2)\sigma(\V_1)^{-1}q^{-l-r}q^{2l\Z_{>0}}\]
respectively, where $r$ ranges over the integers such that $c_{ij}
^{(r)}\neq 0$ for some $i,j\in\bfI$.
\end{enumerate}
\end{thm}

\section{Kohno--Drinfeld theorem for abelian, additive $q$KZ equations}\label{se:KD}

In this section, we prove that, when $\Im\hbar\neq 0$, the monodromy of
the additive $q$KZ equations on $n$ points defined by the commutative
$R$--matrix of the Yangian is given by the commutative $R$--matrix of
the quantum loop algebra. The general case is treated in \ref{ss:KD},
and follows from the $n=2$ case which is treated in \ref{ss:n=2}--\ref
{ss:KD2}. In turn, the latter rests on relating the coefficient matrices
of the difference equations whose solutions are the commutative
$R$--matrix of $\Yhg$ and $\qloop$ respectively, which is done
in \ref{ss:A-B} below.

\subsection{}\label{ss:A-B}

Let $V_1,V_2$ be two \fd representations of $\Yhg$, $\A_{V_1,V_2}(s)$
the meromorphic $GL(V_1\otimes V_2)$--valued function constructed in
\ref{ssec: op-A}, and consider the difference equation\footnote{Note that
\eqref{eq:A 1 ADE} differs from the difference equation considered in \ref
{ssec: R0-prop} since its step is $1$, not $l\hbar$.}
\begin{equation}\label{eq:A 1 ADE}
f(s+1) = \opA{V_1}{V_2}(s)f(s)
\end{equation}

Assume further that $V_1,V_2$ are non--congruent, let $\V_{\ell}=\Fh{}
(V_{\ell})$ be the representations of $\qloop$ obtained by using the
functor $\Fh{}$ of Section \ref{sec: functor}, and $\Aq_{\V_1,\V_2}(\zeta)
\in GL(\V_1\otimes\V_2)$ the operator constructed in \ref{ssec: op-B}.

\begin{prop}\label{prop: A-B}
The operator $\Aq_{\V_1,\V_2}(\zeta)$ is the monodromy of the difference
equation \eqref{eq:A 1 ADE}. That is,
\[
\Aq_{\V_1,\V_2}(\zeta) = 
\left.\prod_{m\in\Z} \A_{V_1,V_2}(s+m)\right|_{\zeta=e^{2\pi\iota s}}
\]
\end{prop}
\begin{pf}
For the purposes of the proof, we restrict ourselves to a typical factor
in the definition of $\A_{V_1,V_2}(s)$. That is, fix $i,j\in\bfI$ and define
\[
\A_{V_1,V_2}^{ij}(s) = \exp\lp\oint_{\calC_1} t_i'(v)\otimes
t_j(v+s)\, dv\rp
\]
where $\calC_1$ encloses the poles of $\xi_i(v)^{\pm 1}$ on $V_1$,
and $s$ is such that $t_j(v+s)$ is analytic within $\calC_1$. Since $
V_1$ is non--congruent, we may further assume that 
no two distinct points in the interior of $\calC_1$ or on $\calC_1$
are congruent modulo $\Z$.

By Theorem \ref{thm: coeff}, $\A_{V_1,V_2}^{ij}(s)$ is a rational
function of the form $1 + O(s^{-2})$. The corresponding monodromy
matrix $M(\zeta)$ is a rational function of $\zeta=\exp(2\pi\iota s)$
which is given by
\[M(\zeta)=
\left.\prod_{m\in\Z}\A_{V_1,V_2}^{ij}(s+m)\right|_{\zeta=e^{2\pi\iota s}}
= \lim_{N\to \infty} \lp \left.\prod_{m=-N}^{N}\A_{V_1,V_2}^{ij}(s+m)\right|_{\zeta=e^{2\pi\iota s}}
\rp
\]
The corresponding factor of $\Aq_{\V_1,\V_2}(\zeta)$ is given by
\[\Aq_{\V_1,\V_2}^{ij}(\zeta)=
\exp\lp\oint_{\wt{\calC}_1}
\Psi_i(w)^{-1}\Psi_i(w)'
\otimes H_j^-(w\zeta)\, dw\rp\]
where $\wt{\calC}_1=\exp(2\pi\iota\calC_1)$, and $H_j^-(w)=\logo
(\Psi^+_{j,0}\Psi_j(w))$ is given by Proposition \ref{pr:matrix log-trig}.
Note that $\wt{\calC_1}$ is again a Jordan curve because of the assumptions
imposed on $\calC_1$.

We wish to show that $M(\zeta) = \Aq_{\V_1,\V_2}^{ij}(\zeta)$. Since
both sides are rational functions of $\zeta$, it suffices to prove this
for $\zeta$ near $0$, that is $\Im(s)\gg 0$. Now
\[M(\zeta)=\lim_{N\to\infty}\exp\left(
\sum_{m=-N}^N 
\oint_{\calC_1} t_i'(v)\otimes t_j(v+s+m)\, dv\right)
\]
Since $t_j(v)=\hbar\xi_{j,0}v^{-1}+O(v^{-2})$, the sum $\sum_{m=-N}^N 
t_j(u+m)$ converges uniformly on compact subsets of $\{|\Im u|> R\}$
for $R$ large enough. 
To see this, we note that $t_j(v)$, as defined using Proposition \ref{pr:matrix log},
is a holomorphic function in a neighborhood of $v=\infty$. Its Taylor series
$t_j(v) = \hbar\sum_{r\geq 0} t_{j,r}u^{-r-1}$ therefore converges uniformly
for $|u|>R$, for some $R>0$. Each partial sum
$f_N(u) = \sum_{m=-N}^N t_j(u+m)$ is a holomorphic function on $\Omega
= \{|\Im(u)| > R\}$ since the set of poles of $f_N$ is contained
in the shifts of the closed disc $D_0(R) := \{|u|\leq R\}$ by integers $m\in \{-N,\ldots,N\}$,
and hence does not intersect $\Omega$.\\

Given a compact subset $K\subset \Omega$, let $R_2>R_1>R$ be such
that $R_1<|u|<R_2$ for each $u\in K$. Take $M>N>R_1+R_2$ and let us
find an upper bound on $\|f_{N+1}(u)-f_M(u)\|$, which goes to $0$ as $N$
goes to $\infty$, uniformly for each $u\in K$. By definition of the radius of
convergence there exists $P>0$ such that $\|\hbar t_{j,r}\| \leq PR_1^{r+1}$
for each $r\geq 0$

{\em First order term.} For each $u\in K$ and $m\geq N$ we have
$|u^2-m^2| > m^2 - R_2^2$. Therefore
\begin{align*}
\left| \sum_{m=N+1}^{M} \frac{1}{u+m} + \frac{1}{u-m}\right| &\leq
\sum_{m=N+1}^{\infty} \frac{2|u|}{|u^2-m^2|} < \sum_{m=N+1}^{\infty} \frac{2R_2}{m^2-R_2^2} \\
&\leq \int_N^{\infty} \frac{2R_2}{x^2-R_2^2}\, dx = \ln\lp\frac{N+R_2}{N-R_2}\rp
\end{align*}

{\em Higher order terms.} Again we have $|u\pm m| > m-R_2$ for each $u\in K$.
Thus, for each $r\geq 1$ we get
\begin{align*}
\left|\sum_{m=N+1}^{M} \frac{1}{(u\pm m)^{r+1}}\right| &\leq
\sum_{m=N+1}^{\infty} \frac{1}{(m-R_2)^{r+1}} \leq \int_{N}^{\infty} \frac{1}{(x-R_2)^{r+1}}\, dx \\
&= \frac{1}{r(N-R_2)^r}
\end{align*}

Hence, for any $M>N$ we obtain the following bound:
\begin{align*}
\|f_{N+1}(u) - f_M(u)\| &< 2PR_1\lp \ln\lp\frac{N+R_2}{N-R_2}\rp
+ \sum_{r\geq 1} \frac{1}{r} \lp\frac{R_1}{N-R_2}\rp^r \rp \\
&= 2PR_1\ln\lp\frac{N+R_2}{N-R_1-R_2}\rp
\end{align*}
Thus given $\epsilon > 0$ we can choose $N$ large enough so that the above bound is less than
$\epsilon$ uniformly for each $u\in K$, as claimed.\\

Now the exponential of $\lim_{N\to\infty} \sum_{m=-N}^N t_j(u+m)$ is 
$\prod_m\xi_j(u+m)=\Psi_j(e^
{2\pi\iota u})$ (see \S \ref{ssec: functor}). By the uniqueness of $\logo$
this implies that, for $\Im s\gg 0$,
\[
\lim_{N\to\infty} \sum_{m=-N}^N t_j(v+s+m) =  
H_j^-(\zeta e^{2\pi\iota v}) + \log\lp e^{-\pi\iota\hbar\xi_{j,0}}\rp
\]
To see this we observe that, for a fixed $v$, both sides of the equation
above have the same exponential and the same value at $\zeta = 0$, or
equivalently $\Im s \to \infty$, in the domain $\{\Im(s)>R-\Im(v)\}$.
For $v$ ranging over a compact set (\eg interior of $\calC_1$
and $\calC_1$ included, which is needed below) we can take $\Im s\gg 0$ so that 
for each $v$ in this compact set, the equation holds. Thus we get
\[M(\zeta)=\exp\left(
\oint_{\calC_1} t_i'(v)\otimes \lp H_j^-(\zeta e^{2\pi\iota v}) + 
\log\lp e^{-\pi\iota\hbar\xi_{j,0}}\rp\rp\, dv\right)
\]

Since $t_i'(v) = O(v^{-2})$, we get $\oint t_i'(v)\otimes \log(e^{-\pi\iota\hbar\xi_{j,0}})\, dv=0$,
which implies that
\[M(\zeta)=
\exp\left( \oint_{\calC_1} \xi_i(v)^{-1}\xi_i(v)'\otimes H_j^-(\zeta e^{2\pi\iota v})\, dv\right)
\]
Noting that, by \eqref{eq:3 terms}
\[\Psi_i(\Exp{v})^{-1}\frac{d\Psi_i(\Exp{v})}{dv}=
g_i^+(v)^{-1}g_i^+(v)'+
\xi_i(v)^{-1}\xi_i(v)'+g_i^-(v)^{-1}g_i^-(v)'\]
and that $g_i^{\pm}(v)$ are analytic and invertible within $\calC_1$ by the
non--congruence assumption, so that 
\[\oint_{\calC_1} g_i^\pm(v)^{-1}g_i^\pm(v)'
\otimes 
H_j^-(\zeta e^{2\pi\iota v})\, dv = 0\]
we get
\begin{align*}
M(\zeta)
&= \exp\left(\oint_{\calC_1}
\Psi_i(\Exp{v})^{-1}\frac{d\Psi_i(\Exp{v})}{dv}
\otimes H_j^-(\zeta e^{2\pi\iota v})\, dv\right)\\
&= \exp\left(\oint_{\wt{\calC}_1}
\Psi_i(w)^{-1}\frac{d\Psi_i(w)}{dw}
\otimes H_j^-(\zeta w)\, dw\right)
\end{align*}
as claimed.
\end{pf}

\subsection{The (reduced) $q$KZ equations on $n=2$ points}\label{ss:n=2}

Assume henceforth that $\Im\hbar\neq 0$. 
Fix $\veps\in\{\pm\}$, let $V_1,V_2\in\Rep_{\ffd}(\Yhg)$, and consider
the abelian $q$KZ equation
\[f(s+1)=\Reps_{V_1,V_2}(s)f(s)\]
with values in $\End(V_1\otimes V_2)$.

By Proposition \ref{pr:left right}, this equation admits both right and
left canonical solutions $\Phi^\veps_\pm(s)$. The corresponding
connection matrix is given by
\[\begin{split}
\calS_{V_1,V_2}^\veps(s)
&=
\Phi^\veps_+(s)^{-1}\Phi^\veps_-(s)\\
&=
\lim_{N\to\infty}\Reps_{V_1,V_2}(s+N)\cdots
\Reps_{V_1,V_2}(s)\cdots\Reps_{V_1,V_2}(s-N)
\end{split}\]
and is a meromorphic function of $\zeta=\Exp{s}$ which admits
a limit as $\Im s\to\pm\infty$, depending on whether $\Im(\veps
\hbar)\gtrless 0$. In particular, $\calS^\veps_{V_1,V_2}(\zeta)$
is regular at $\zeta=q^{\veps\infty}$.

\begin{lem}\label{lem: S-limits}
\[\calS_{V_1,V_2}^\veps(q^{\veps\infty})=
\left\{\begin{array}{ll}
q^{-\Omega_\h}		&\text{if $q^{\veps\infty}=0$}\\
q^{\Omega_\h}	&\text{if $q^{\veps\infty}=\infty$}
\end{array}\right.\]
\end{lem}
\begin{pf}
Let us assume $\Im(\hbar)>0$ and $\veps = +$,
for definiteness. Then, by Proposition \ref{pr:left right},
$\Phi^+_+(s)$ has the asymptotic expansion of the form $(1+O(s^{-1}))s^{\hbar\Omega_{\h}}$
in any right half--plane, while $\Phi^+_-(s)\sim (1+O(s^{-1}))(-s)^{\hbar\Omega_{\h}}$ only in an obtuse sector
shown in Figure \ref{fig: zones}. Thus we can find a common domain
for both, where the limit $\Im(s)\to \infty$ can be taken. Now we have
\begin{align*}
\calS^+_{V_1,V_2}(0) &= \lim_{\Im(s)\to\infty} \lp \Phi^+_+(s)\rp^{-1}\Phi^+_-(s) 
= \lim_{\Im(s)\to\infty} s^{-\hbar\Omega_{\h}}(-s)^{\hbar\Omega_{\h}} \\
&= \lim_{\Im(s)\to\infty} e^{\hbar\Omega_{\h}(\ln(-s)-\ln(s))} 
= e^{-\pi\iota\hbar\Omega_{\h}}
\end{align*}
\end{pf}

\begin{rem}
Note that in the proof above, there is no common domain where both $\Phi^+
_{\pm}(s)$ admit the claimed asymptotic expansions and $\Im(s)$ can go to $-\infty$. Consequently, the
computation above cannot be carried out for $\calS^+_{V_1,V_2}(\infty)$. This is
in contrast with the computation of the monodromy of an additive
difference equation when the coefficient matrix is rational, given,
for example, in \cite[Prop. 4.8]{sachin-valerio-2}.
\end{rem}

\subsection{Kohno--Drinfeld theorem for abelian $q$KZ equations on 2 points.}\label{ss:KD2}


The following equates the monodromy of the abelian $q$KZ
equations with the commutative $R$--matrix of $\qloop$
constructed in \ref{ssec: R0-qla}.

\begin{thm}\label{th: n=2KD}
If $V_1,V_2$ are non--congruent, $\V_\ell=\Fh{}(V_\ell)$ are
the corresponding representations of $\qloop$, and $\Rq_{\V
_1,\V_2}^{0,\veps}(\zeta)$ is the commutative $R$--matrix of
$\qloop$, then
\[\calS_{V_1,V_2}^\veps(\zeta)=\Rq_{\V_1,\V_2}^{0,\veps}(\zeta)\]
\end{thm}

\begin{pf}
Let $\A_\pm(s)$ be the right and left fundamental solutions of
the difference equation $f(s+1)=\A_{V_1,V_2}(s)f(s)$ considered
in \ref{ss:A-B}. We claim that
\begin{equation}\label{eq:intermediate}
\Phi^\veps_\pm(s+l\hbar)\Phi^\veps_\pm(s)^{-1}=\A_\pm(s)
\end{equation}

Assuming this for now, we see that $\calS^\veps_{V_1,V_2}(\zeta)$
and $\Rq^{0,\veps}_{\V_1,\V_2}(\zeta)$ satisfy the same $q$--difference
equation. Indeed,
\[\begin{split}
\calS_{V_1,V_2}^\veps(q^{2l}\zeta)\calS_{V_1,V_2}^\veps(\zeta)^{-1}
&=
\Phi^\veps_+(s+l\hbar)^{-1}\Phi^\veps_-(s+l\hbar)
\Phi^\veps_-(s)^{-1}\Phi^\veps_+(s)\\
&=
\A_+(s)^{-1}\A_-(s)\\
&=
\Aq_{\V_1,\V_2}(\zeta)\\
&=
\Rq^{0,\veps}_{\V_1,\V_2}(q^{2l}\zeta)\Rq^{0,\veps}_{\V_1,\V_2}(\zeta)^{-1}
\end{split}\]
where the third equality follows by Proposition \ref{prop: A-B},
and the last one by definition of $\Rq^{0,\veps}_{\V_1,\V_2}$.
Note that the reordering of factors in the calculation above is 
permissible since all the meromorphic functions involved take
values in a commutative subalgebra of $\End(V_1\otimes V_2)$.
Since both $\calS_{V_1,V_2}^\veps$ and $\Rq^{0,\veps}_{\V_1,
\V_2}$ are holomorphic near, and have the same value at $
\zeta=q^{\veps\infty}$, they are equal.

Returning to the claim, let $L^\veps_\pm(s)$ denote the \lhs
of \eqref{eq:intermediate}. Then,
\[\begin{split}
L^\veps_\pm(s+1)L^\veps_\pm(s)^{-1}
&=
\Phi^\veps_\pm(s+l\hbar+1)\Phi^\veps_\pm(s+1)^{-1}
\Phi^\veps_\pm(s)\Phi^\veps_\pm(s+l\hbar)^{-1}\\
&=
\Reps_{V_1,V_2}(s+l\hbar)\Reps_{V_1,V_2}(s)^{-1}\\
&=
\A_{V_1,V_2}(s)
\end{split}\]
Thus, $L^\veps_\pm(s)$ and $\A_\pm(s)$ satifsy the same
difference equation. Since they also have the same asymptotics
as $s\to\infty$ by Proposition \ref{pr:left right}, it follows that they
are equal.
\end{pf}

\begin{rem}\label{rem:KDimpliesBraiding}
The monodromy $\calS_{V_1,V_2}^\veps(\zeta)$ may be written
in terms of the tensor structures ${\mathcal J}^{\pm}_{V_1,V_2}(s)$
constructed in \ref{ssec: twist-defn} as
\[\calS_{V_1,V_2}^\veps(\zeta)=
\Jeps_{V_1,V_2}(s)\Reps_{V_1,V_2}(s)\left(\Jmeps_{V_2,V_1}(-s)\right)_{21}^{-1}\]
where we used the unitarity constraint (iii) of Theorem \ref{thm: R0}.
We can rearrange the factors of the triple product in the \rhs above, again using
the fact that all relevant meromorphic functions take values in
a commutative subalgebra of $\End(V_1\otimes V_2)$.
This, and Theorem \ref{th: n=2KD} imply the following equation
\[
\sigma\circ \Rq^{0,\veps}_{\V_1,\V_2}(\zeta) = 
\Jmeps_{V_2,V_1}(-s)^{-1}\circ\lp \sigma\circ \Reps_{V_1,V_2}(s)\rp\circ
\Jeps_{V_1,V_2}(s)
\]
Here $\V_{\ell} = \Fh{}(V_{\ell})$ for $\ell = 1,2$. This equation
implies the commutativity of the following diagram, which means that
the tensor structures ${\mathcal J}^{\pm}_{V_1,V_2}(s)$ are compatible
with the meromorphic braidings on $\Ryang$ and $\Rloop$
given by $\Reps_{V_1,V_2}(s)$ and $\Rq_{\V_1,\V_2}^{0,\veps}(\zeta)$.
\[\xymatrix{
\Fh{}(V_1)\otimes_{\zeta}\Fh{}(V_2) \ar@{=}[d]
\ar[rrr]^{\Jeps_{V_1,V_2}(s)} &&& \Fh{}(V_1\otimes_s V_2)\ar@{=}[d] \\
\Fh{}(V_1)(\zeta)\otimes_1 \Fh{}(V_2) 
\ar[dd]_{\sigma\circ\Rq^{0,\veps}_{\V_1,\V_2}(\zeta)}
&&& \Fh{}(V_1(s)\otimes_0 V_2) \ar[dd]^{\Fh{}(\sigma\circ\Reps_{V_1,V_2}(s))} \\
&&& \\
\Fh{}(V_2)\otimes_1\Fh{}(V_1)(\zeta) \ar@{=}[d]
&&& \Fh{}(V_2\otimes_0 V_1(s))\ar@{=}[d]\\
\lp\Fh{}(V_2)\otimes_{\zeta^{-1}}\Fh{}(V_1)\rp(\zeta)
\ar[rrr]_{\Jmeps_{V_2,V_1}(-s)}
&&&
\lp\Fh{}(V_2\otimes_{-s}V_1)\rp(\zeta)
}
\]
\end{rem}

\subsection{The abelian $q$KZ equations}\label{ssec: qKZ-recall}

Fix $\veps\in\{\pm \}$ and $n\geq 2$, and let $V_1,\ldots,V_n\in\Ryang$.

The following system of difference equations for a meromorphic function
of $n$ variables $\Phi:\C^n\to\End(V_1\otimes\cdots\otimes V_n)$ is an
abelian version of the $q$KZ equations \cite{frenkel-reshetikhin,smirnov}
\begin{equation}\label{eq: qKZ-n}
\Phi(\ul{s}+e_i) = A_i(\ul{s})\Phi(\ul{s})
\end{equation}
where $\ul{s} = (s_1,\ldots,s_n)$, $\{e_i\}_{i=1}^n$ is the standard
basis of $\C^n$, and
\begin{multline*}
A_i(\ul{s}) = \Reps_{i-1,i}(s_{i-1}-s_i-1)^{-1}\cdots \Reps_{1,i}(s_1-s_i-1)^{-1}\\
\cdot\Reps_{i,n}(s_i-s_n)\cdots \Reps_{i,i+1}(s_i-s_{i+1})
\end{multline*}
with $\Reps_{i,j} = \Reps_{V_i,V_j}$.

The above system is integrable, that is it satisfies
\[A_i(\uls+e_j)A_j(\uls)=A_j(\uls+e_i)A_i(\uls)\]

\subsection{Canonical fundamental solutions}

The above system admits a set of canonical fundamental solutions
which are parametrised by permutations $\sigma\in\Sym_n$, and
correspond to the right/left solutions in the case $n=2$.

To describe them, let $\Sigma^{\ve}_{\pm,ij}\subset\C^n$ denote
the asymptotic zones given in Proposition \ref{pr:left right} with $s
=s_i-s_j$, where $1\leq i\neq j\leq n$. Thus, 
\[\Sigma^\veps_{\pm,ij}=
\{\uls\in\C^n|\pm\Re(s_i-s_j)\gg 0\quad\text{and}\quad\pm\Re((s_i-s_j)/n)\gg 0\}\]
where $n\in\C^\times$ is perpendicular to $\hbar$ and such
that $\Re(n)\geq 0$, and the second condition in the definition of
$\Sigma^\veps_{\pm,ij}$ is required only if $\pm\Re(\veps\hbar)<0$.

For a permutation $\sigma\in\Sym_n$, set
\[C^\pm(\sigma)=\{i<j|\,\sigma^{-1}(i)\lessgtr\sigma^{-1}(j)\}\]
and define $\Sigma^{\ve}(\sigma)\in\C^n$ by
\[\Sigma^{\ve}(\sigma) =
\bigcap_{(i,j)\in C^+(\sigma)} \Sigma^{\ve}_{+,ij}
\cap
\bigcap_{(i,j)\in C^-(\sigma)} \Sigma^{\ve}_{-,ij}\]

\begin{prop}\label{pr:n var can}
For any $\sigma\in\Sym_n$, the equation \eqref{eq: qKZ-n} admits a
fundamental solution $\Phi^{\ve}_{\sigma}$ which is uniquely determined
by the following requirements
\begin{enumerate}
\item $\Phi^{\ve}_{\sigma}$ is holomorphic and invertible in $\Sigma^{\ve}(\sigma)$.
\item $\Phi^{\ve}_{\sigma}$ has an asymptotic expansion of the form
\[
\Phi^{\ve}_{\sigma}(\ul{s}) \sim (1 + o(1)) \prod_{(i,j)\in C^+(\sigma)}
(s_i-s_j)^{\hbar\Omega_{\h}} \prod_{(i,j)\in C^-(\sigma)} (s_j-s_i)^{\hbar\Omega_{\h}}
\]
for $\uls\in\Sigma^\veps(\sigma)$, with $s_i-s_j\to\infty$ for any $i\neq j$.
\end{enumerate}
\end{prop}
\begin{pf}
The solution $\Phi^{\ve}_{\sigma}$ is constructed as follows. For each $i<j$,
let $\Phi^{\ve}_{\pm, ij}$ be the right and left canonical solutions of the abelian
$q$KZ equation $\Phi_{ij}(s+1) = \Reps_{i,j}(s)\Phi_{ij}(s)$ given in Proposition
\ref{pr:left right}. Then,
\[\Phi^{\ve}_{\sigma}(\ul{s}) = \prod_{(i,j)\in C^+(\sigma)} \Phi^{\ve}_+(s_i-s_j)
\prod_{(i,j)\in C^-(\sigma)} \Phi^{\ve}_-(s_i-s_j)\]

We now prove the uniqueness (see, \eg \cite[\S 4.3]{sachin-valerio-2} for the one
variable case). The ratio $\Xi^\ve_\sigma=(\Phi^{\ve}_{\sigma})^{-1}\Psi^
{\ve}_{\sigma}$ of two solutions is holomorphic for $\uls\in\Sigma^\veps(\sigma)$,
and periodic under the lattice $\Z^n\subset\C^n$. It therefore descends to a
holomorphic function on the torus $T=\C^n/\Z^n=(\nC)^n$. 
We claim that $\Xi^\ve_\sigma(\ul{\zeta})=1$ for any $\ul{\zeta}\in(\nC)^n$.
Note that $\Xi^\ve_\sigma(\ul{\zeta}) = \Phi^{\ve}_{\sigma}(\uls)^{-1}\Psi^{\ve}_{\sigma}(\uls)$
for any $\uls \in \Sigma^{\ve}_{\sigma}$ such that $\zeta_j = e^{2\pi\iota s_j}$
for every $j$. By definition of the asymptotic zone $\Sigma^{\ve}_{\sigma}$,
we can find a sequence of points $\{\uls^{(1)}, \uls^{(2)}, \cdots\}$
in $\Sigma^{\ve}_{\sigma}$ such that
\begin{itemize}
\item[(a)] For every $j=1,\cdots, n$, and $N\geq 1$, $e^{2\pi\iota s_j^{(N)}}=\zeta_j$.
\item[(b)] For $i\neq j$, $s^{(N)}_i - s^{(N)}_j \to \infty$ as $N\to \infty$.
\end{itemize}
Property (a) ensures that we have the following for each $N\geq 1$
\[
\Xi^\ve_\sigma(\ul{\zeta}) = \Phi^{\ve}_{\sigma}(\uls^{(N)})^{-1}\Psi^{\ve}_{\sigma}(\uls^{(N)})
\]
The asymptotics of $\Phi^{\ve}_{\sigma}$ and $\Psi^{\ve}_{\sigma}$, and property (b)
above then imply that, as we let $N\to \infty$, the ratio goes to $1$. Note that,
because of the abelian nature of the difference equations,
the multivalued factors in the asymptotics from (ii) of the statement of the proposition 
cancel out. Thus $\Xi^\ve_\sigma(\ul{\zeta})=1$ for every $\ul{\zeta}\in(\nC)^n$
and we are done.
\end{pf}

\subsection{Kohno--Drinfeld theorem for abelian $q$KZ equations}\label{ss:KD}

Assume now that $V_1,\ldots,V_n$ are non--resonant, and let $\V_i=\Fh
{}(V_i)$ be the corresponding representations of $\qloop$. The following
computes the monodromy of the abelian $q$KZ equations on $V_1\otimes
\cdots\otimes V_n$ in terms of the commutative $R$--matrix of $\qloop$
acting on $\V_1\otimes\cdots\otimes\V_n$.

\begin{thm}
Let $\sigma\in\Sym_n$, and set $\sigma_i=(i\ i+1)$. Then,
\[
(\Phi^{\ve}_{\sigma}(\ul{s}))^{-1}\Phi^{\ve}_{\sigma_i\sigma}(\ul{s}) = 
\Rq^{0,\veps}_{\V_i,\V_{i+1}}(\zeta_i\zeta_{i+1}^{-1})^{\pm 1}
\]
if $(i,i+1)\in C^{\pm}(\sigma)$, where $\zeta_j=\Exp{s_j}$.
\end{thm}
\begin{pf}
This follows from the explicit form of the canonical fundamental
solutions given by Proposition \ref{pr:n var can} and Theorem
\ref{th: n=2KD}.
\end{pf}

\appendix
\section{The inverse of the $q$--Cartan matrix of $\g$}\label{se:KT}

\subsection{}

Let $\bfA = (a_{ij})_{i\in\bfI}$ be a Cartan matrix of finite type,
and $d_i\in\Z_{>0}$ ($i\in\bfI$) be relatively prime symmetrising
integers, \ie $d_ia_{ij} = d_ja_{ji}$ for every $i,j\in\bfI$. Consider
the symmetrised Cartan matrix $\mathbf{B} = (d_ia_{ij})$, and its
$q$--analog $\mathbf{B}(q) = (\qnum{d_ia_{ij}}{})$. The latter
defines a $\C(q)$--valued, symmetric bilinear form on $\bigoplus_{j\in\bfI}
\Q(q) \alpha_j$ by
\[
(\alpha_i,\alpha_j)_q = \qnum{d_ia_{ij}}{}
\]

We give below explicit expressions for the fundamental coweights
$\{\cowt{i}\}_{i\in\bfI}$ in terms of $\{\alpha_i\}$. That is,
we compute certain elements $\cowt{i} \in \bigoplus_{j\in\bfI} \mathbb{Q}(q)\alpha_j$
such that
$(\cowt{i},\alpha_j)_q = \delta_{ij}$ for every $i,j\in\bfI$. The main result
of these calculations is the following. 

\begin{thm}\label{thm: KT-main}
Let 
$l = mh^{\vee}$ where $m = 1, 2, 3$ for types $\mathsf{ADE}, 
\mathsf{BCF}$ and $\mathsf{G}$ respectively, and $h^{\vee}$
is the  dual Coxeter number. Then, for each $i\in\bfI$
\[
\qnum{l}{}\cowt{i} \in \oplus_{j\in\bfI} \N[q,q^{-1}]\alpha_j
\]
\end{thm}

\subsection{}

Below we follow Bourbaki's conventions, especially for the labels of 
the Dynkin diagrams. Recall the standard notations for $q$--numbers
introduced in Section \ref{ssec: qla}:
$\ds \qnum{m}{} = \frac{q^m-q^{-m}}{q-q^{-1}}$. For $m\geq 0$, 
$[m]_q = \sum_{i=0}^{m-1} q^{m-1-2i} \in \N[q,q^{-1}]$.
Moreover, define $\whacky{m} := q^m+q^{-m}$. The following identity is immediate
and will be needed later:

\begin{equation}\label{eq: w-q}
\qnum{a}{}\whacky{b} = \qnum{a+b}{} + \qnum{a-b}{}
\end{equation}
which belongs to $\N[q,q^{-1}]$ if $a\geq b\geq 0$.\\

Also we note that for $a,b\in\N$, with $a\neq 0$, we have
\[
\frac{\qnum{ab}{}}{\qnum{a}{}} = \qnum{b}{a} \in \N[q,q^{-1}]
\]

\subsection{$\mathsf{A}_n$}

In this case $l = n+1$. We have
\[
\cowt{i} = \frac{1}{\qnum{n+1}{}} \lp 
\qnum{n-i+1}{} \lp \sum_{j=1}^{i-1} \qnum{j}{} \alpha_j\rp
+ \qnum{i}{} \lp \sum_{j=i}^n \qnum{n-j+1}{} \alpha_j\rp\rp
\]

Thus the assertion of Theorem \ref{thm: KT-main} holds in this case.

\subsection{$\mathsf{B}_n$}
In this case $l = 2(n+1)$. 
For $1\leq i\leq n-1$ we have
\[
\cowt{i} = \frac{1}{\whacky{n+1}} \lp
\whacky{n-i+1}\lp \sum_{j=1}^{i-1} \qnum{j}{} \alpha_j\rp
+\qnum{i}{}\lp \lp\sum_{j=i}^{n-1} \whacky{n-j+1}\alpha_j\rp + \alpha_n\rp
\rp
\]
and
\[
\cowt{n} =\frac{1}{\whacky{n+1}}\lp
\lp\sum_{j=1}^{n-1} \qnum{j}{}\alpha_j \rp + \frac{\qnum{n}{}}{\qnum{2}{}}
\alpha_n\rp
\]
The statement of Theorem \ref{thm: KT-main} in this case follows
for $1\leq i\leq n-1$ from the identity $\qnum{m}{}\whacky{m} = 
\qnum{2m}{}$. For $\cowt{n}$, we can write (using the same identity)
\[
\cowt{n} =\frac{1}{\qnum{2(n+1)}{}}\lp
\qnum{n+1}{}\lp\sum_{j=1}^{n-1} \qnum{j}{}\alpha_j \rp + 
\frac{\qnum{n+1}{}\qnum{n}{}}{\qnum{2}{}}
\alpha_n\rp
\]
Now it is clear that the coeffient of $\alpha_n$ is a Laurent polynomial
in $q$ with positive integer coefficients.

\subsection{$\mathsf{C}_n$}

In this case $l = 2(2n-1)$. We have the following for each $1\leq i\leq n-1$
\[
\cowt{i} = \frac{1}{\qnum{2}{}\whacky{2n-1}}
\lp
\whacky{2n-2i-1}\lp \sum_{j=1}^{i-1} \qnum{j}{2} \alpha_j \rp
+ \qnum{i}{2}\lp \sum_{j=i}^{n-1} \whacky{2n-2j-1}\alpha_j\rp
+ \qnum{2i}{} \alpha_n\rp
\]
and
\[
\cowt{n} = \frac{1}{\qnum{2}{}\whacky{2n-1}}
\sum_{j=1}^n \qnum{2j}{} \alpha_j
\]
The statement of Theorem \ref{thm: KT-main} follows for $\cowt{n}$. For
$1\leq i\leq n-1$ we will have to use the following variant of
\eqref{eq: w-q}:
\[
\frac{\qnum{2n-1}{}\whacky{2n-2j-1}}{\qnum{2}{}} = 
\frac{\qnum{4n-2j-2}{} + \qnum{2j}{}}{\qnum{2}{}} \in \N[q,q^{-1}]
\]

\subsection{$\mathsf{D}_n$}

In this case $l = 2n-2$. We have the following for $1\leq i\leq n-2$:
\[
\cowt{i} = \frac{1}{\whacky{n-1}} \lp
\whacky{n-i-1} \lp \sum_{j=1}^{i-1} \qnum{j}{} \alpha_j\rp + 
\qnum{i}{}\lp\lp \sum_{j=i}^{n-2} \whacky{n-j-1}\alpha_j\rp
+\alpha_{n-1}+\alpha_n\rp\rp
\]
and
\begin{align*}
\cowt{n-1} &= \frac{1}{\whacky{n-1}}\lp \lp\sum_{j=1}^{n-2} \qnum{j}{}
\alpha_j \rp + \frac{\qnum{n}{}}{\qnum{2}{}}\alpha_{n-1} + 
\frac{\qnum{n-2}{}}{\qnum{2}{}}\alpha_n\rp \\
\cowt{n} &= \frac{1}{\whacky{n-1}}\lp \lp\sum_{j=1}^{n-2} \qnum{j}{}
\alpha_j \rp + \frac{\qnum{n-2}{}}{\qnum{2}{}}\alpha_{n-1} + 
\frac{\qnum{n}{}}{\qnum{2}{}}\alpha_n\rp
\end{align*}

Again we obtain Theorem \ref{thm: KT-main} by the same argument 
as for $\mathsf{B}_n$.

\subsection{$F_4$}

In this case $l = 18$. We get the following
\begin{align*}
\cowt{1} &= \frac{\whacky{3}}{\whacky{9}}\lp \whacky{5}\alpha_1
+ \qnum{3}{2}\alpha_2 + \whacky{2}\alpha_3+
\alpha_4\rp \\
\cowt{2} &= \frac{\whacky{3}}{\whacky{9}}\lp \qnum{3}{2}\alpha_1
+\qnum{6}{}\alpha_2 + \qnum{4}{}\alpha_3 + \qnum{2}{}\alpha_4\rp \\
\cowt3 &= \frac{1}{\whacky{9}}\lp \whacky{2}\whacky{3}\alpha_1
+ \qnum{4}{}\whacky{3}\alpha_2 + \qnum{3}{2}(\whacky{2}\alpha_3
+\alpha_4)\rp \\
\cowt4 &= \frac{1}{\whacky{9}}\lp \whacky{3}\alpha_1 + 
\qnum{2}{}\whacky{3}\alpha_2 + \qnum{3}{2}\alpha_3 + 
\frac{\whacky{3}\whacky{4}}{\qnum{2}{}} \alpha_4\rp
\end{align*}

Again the statement of Theorem \ref{thm: KT-main} is clearly true, except
for the coefficient of $\alpha_4$ in $\cowt{4}$. For that entry we have
\[
\frac{\qnum{9}{}\whacky{3}}{\qnum{2}{}} = \frac{\qnum{12}{} + \qnum{6}{}}{\qnum{2}{}}
\in\N[q,q^{-1}]
\]

\subsection{$G_2$}

In this case $l=12$. We have the following answer
\begin{align*}
\cowt1 &= \frac{\whacky{2}}{\whacky{6}} \lp \frac{\qnum{2}{}}
{\qnum{3}{}} \alpha_1 + \alpha_2\rp & 
\cowt2 &= \frac{\whacky{2}}{\whacky{6}} (\alpha_1 + \whacky{3}\alpha_2)
\end{align*}

As before we multiply and divide these expressions by $\qnum{6}{}$ to
get the denominator $\qnum{12}{}$. Then it is easy to see the coefficients
of $\alpha_1,\alpha_2$ are in $\N[q,q^{-1}]$ as claimed.

\subsection{$\sfE$ series}

The computations below were carried out using \texttt{sage}.

\subsection{$\mathsf{E}_6$}

In this case $l=12$. We have the following expressions:

\begin{align*}
\qnum{12}{}\cowt{1} &= \whacky{3}\qnum{8}{}\alpha_1 
+\whacky{2}\qnum{6}{}\alpha_2
+\whacky{2}\whacky{3}\qnum{5}{}\alpha_3
+\qnum{4}{}\qnum{6}{}\alpha_4\\
&\phantom{=}+\qnum{2}{}\whacky{3}\qnum{4}{}\alpha_5
+\whacky{3}\qnum{4}{}\alpha_6\\
\qnum{12}{}\cowt{2} &= \whacky{2}\qnum{6}{}\alpha_1 
+\whacky{2}\whacky{3}\qnum{6}{}\alpha_2
+\qnum{4}{}\qnum{6}{}\alpha_3
+\whacky{2}\qnum{3}{}\qnum{6}{}\alpha_4\\
&\phantom{=}+\qnum{4}{}\qnum{6}{}\alpha_5
+\whacky{2}\qnum{6}{}\alpha_6\\
\qnum{12}{}\cowt{3} &= \whacky{2}\whacky{3}\qnum{5}{}\alpha_1 
+\qnum{4}{}\qnum{6}{}\alpha_2
+\whacky{3}\qnum{4}{}\qnum{5}{}\alpha_3
+\whacky{1}\qnum{4}{}\qnum{6}{}\alpha_4\\
&\phantom{=}+\qnum{2}{}^2\whacky{3}\qnum{4}{}\alpha_5
+\qnum{2}{}\whacky{3}\qnum{4}{}\alpha_6\\
\qnum{12}{}\cowt{4} &= \qnum{4}{}\qnum{6}{}\alpha_1 
+\whacky{2}\qnum{3}{}\qnum{6}{}\alpha_2
+\qnum{2}{}\qnum{4}{}\qnum{6}{}\alpha_3
+\qnum{3}{}\qnum{4}{}\qnum{6}{}\alpha_4\\
&\phantom{=}+\qnum{2}{}\qnum{4}{}\qnum{6}{}\alpha_5
+\qnum{4}{}\qnum{6}{}\alpha_6\\
\qnum{12}{}\cowt{5} &= \qnum{2}{}\whacky{3}\qnum{4}{}\alpha_1 
+\qnum{4}{}\qnum{6}{}\alpha_2
+\qnum{2}{}^2\whacky{3}\qnum{4}{}\alpha_3
+\qnum{2}{}\qnum{4}{}\qnum{6}{}\alpha_4\\
&\phantom{=}+\whacky{3}\qnum{4}{}\qnum{5}{}\alpha_5
+\whacky{2}\whacky{3}\qnum{5}{}\alpha_6\\
\qnum{12}{}\cowt{6} &= \whacky{3}\qnum{4}{}\alpha_1
+\whacky{2}\qnum{6}{}\alpha_2
+\qnum{2}{}\whacky{3}\qnum{4}{}\alpha_3
+\qnum{4}{}\qnum{6}{}\alpha_4 \\
&\phantom{=} + \whacky{2}\whacky{3}\qnum{5}{}\alpha_5
+\whacky{3}\qnum{8}{}\alpha_6
\end{align*}

\subsection{$\mathsf{E}_7$}

In this case $l = 18$ and we have the following expressions:

\vspace*{-.1in}

\begin{align*}
\whacky{9}\cowt{1} &= \whacky{3}\whacky{5}  \alpha_1
+ \whacky{2}\whacky{3} \alpha_2
+ \whacky{3}\qnum{3}{2} \alpha_3
+ \whacky{3}\qnum{4}{} \alpha_4 \\
&\phantom{=} + \qnum{6}{} \alpha_5
+ \qnum{2}{}\whacky{3}  \alpha_6
+ \whacky{3} \alpha_7 \\
\whacky{9}\cowt{2} &= \whacky{2}\whacky{3}  \alpha_1
+ \frac{\whacky{3}\qnum{7}{}}{\qnum{2}{}}  \alpha_2
+ \whacky{3}\qnum{4}{} \alpha_3
+ \whacky{2}\qnum{6}{} \alpha_4 \\
&\phantom{=} + \qnum{3}{}\qnum{3}{2} \alpha_5
+ \qnum{6}{} \alpha_6
+ \qnum{3}{2} \alpha_7 \\
\whacky{9}\cowt{3} &= \whacky{3}\qnum{3}{2} \alpha_1
+ \whacky{3}\qnum{4}{} \alpha_2
+ \whacky{3}\qnum{6}{} \alpha_3
+ \qnum{2}{}\whacky{3}\qnum{4}{} \alpha_4 \\
&\phantom{=} + \qnum{2}{}\qnum{6}{} \alpha_5
+ \qnum{2}{}^2\whacky{3} \alpha_6
+ \qnum{2}{}\whacky{3} \alpha_7 \\ 
\whacky{9}\cowt{4} &= \whacky{3}\qnum{4}{} \alpha_1
+ \whacky{2}\qnum{6}{} \alpha_2
+ \qnum{2}{}\whacky{3}\qnum{4}{} \alpha_3
+ \qnum{4}{}\qnum{6}{} \alpha_4 \\
&\phantom{=} + \qnum{3}{}\qnum{6}{} \alpha_5
+ \qnum{2}{}\qnum{6}{} \alpha_6
+ \qnum{6}{} \alpha_7 \\
\whacky{9}\cowt{5} &= \qnum{6}{}  \alpha_1
+ \qnum{3}{}\qnum{3}{2} \alpha_2
+ \qnum{2}{}\qnum{6}{} \alpha_3
+ \qnum{3}{}\qnum{6}{} \alpha_4 \\
&\phantom{=} + \qnum{3}{2}\qnum{5}{} \alpha_5
+ \whacky{3}\qnum{5}{} \alpha_6
+ \frac{\whacky{3}\qnum{5}{}}{\qnum{2}{}} \alpha_7 \\
\whacky{9}\cowt{6} &= \qnum{2}{}\whacky{3} \alpha_1
+ \qnum{6}{} \alpha_2
+ \qnum{2}{}^2\whacky{3} \alpha_3
+ \qnum{2}{}\qnum{6}{} \alpha_4 \\
&\phantom{=} + \whacky{3}\qnum{5}{} \alpha_5
+ \qnum{2}{}\whacky{3}\whacky{4} \alpha_6
+ \whacky{3}\whacky{4} \alpha_7 \\
\whacky{9}\cowt{7} &= \whacky{3} \alpha_1
+ \qnum{3}{2} \alpha_2
+ \qnum{2}{}\whacky{3} \alpha_3
+ \qnum{6}{} \alpha_4 \\
&\phantom{=} + \frac{\whacky{3}\qnum{5}{}}{\qnum{2}{}} \alpha_5
+ \whacky{3}\whacky{4} \alpha_6
+ \qnum{3}{4} \alpha_7 \\
\end{align*}

It only remains to observe that 
\[
\frac{\qnum{9}{}\whacky{3}}{\qnum{2}{}} = \frac{\qnum{12}{} + \qnum{6}{}}{\qnum{2}{}}
 = \qnum{6}{2} + \qnum{3}{2} \in \N[q,q^{-1}]
\]
\vspace*{-.20in}

\subsection{$\mathsf{E}_8$}

In this case $l = 30$ and we have the following expression:

\vspace*{-.1in}

\begin{align*}
\whacky{15}\cowt{1} &= \whacky{5}\qnum{4}{3} \alpha_1
+\whacky{3}\qnum{5}{2} \alpha_2
+\qnum{2}{3}\frac{\whacky{5}\qnum{7}{}}{\qnum{2}{}} \alpha_3
+\whacky{3}\qnum{10}{} \alpha_4\\
&\phantom{=} + \whacky{3}\qnum{4}{}\whacky{5}\alpha_5
+ \whacky{5}\qnum{6}{}\alpha_6
+ \qnum{2}{}\whacky{3}\whacky{5}\alpha_7
+ \whacky{3}\whacky{5}\alpha_8\\
\whacky{15}\cowt{2} &= \whacky{3}\qnum{5}{2} \alpha_1
+ \whacky{3}\whacky{5}\qnum{4}{2}\alpha_2
+ \whacky{3}\qnum{10}{}\alpha_3
+ \qnum{3}{2}\qnum{10}{}\alpha_4\\
&\phantom{=} + \whacky{2}\whacky{5}\qnum{6}{} \alpha_5
+ \qnum{3}{}\qnum{3}{2}\whacky{5} \alpha_6
+ \whacky{5}\qnum{6}{} \alpha_7
+ \whacky{5}\qnum{3}{2}\alpha_8\\
\whacky{15}\cowt{3} &= \qnum{2}{3}\frac{\whacky{5}\qnum{7}{}}{\qnum{2}{}} \alpha_1
+ \whacky{3}\qnum{10}{} \alpha_2
+ \qnum{2}{3}\whacky{5}\qnum{7}{} \alpha_3
+ \qnum{2}{}\whacky{3}\qnum{10}{} \alpha_4\\
&\phantom{=} + \qnum{2}{}\whacky{3}\qnum{4}{}\whacky{5} \alpha_5
+ \qnum{2}{}\whacky{5}\qnum{6}{}\alpha_6
+ \qnum{2}{}^2\whacky{3}\whacky{5}\alpha_7
+ \qnum{2}{}\whacky{3}\whacky{5} \alpha_8\\
\whacky{15}\cowt{4} &= \whacky{3}\qnum{10}{} \alpha_1
+ \qnum{3}{2}\qnum{10}{} \alpha_2
+ \qnum{2}{}\whacky{3}\qnum{10}{} \alpha_3
+ \qnum{6}{}\qnum{10}{} \alpha_4\\
&\phantom{=} + \qnum{4}{}\whacky{5}\qnum{6}{}\alpha_5
+ \qnum{3}{}\whacky{5}\qnum{6}{}\alpha_6
+ \qnum{2}{}\whacky{5}\qnum{6}{} \alpha_7
+ \whacky{5}\qnum{6}{}\alpha_8\\
\whacky{15}\cowt{5} &= \whacky{3}\qnum{4}{}\whacky{5} \alpha_1
+ \whacky{2}\whacky{5}\qnum{6}{} \alpha_2
+ \qnum{2}{}\whacky{3}\qnum{4}{}\whacky{5} \alpha_3
+ \qnum{4}{}\whacky{5}\qnum{6}{} \alpha_4\\
&\phantom{=} + \whacky{2}\whacky{3}\qnum{10}{} \alpha_5
+ \qnum{3}{2}\qnum{10}{} \alpha_6
+ \whacky{3}\qnum{10}{} \alpha_7
+ \whacky{3}\qnum{5}{2} \alpha_8\\
\whacky{15}\cowt{6} &= \whacky{5}\qnum{6}{} \alpha_1
+ \qnum{3}{}\qnum{3}{2}\whacky{5} \alpha_2
+ \qnum{2}{}\whacky{5}\qnum{6}{} \alpha_3
+ \qnum{3}{}\whacky{5}\qnum{6}{} \alpha_4\\
&\phantom{=} + \qnum{3}{2}\qnum{10}{} \alpha_5
+ \whacky{4}\whacky{5}\qnum{6}{} \alpha_6
+ \qnum{2}{}\qnum{2}{3}\whacky{4}\whacky{5}\alpha_7
+ \qnum{2}{3}\whacky{4}\whacky{5} \alpha_8\\
\whacky{15}\cowt{7} &= \qnum{2}{}\whacky{3}\whacky{5} \alpha_1
+ \whacky{5}\qnum{6}{} \alpha_2
+ \qnum{2}{}^2\whacky{3}\whacky{5} \alpha_3
+ \qnum{2}{}\whacky{5}\qnum{6}{} \alpha_4\\
&\phantom{=} + \whacky{3}\qnum{10}{} \alpha_5
+ \qnum{2}{}\qnum{2}{3}\whacky{4}\whacky{5} \alpha_6
+ \qnum{2}{}\qnum{3}{4}\whacky{5} \alpha_7
+ \qnum{3}{4}\whacky{5} \alpha_8\\
\whacky{15}\cowt{8} &= \whacky{3}\whacky{5} \alpha_1
+ \whacky{5}\qnum{3}{2} \alpha_2
+ \qnum{2}{}\whacky{3}\whacky{5} \alpha_3
+ \whacky{5}\qnum{6}{} \alpha_4\\
&\phantom{=} + \whacky{3}\qnum{5}{2} \alpha_5
+ \qnum{2}{3}\whacky{4}\whacky{5} \alpha_6
+ \qnum{3}{4}\whacky{5} \alpha_7
+ \whacky{5}\whacky{9} \alpha_8
\end{align*}
It only remains to observe that 
\[
\frac{\qnum{15}{}\whacky{5}}{\qnum{2}{}} = \frac{\qnum{20}{} + \qnum{10}{}}{\qnum{2}{}}
 = \qnum{10}{2} + \qnum{5}{2} \in \N[q,q^{-1}]
\]

\bibliographystyle{amsplain}
\bibliography{slayer}

\end{document}